\documentclass[11pt,equ]{article}

\usepackage{graphicx}
\usepackage{hyperref}
\newcommand{\ord}{\frak o}
\newcommand{\pp}{\frak p}
\newcommand{\kk}{k}
\newcommand{\n}{n}
\newcommand{\cM}{{\cal M}}

\newcommand{\cD}{{\cal D}}
\newcommand{\cX}{{\cal X}}
\newcommand{\cY}{{\cal Y}}
\newcommand{\gb}{\beta}
\newcommand{\gd}{\delta}
\newcommand{\gl}{\lambda}
\newcommand{\gw}{\omega}
\newcommand{\g}{\gamma}
\newcommand{\gt}{\theta}
\newcommand{\s}{\frak s}
\newcommand{\gS}{\Sigma}
\newcommand{\gs}{\sigma}
\newcommand{\vep}{\varepsilon}
\newcommand{\EE}{{\Bbb E}}
\newcommand{\RW}{RRW}
\newtheorem{theorem}{Theorem}[section]
\newtheorem{corollary}{Corollary}[section]

\newtheorem{defi}{Definition}[section]
\newtheorem{conj}{Conjecture}[section]

\newcommand{\eqn}[1]{(\ref{#1})}

\newcommand{\df}{\displaystyle\frac}

\newcommand{\hsp}{\hspace{\parindent}}

\newcommand{\bsq}{\vrule height .9ex width .8ex depth -.1ex}

\newcommand{\DD}{{\Bbb D}}

\newcommand{\NN}{{\Bbb N}}

\newcommand{\PP}{{\Bbb P}}
\newcommand{\QQ}{{\Bbb Q}}
\newcommand{\RR}{{\Bbb R}}

\newcommand{\ZZ}{{\Bbb Z}}

\newcommand{\sR}{{\cal R}}

\newcommand{\sE}{{\cal E}}

\newcommand{\sT}{{\cal T}}
\newcommand{\sB}{{\cal B}}

\newcommand{\beql}[1]{\begin{equation}\label{#1}}
\newcommand{\eeq}{\end{equation}}

\catcode`\@=11
\renewcommand{\section}{
        \setcounter{equation}{0}
        \@startsection {section}{1}{\z@}{-3.5ex plus -1ex minus
        -.2ex}{2.3ex plus .2ex}{\large\bf}
        }
\catcode`@=12
\makeatletter
\def\eqalignno#1{\displ@y \ta {\bf s} kip\@centering
  \halign to\displaywidth{\hfil$\@lign\displaystyle{##}$\ta {\bf s} kip\z@skip
    & $\@lign\displaystyle{{}##}$\hfil\ta {\bf s} kip\@centering
    & \llap{$\@lign##$}\ta {\bf s} kip\z@skip\crcr
    #1\crcr}}
\makeatother

\makeatletter
\def\@sect#1#2#3#4#5#6[#7]#8{\ifnum #2>\c@secnumdepth
     \def\@svsec{}\else 
     \refstepcounter{#1}\edef\@svsec{\csname the#1\endcsname.\hskip .75em }\fi
     \@tempskipa #5\relax
      \ifdim \@tempskipa>\z@ 
        \begingroup #6\relax
          \@hangfrom{\hskip #3\relax\@svsec}{\interlinepenalty \@M #8\par}%
        \endgroup
       \csname #1mark\endcsname{#7}\addcontentsline
         {toc}{#1}{\ifnum #2>\c@secnumdepth \else
                      \protect\numberline{\csname the#1\endcsname}\fi
                    #7}\else
        \def\@svsechd{#6\hskip #3\@svsec #8\csname #1mark\endcsname
                      {#7}\addcontentsline
                           {toc}{#1}{\ifnum #2>\c@secnumdepth \else
                             \protect\numberline{\csname the#1\endcsname}\fi
                       #7}}\fi
     \@xsect{#5}}
\def\@theorem#1#2{\it \trivlist \item[\hskip \labelsep{\bf #1\ #2.}]}
\input amssym.def
\input amssym.tex

\makeatother
\begin{document}
\begin{center}
{\large {\bf Stochastic Models for the  $3x+1$ and $5x+1$ Problems
}}\\ 

\bigskip

{\large {\em Alex V. Kontorovich}
\footnote{AVK received support from an NSF Postdoc, grant DMS 0802998.}} \\
Department of Mathematics\\
Brown University\\
Providence, RI \\
{\tt alexk@math.brown.edu}\\

\ 

\

{\large {\em Jeffrey C. Lagarias}
\footnote{JCL received support from  NSF Grants DMS-0500555 and DMS-0801029.}} \\
Department of Mathematics\\
University of Michigan \\
Ann Arbor, MI 48109-1109 \\
{\tt lagarias@umich.edu} \bigskip \\


\setcounter{tocdepth}{1}

\vspace*{2\baselineskip}
(October 7, 2009%
) \\

\vspace*{1\baselineskip}

{\em Abstract}
\end{center}

This paper discusses stochastic models for predicting the long-time behavior of 
the trajectories of orbits of the $3x+1$ problem and, for comparison, the $5x+1$
problem.
The stochastic models are rigorously analyzable, and yield heuristic predictions
(conjectures) for the behavior of $3x+1$ orbits and $5x+1$ orbits. 
 \vspace*{1\baselineskip}



\noindent

\setlength{\baselineskip}{1.0\baselineskip}
%
%
%
\section{Introduction}

The $3x+1$ problem concerns
the following operation on integers: if an integer is odd
``multiply by three and add one," while if it is even  ``divide by two."
This operation is given by the  {\em Collatz function}
\beql{101}
C(\n) = 
\left\{
\begin{array}{cl}
3n+1 & \mbox{if}~ \n \equiv 1~~ (\bmod ~2 ) ~,  \\\
~~~ \\
\df{\n}{2} & \mbox{if} ~~\n \equiv 0~~ (\bmod~2) ~.
\end{array}
\right.
\eeq
The $3x+1$ problem concerns  what happens if one iterates this operation
starting from a given positive integer $\n$.
The  unsolved {\em $3x+1$ Problem}  or {\em Collatz problem} is to prove
(or disprove) that such 
iterations always eventually  reach the number $1$
(and therefter cycle, taking values $1, 4, 2,1$). 
This problem  goes under many other names, including:
 {\em Syracuse Problem}, {\em Hasse's Algorithm},
 {\em Kakutani's Problem}
and {\em Ulam's Problem}.  \\

The $3x+1$ Conjecture has now been verified for all $\n \le 
5.67  \times 10^{18}
$ by
computer experiments \cite{OeS09}. \\

%
%
%

 \subsection{$3x+1$ Function}

There are a number of different functions that encode the $3x+1$ problem, which
proceed through the iteration at different speeds. The following two functions
prove to be  more convenient
for probabilistic analysis than the Collatz function. 
The first  of these is the {\em $3x+1$ function} $T(\n)$ (or {\em $3x+1$ map})
\beql{102}
T(\n) = \left\{
\begin{array}{cl}
\df{3\n+1}{2}  & \mbox{if} ~~\n \equiv 1~~ (\bmod ~ 2)~, \\
~~~ \\
\df{\n}{2} & \mbox{if}~~ \n \equiv 0~~ (\bmod ~2 )~.
\end{array}
\right.
\eeq
This function divides out one power of $2$, after an odd input is encountered;
it is defined on the domain  of all integers. \\

The second
function, the {\em accelerated $3x+1$  function} $U(\n)$, is defined on the domain of 
all odd integers, and removes all powers of $2$ at each step. It is given by 
\beql{103}
U(\n) = \frac{3\n+1}{2^{{\rm ord}_2(3\n+1)}},
\eeq
in which ${\rm ord}_2(\n)$ counts the number of powers of $2$ dividing $\n$.
The function $U(\n)$ was studied by Crandall \cite{Cra78} in 1978.\\

The long-term dynamics under iteration of the $3x+1$ map has proved
resistant to rigorous analysis. It is conjectured that there 
is a finite positive  constant $C$ so that all trajectories eventually
enter and stay in the region $-C \le n \le C.$ In particular, there 
are finitely many
periodic orbits and 
all trajectories eventually enter one of these periodic orbits. 
On the  domain of positive integers it is conjectured there is 
is a single periodic orbit $\{ 1, 2\}$; this is part of the $3x+1$ Conjecture.
 On the domain of negative integers, the known periodic orbits are the three orbits
$\{ -1\} $,   $\{-5, -7, -10\}$ and
 $\{ -17, -25, -37, -55, -82, -41, -61, -91, -136, -68, -34\}$.\\

%
%
%

 \subsection{$5x+1$ Problem}

For comparison purposes, we also consider the {\em $5x+1$ problem}, which 
concerns iterates of 
the {\em Collatz $5x+1$ function}
\beql{111}
C_5(\n) = 
\left\{
\begin{array}{cl}
5\n+1 & \mbox{if}~ \n \equiv 1~~ (\bmod ~2 ) ~,  \\\
~~~ \\
\df{\n}{2} & \mbox{if} ~~\n \equiv 0~~ (\bmod~2) ~.
\end{array}
\right.
\eeq
For this function we also have analogues of the
other two functions above. We define the {\em $5x+1$ function} $T_5(\n)$ (or {\em $5x+1$ map}), given by 
\beql{102a}
T_5(\n) = \left\{
\begin{array}{cl}
\df{5\n+1}{2}  & \mbox{if} ~~\n \equiv 1~~ (\bmod ~ 2)~, \\
~~~ \\
\df{\n}{2} & \mbox{if}~~ \n \equiv 0~~ (\bmod ~2 )~.
\end{array}
\right.
\eeq
It is defined on the set of all integers.\\
 
The second
function, the {\em accelerated $5x+1$  function} $U_5(\n)$, is defined on the domain of 
all odd integers, and removes all powers of $2$ at each step. It is given by 
\beql{103a}
U_5(\n) = \frac{5\n+1}{2^{ord_2(5\n+1)}},
\eeq
in which $ord_2(\n)$ counts the number of powers of $2$ dividing $\n$.\\



The long-term dynamics under iteration of the $5x+1$ map on the integers is conjecturally
quite different from the $3x+1$ map. It is conjectured that  a  density one set of integers  belong to
divergent trajectories, ones with  $|T^{(k)}(n)| \to \infty$.
It is also conjectured that there are a finite number of periodic orbits, which include
the orbits $\{1, 3, 8, 4, 2\}$ and $\{ 13, 33, 83, 208, 104, 52, 26\}$ on the positive integers
and the orbit  $\{-1, -2\}$ on the negative integers. An infinite number of trajectories eventually enter
one of these orbits, but the set of all integers entering each of these orbits is believed to have
density zero. \\
 \\

%
%
%

 \subsection{Stochastic models }

This paper is concerned with probabilistic models for the behavior
of the $3x+1$ function iterates, and for comparison, the $5x+1$ function iterates.
The absence of rigorous analysis  of the long-term behavior under
iteration of these functions provides 
one motivation  to formulate probabilistic models of the behavior of the
$3x+1$ map and $5x+1$ map. These models can  make
predictions that can be 
compared to empirical data, which,  by uncovering discrepancies,   may 
lead to the discovery of new hidden
regularities in their behavior under iterations. 
Note that both the $3x+1$ map and the $5x+1$ map have the positive integers and 
negative integers as invariant subsets; thus their dynamics can be studied 
separately on these domains. The original problems concern 
their dynamics 
restricted to the positive integers.\\

Here we survey what is known about iteration of these maps, 
in  frameworks  which have a probabilistic interpretation.
A great deal is known about the initial behavior of the iteration of the
$3x+1$ map and $5x+1$ map; such results are summarized in \S
\ref{sec2} and \S\ref{sec7},
respectively. Here some models for the $5x+1$ problem are new, developed
in parallel with models in Lagarias and Weiss \cite{LW92}. 
The major unsolved questions have to do with the  behavior
of long term aspects of the iterations. It  is here that stochastic models
have an important role to play. 
We present  models for forward iteration of the map 
which are of
random walk or Markov process type, and models for backwards iteration of the map, which
are  branching processes or branching random walks.
Such models can address   how the iteration behaves for a randomly selected input 
value $n$.  More sophisticated models address behavior of ``extremal" input values. 
Analysis of these latter  models typically uses some variant of the theory of large deviations. \\

We are interested in using these stochastic models to
explore similarities and differences between the iteration behavior
of the $3x+1$ and $5x+1$  functions. There are many similarities which
are exact parallels, listed
in the concluding 
\S\ref{sec11}. The main differences are: in short term iteration on the
integers $\ZZ$,
$3x+1$ iterates tend to get smaller, while $5x+1$ iterates tend to get larger
(in absolute value). For long term iteration it is conjectured that all $3x+1$ trajectories
eventually enter finite cycles; it is conjectured that almost all $5x+1$ trajectories
diverge. Stochastic models permit making some quantitative versions of this
behavior. These include the following  (conjectural) predictions.

\begin{enumerate}
\item
The number of integers $1 \le n \le x$ whose $3x+1$  forward orbit reaches
$1$ is about $x^{\eta_3 + o(1)}$, where $\eta_3=1.$
\item
Restricting to those integers $1 \le n \le x$ whose $3x+1$ map forward orbit includes $1$, 
the trajectories of most such $n$  reach $1$ after about $6.95212 \log n$ steps.

\item
Only finitely many  $3x+1$ map trajectories  starting 
at $x$ reach $1$
after more than $(\gamma_{3}+ \epsilon) \log x$ steps, while infinitely
many positive $x$ reach $1$ after more than  $(\gamma_{3}- 
\epsilon) \log x$
steps, where $\gamma_3 \approx 41.67765$.

\item
The number of integers $1 \le n \le x$ whose $5x+1$ map forward orbit includes
$1$ is about $x^{\eta_5 + o(1)}$, where $\eta_5 \approx 0.65049.$

\item
Restricting to those integers $1\le n \le x$ whose $5x+1$ map forward orbit includes
$1$, the trajectories of most such $n$  reach $1$ after about $9.19963 \log n$ steps. 
\item
Only finitely many $5x+1$ map trajectories starting at $x$ reach $1$
after more than $(\gamma_{5}+ \epsilon) \log x$ steps, while infinitely
many positive $x$ reach $1$ after more than  $(\gamma_{5}- 
\epsilon) \log x$
steps, where $\gamma_5 \approx 84.76012$.
\end{enumerate}

In the case of the $3x+1$ map, extensive numerical evidence supports these predictions.
There has been much less computational testing of the $5x+1$ map, and the 
predictions above are less tested in these cases. 

We also survey a number of rigorous results that fit in this framework:
these results describe aspects of the initial part of the iteration.
These include symbolic dynamics for accelerated iteration, given in \S\ref{sec6},  which were
used by Kontorovich and Sinai \cite{KS02} to show that suitably scaled versions of initial trajectories
converge in a limit to geometric Brownian motion. 
These also include results on Benford's law for the initial base B digits of the
initial iterates 
of the functions above,
given in \S\ref{sec9}. \\

%
%
%

 \subsection{Contents of the paper }

In 
\S\ref{sec2} through \S\ref{sec6} we first consider the $3x+1$ function. Then in
\S\ref{sec7} and \S\ref{sec8} we give comparison results for the $5x+1$ problem.
In 
\S\ref{sec9} and \S\ref{sec10} we give results
on 
Benford's law and for
$2$-adic generalizations, in parallel for 
both the $3x+1$ function and
$5x+1$ function.\\

In 
\S\ref{sec2}  we discuss the iteration of
the $3x+1$ map. We describe its symbolic dynamics, and formulate
several statistics of orbits, which will be 
studied via
stochastic models in later sections. We  state various rigorously proved results
about these statistics. 
 For a given starting value $n$, 
 these  statistics include  the {\em $\lambda$-stopping time} $\sigma_{\lambda}(n)$,
the {\em total stopping time} $\sigma_{\infty}(n),$ the {\em maximum excursion value} $t(n)$, and 
{\em counting functions }$N_k(n)$ and $N_k^{\ast} (n)$,
for the number of backward iterates at depth $k$ of a given
integer $a$, with the latter only counting  iterates that are not divisible by $3$.
We also review what has been rigorously proved about these statistics, 
and give tables of empirical results known about these statistics, found by large scale
computations. Further
 data appears in the paper of Oliveira e Silva \cite{OeS09}
(in this volume). \\

In 
\S\ref{sec3} we discuss stochastic models for a single
orbit under forward iteration of the $3x+1$ map. 
These include a multiplicative random product model (MRP model) and a
logarithmic rescaling giving an additive random walk model
taking unequal steps (BRW model), which has a negative drift. These 
models predict that all orbits converge to a bounded set,
and that the total stopping time $\sigma_{\infty}(n)$ for the $3x+1$ map of a random starting point $n$
should be about $6.95212 \log n$ steps, and as $n \to \infty$ have
 a Gaussian distribution around this value, with standard deviation 
proportional to $\sqrt{\log n}$.\\

In 
\S\ref{sec4} we discuss  models for extreme values of the 
total stopping time of the $3x+1$ map.
We introduce  a repeated random walk model (RRW model)
which produces a random trajectory separately for each integer $n$.
We present  results obtained using  the theory of large deviations which rigorously
determine behavior in this model of
a statistic  which is an analogue of 
 the scaled total stopping time $\gamma(n):= \frac{\sigma_{\infty}(n)}{\log n}$.
The model predicts that the   limit superior of these values
should be a constant $\gamma_{\RW} \approx 41.67765$, which is larger
than the average value $6.95212$ this variable takes.
This prediction agrees fairly well  with the empirical data given in 
\S\ref{sec2}.\\

In 
\S\ref{sec5} we survey results concerning 
forward iteration of the accelerated $3x+1$ map. These include
a complete description of its symbolic dynamics. We also show that
a  suitable scaling limit of these trajectories is a 
geometric Brownian motion, and discuss the equidistribution of various images via entropy.\\

In  
\S\ref{sec6} we describe stochastic models 
simulating  backward
iteration of the $3x+1$ function. These models grow random labelled trees,
whose levels describe branching random walks.
These models give exact answers for the expected number of
leaves at a given depth  $k$, analogous to the number of integers
having total stopping time $k$, and also predict the extremal behavior of
the scaled total stopping time function $\gamma(n):= \frac{\sigma_{\infty}(n)}{\log n}$.
It yields a prediction for the limit superior of these values to be 
$\gamma_{BP} \approx 41.677647$, the same value as for the repeated random
walk process above.\\

In \S\ref{sec7} and \S\ref{sec8} we present analogous results for the $5x+1$ map. 
Much less empirical study has been made for iteration
of the $5x+1$ function, so there is less empirical data available
for comparison.\\

In  
\S\ref{sec7} we define $5x+1$ statistics of orbits.
These are analogues of the $3x+1$ statistics given in \S\ref{sec3}, but some require
modification to reflect the fact that  $5x+1$ orbits grow on average. We also 
review what is known rigorously about the behavior of this function; in particular 
the symbolic dynamics of the forward iteration of the $5x+1$ map is exactly the same as
that for the $3x+1$ map.
The statistics introduced include a reverse
analogue of the stopping time, {\em  the $\lambda^{+}$-stopping time $\sigma_{\lambda}^{+}(n)$,
and also the {\em total stopping time $\sigma_{\infty}(n; T_5)$,
Since most trajectories are believed to be unbounded, the total stopping time
is believed to 
take the value $+\infty$ for almost all initial conditons.
In place of the maximum excursion value, we consider
the {\em minimum excursion value $t^{-}(n)$}! \\

In  
\S\ref{sec8} we present results on stochastic models for the $5x+1$ iteration.
These include repeated random walk models for the forward iteration of
this function, paralleling results of \S\ref{sec4}; 
the convergence to Brownian motion of appropriately scaled trajectories, paralleling
results of \S\ref{sec5}; 
and  branching random walk models for inverse iteration,
paralleling results of \S\ref{sec6}. In the latter case we present some new results.
The most interesting results of the analysis of these models is the prediction
that the number of integers below $x$ which iterate under the $5x+1$ to $1$
should be about $x^{\delta_5 +o(1)}$ with $\delta_5 \approx 0.65041$, and
that all integers below $x$ that eventually iterate to $1$ necessarily do it
in at most $(\gamma_{5, BP} +o(1)) \log x$ steps, where 
$\gamma_{5, BP} \approx 84.76012$. \\

In 
\S\ref{sec9} we discuss another property of $3x+1$ iterates and $5x+1$ iterates:
Benford's law.  In this context ``Benford's law" asserts that 
the 
distribution of the initial decimal digits of
numbers in a 
trajectory  $\{ T^{(k)}(n): 1 \le k \le m \}$ approaches  a 
particular non-uniform probability distribution,
the Benford distribution, in which an  initial digit less than $k$  occurs with probability
 $\log_{10} k$, so that $1$ is the most likely initial digit.
We summarize results showing that 
most initial starting values  of both the
$3x+1$ map and the 
$5x+1$ map have initial iterates exhibiting 
 Benford-like
behavior; this property holds for any fixed finite set of initial iterates.\\

In 
\S\ref{sec10} we review results on  the extensions to the 
domain of  $2$-adic integers $\ZZ_2$
of the 
 functions $T_3(\n)$ and $T_5(\n)$.
 These functions have  the pleasant property that their definition makes
sense $2$-adically, and each function has a unique 
continuous $2$-adic extension, which we denote 
 $\tilde{T}_3: \ZZ_2 \to \ZZ_2$
and $\tilde{T}_5:\ZZ_2 \to \ZZ_2$, respectively.
These extended maps are measure-preserving for the $2$-adic Haar measure, and are 
 ergodic in  a very strong sense. The interesting feature
is that at the level of $2$-adic extensions the $3x+1$ map and $5x+1$ 
map are identical maps from the perspective of
measure theory.  They are both topologically and measurably conjugate
to the full shift on the $2$-adic integers, hence they are topologically and
measurably conjugate to each other! Thus their dynamics are ``the same."
This contrasts with the great difference between these maps view on
the domain of integers. \\ 

In 
\S\ref{sec11} we present concluding remarks, summarizing this paper, 
 comparing properties under iteration of the $3x+1$ map and $5x+1$ map.
The short-run behavior under iteration of these maps have some strong similarities. 
However all evidence indicates that the long-run behavior of iteration 
for the $3x+1$ map and the $5x+1$ map
on the integers $\ZZ$ is very different. We also list a set of insights and topics
for further investigation.

\paragraph{Notation.}
For convenience, when comparing the $3x+1$ maps with the corresponding
$5x+1$ maps, we may write $C_3(\n), T_3(\n), U_3(\n)$ in place of $C(\n), T(\n), U(\n)$ above.\\

\paragraph{Acknowledgments.} 
The authors thank Steven J. Miller for a careful reading of and many corrections to an earlier draft of 
this manuscript.
 AVK wishes to thank the hospitality of  Dorian Goldfeld and Columbia University  during this project.

%
%

\section{The $3x+1$ Function: Symbolic Dynamics and Orbit Statistics}\label{sec2}

In this section we consider  the $3x+1$ map $T(\n)$.
We recall basic properties of its symbolic dynamics. We also define
several different  statistics for describing  its behavior on individual trajectories, and summarize
what is rigorously proved about these statistics. In later sections
we will present probabilistic models which are intended to model the
behavior of these statistics.

%
%

\subsection{$3x+1$ Symbolic Dynamics: Parity Sequence}

The behavior of the map $T(n)$ under iteration is completely described by
the parities of the successive iterates. 

\begin{defi}~\label{de21}
{\em 
(i) For a function $T: \ZZ \to \ZZ$ and input value $\n \in \ZZ$ define
the  {\em parity sequence} of $n$ to be
\beql{201}
S(\n) := ( \n ~(\bmod~2), T(\n)~(\bmod~2), T^{(2)}(\n) ~(\bmod~2), ...)
\eeq
in which $T^{(k)}(\n)$ denotes the $k$-th iterate, so that $T^{(2)}(\n) := T( T(\n))$.
This is an infinite vector of zeros and ones.\\

(ii) For $k \ge 1$ its {\em $k$-truncated parity sequence} is 
a vector giving the  initial segment of $k$ terms of $S(\n)$, i.e.
\beql{202}
S^{[k]}(\n):= ( \n ~(\bmod~2), T(\n)~(\bmod~2), T^{(2)}(\n) ~(\bmod~2),\cdots ,T^{(k-1)} (\n) ~(\bmod~2)).
\eeq 
}
\end{defi}

A basic result on the iteration is as follows.
%
%

\begin{theorem}~\label{th21}
{\em (Parity Sequence Symbolic Dynamics)}
The $k$-truncated parity sequence $S^{[k]}(n)$ of the first $k$ iterates of the $3x+1$ map $T(x)$
is periodic  in $n$ with period $2^k$. Each of the $2^k$ possible $0-1$ vectors occurs
exactly once in the initial segment $1 \le n \le 2^k$.
\end{theorem}

\paragraph{Proof.} This result is
due to Terras \cite{Te76} in 1976
and Everett \cite{Ev77} in 1977.
A proof is given as Theorem B in Lagarias \cite{Lag85}.
$~~~\bsq$  \\

An immediate consequence is that an integer $n$ is uniquely
determined by the parity sequence $S(n)$ of its forward orbit.
To see this, note that 
any two distinct integers fall in different residue classes
$(\bmod~ 2^k)$ for large enough $k$, so will have different parity sequences.
The parity sequence thus provides a {\em symbolic dynamics} which keeps track of
the orbit. Taken on the integers, only countably many different parity sequences occur
out of the uncountably many possible infinite $0-1$ sequences.

%
%

\subsection{$3x+1$ Stopping Time Statistics: $\lambda$-stopping times}

  The initial statistic we consider is the number of iteration steps needed
  to observe a fixed amount of decrease of size in the iterate.

%
%
%

\begin{defi}~\label{de22}
{\em
For fixed $\lambda >0$ the {\em $\lambda$-stopping time} 
$\sigma_{\lambda}(n)$ of a map $T: \ZZ \to \ZZ$ from input $n$ is  the
minimal value of $k \ge 0$ such that
$T^{(k)}(n) < \lambda n$, e.g.
\beql{231}
\sigma_{\lambda}(n) := \inf 
\left\{
k \ge 0:  \frac{T^{(k)}(n)}{n}  < \lambda
\right\}.
\eeq
If no such value $k$ exists, we set  $\sigma_{\lambda}(n) = + \infty.$
}
\end{defi}

This notion for $\lambda=1$ was introduced in 1976 by Terras \cite{Te76}
who called it the {\em stopping time}, and denoted it   $\sigma(n)$.
The more general $\lambda$-stopping time  is interesting in the range  $0 <\lambda \le 1$;
it satisfies  $\sigma_{\lambda}(n)=0$ for all $\lambda >1$.

 Terras \cite{Te76}  
studied the set of numbers having stopping time at most $k$, denoted
\beql{232}
S_{1}(k) :=\{ n:  ~\sigma_{1}(n) \le k\}.
\eeq 
He used Theorem~\ref{th21} to show  (\cite{Te76}, \cite{Te79})
 that this set of integers
has a natural density, as defined below, and that this
density approaches $1$ as $k \to \infty$.  

 Later this result was generalized. Rawsthorne \cite{Raw85} in 1985  
introduced the case  of general
$\lambda$, and Borovkov and Pfeifer \cite[Theorem 2]{BP00} in 2000 considered
criteria with several  stopping time conditions. \\

There are several notions of density of a set $\Sigma$
of the natural numbers $\NN=\{ 1,2, 3, ...\}$. %
The  {\em lower asymptotic density} $\underline{\DD}(\Sigma)$
 is defined for all infinite sets $\Sigma$, and is  given by 
\beql{220}
\underline{\DD}(\Sigma) := \liminf_{t \to \infty} \frac{1}{t} | \{n \in \Sigma: ~n \le t\}|.
\eeq
The assertion that an  infinite set $\Sigma \subset \NN$ of natural numbers
has a  {\em natural density $\DD(\Sigma)$} is the assertion that the following limit exists: 
\beql{220b}
\DD(\Sigma) := \lim_{t \to \infty} \frac{1}{t} | \{n \in \Sigma: ~n \le t\}|.
\eeq
Sets with a natural density automatically have $\DD(\Sigma)=\underline{\DD}(\Sigma).$

%
%

\begin{theorem}~\label{th22} {\em ($\lambda$-Stopping Time Natural  Density)}

(i) For the $3x+1$ map $T(n)$, and any fixed $0< \lambda \le 1$ and $k\ge1$,  the set $S_{\gl}(k)$ of 
integers having
$\lambda$-stopping time at most $k$ 
has a well-defined natural density 
$\DD( S_{\lambda}(k))$. \\

(ii) 
For $\gl$ fixed and $k\to\infty$, this
natural density satisfies 
\beql{221}
\DD (S_{\lambda}(k)) 
\to
1.
\eeq
In particular, the set of numbers with finite $\lambda$-stopping time has
natural density $1$.
\end{theorem}

\paragraph{Proof.} 
For the
special case $\lambda=1$, that is  
the stopping time,  
this 
is the basic result of   Terras \cite{Te76}, \cite{Te79},
obtained also by Everett \cite{Ev77}.
A proof for $\lambda=1$ is given as 
Theorem A in Lagarias \cite{Lag85}. 
The  idea 
is that 
by
Theorem ~\ref{th21}, 
each arithmetic progression $(\bmod~ 2^k)$
has iterates that multiply by a certain pattern of 
 $\frac{1}{2}$ or $\frac{3}{2}$
for the first $k$ steps. A certain 
subset
of the $2^k$-arithmetic
progressions $(\bmod 2^k)$ will have the product of these
numbers fall below $\lambda$, and these
arithmetic progressions give the density. To see that the density goes to $1$ 
as $k \to \infty$, one  must show that most arithmetic progressions $(\bmod 2^k)$ have
a product smaller than one. Theorem~\ref{th21}
says that all products occur equally likely, and since the 
 geometric mean of these products is 
 $\left(\frac{3}{4}\right)^{\frac{1}{2}} < 1$, 
one can establish 
that
such a decrease occurs for all but an exponentially small
set of patterns, of size $O ( 2^{0.94995k})$ out of $2^k$ possible patterns.
One can show a similar result for decrease by a factor of  any fixed $\lambda$, 
and  a proof of natural density  for general $\lambda >0$ 
is given in  Borovkov and Pfeifer \cite[Theorem 3]{BP00}. $~~~\bsq$\\


The results above are rigorous results, and therefore we have no compelling need to 
find stochastic models to model the behavior of stopping times. 
Nevertheless stochastic models intended to analyze  other statistics produce 
in passing 
models 
for stopping time distributions. In 
\S\ref{sec3p1} we present such a
model, which gives an interpretation of these stopping time densities as exact
probabilities of certain events. 


\paragraph{Remark.}
The  analysis in Theorem~\ref{th22}  treats $\lambda$ as fixed. In fact one can also 
prove rigorous results which allow  $\lambda$ to vary  slowly (as a function of $n$), 
 under the restriction that  $\lambda \le \log_2 n.$

%
%

\subsection{$3x+1$ Stopping Time Statistics: Total Stopping Times}

 The following concept  concerns the speed at which positive integers iterate
to $1$ under the map $T$, assuming  they eventually get there.

%
%
%

\begin{defi}\label{de231}
{\em 
The {\em total stopping time} $\sigma_{\infty}(n)$ for iteration of the $3x+1$ map
$T(\n)$ is defined for positive integers $\n$ by 
$$
\sigma_{\infty}(\n) : = \inf\{ k\ge 0:~T^{(k)}(n) = 1\}.
$$
We set $\sigma_{\infty}(n) = +\infty$ if no finite $k$ has this property.
}
\end{defi}

The $3x+1$ Conjecture asserts that all positive integers have a finite total stopping time.\\

Concerning lower bounds for this statistic, there are some rigorous results.
First, since each step decreases $n$ by at most a factor of $2$, we 
trivially
have
$$
\sigma_{\infty}(n) \ge \frac{\log n}{\log 2} \approx 1.4426 \log n.
$$

The strongest result on the existence of integers having a large total stopping time
is the following result of Applegate and Lagarias \cite[Theorem 1.1]{AL03}.

%
%

\begin{theorem}\label{th23}
{\em (Lower Bound for $3x+1$ Total Stopping Times)}
There are infinitely many $n$ whose total stopping time satisfies
\beql{239}
\sigma_{\infty}(n) \ge \left(\frac{29}{29 \log 2 - 14 \log 3}\right) \log n \approx 6.14316 \log n.
\eeq
\end{theorem}


Nothing has been rigorously proved about either 
 the average size of the total stopping time,
or about upper bounds for the total stopping time 
 (since such would imply the main conjecture!).
 This provides motivation 
 to study  stochastic models for this statistic, to make guesses
how it may behave. \\

The  various stochastic models discussed
in \S
\ref{sec3}, as well as empirical evidence
given below, suggest that the size of this statistic will always be proportional to 
$\log n$. This motivates the following definition.

%
%
%

 \begin{defi}~\label{de232}
 {\em 
For $n \ge 1$ the {\em scaled total stopping time}
$\gamma_{\infty}(n)$ of the
$3x+1$ function is given by 
\beql{401}
\gamma_{\infty}(n)  := \frac{\sigma_{\infty}(n)}{\log n}.
\eeq
This value will be finite for all positive $n$ only if  the $3x+1$ conjecture  is true. 
}
\end{defi}

 A stochastic model in 
 \S\ref{sec3}  makes strong predictions about the distribution of
 scaled total stopping times: 
 they should have a Gaussian distribution with
 mean
 $$
 \mu := \left(\frac{1}{2} \log \frac{4}{3}\right)^{-1} \approx 6.95212
$$
and variance 
$$
\sigma := \frac{1}{2} \log 3 \left(\frac{1}{2} \log \frac{4}{3}\right)^{\frac{3}{2}} , 
$$
cf. 
Theorem~\ref{th31}. In particular, half of all integers ought to 
have a total stopping time $\sigma_{\infty}(n) \ge \mu \log n \approx  6.95212 \log n.$
It seems {\bf scandalous} that there is no unconditional proof that infinitely
many $n$ have a stopping time at least this large, compared to the bound
\eqn{239} in Theorem~\ref{th23} above!\\

We next define a 
  limiting constant associated with extremal values of the scaled total
  stopping time for the $3x+1$ map.\\

%
%
%

\begin{defi}~\label{de233}
{\em  The $3x+1$ {\em scaled stopping constant} is the quantity
 \beql{404}
 \gamma = \gamma_3
 := 
 \limsup_{n \to \infty} \g_{\infty}(n)
 = 
 \limsup_{n \to \infty} \frac{\sigma_{\infty}(n)}{\log n}
 .
\eeq
}
\end{defi}

We now give empirical data  about these extremal values. 
Table \ref{tab21} presents empirical data on record holders for the function $\gamma_{\infty}(n)$,
compiled by  Roosendaal \cite{Roos}. 
This table also includes data on
another statistic called the {\em ones-ratio} (or {\em completeness}), taken
from Roosendaal \cite[Completeness and Gamma Records]{Roos}. 
The function  $\mbox{ones}$($n$) counts the
number of odd iterates of the $3x+1$ function to reach $1$ starting from $n$ (including $1$), and
\beql{151a} 
\mbox{ones-ratio}(n)  :=\mbox{ones}(n)/ \sigma_{ \infty}(n). 
\eeq
  Table \ref{tab21} shows that the function $\gamma(n)$ is not a monotone 
increasing function of the ones-ratio, compare rows 9 and 10. The values with
question marks mean that all intermediate values have not been searched, so
these values are not known to be record holders.\\

%
%
%

 \begin{table}\centering
\renewcommand{\arraystretch}{.85}
\begin{tabular}{|r|r|r|r|r|r|}
\hline
\multicolumn{1}{|c|}{$k$} &
\multicolumn{1}{c|}{$\mbox{\#k-th  record} ~~n_k$}
&\multicolumn{1}{c|}{$\sigma_{\infty}(n_k)$} &
\multicolumn{1}{c|}{$ones$($n_k$)} &
\multicolumn{1}{c|}{$ones-ratio$} &
\multicolumn{1}{c|}{$\gamma_{\infty}(n_k) $} \\ \hline
 1 &         3 &    5 &  2 & 0.400000 & 4.551196 \\ \hline
 2 &          7&    11&    5 &0.454545 &    5.652882 \\ \hline
 3 &           9 &    13 & 6& 0.461358 & 5.916555 \\ \hline
 4 &        27 &    70 &  41& 0.585714 & 21.238915\\ \hline
 5 &        230~631 &   278  & 164 &   0.589928 & 22.512720\\ \hline
 6 &        626~331 &   319 &  189&  0.592476 & 23.899366 \\ \hline
 7 &       837~799 &   329  &   195&  0.592705 & 24.122828 \\ \hline
 8 &       1~723~519&    349  & 207& 0.593123 & 24.303826 \\ \hline
 9 &       3~732~423&    374  &  222& 0.593583 & 24.714906 \\ \hline
10 &       5~649~499 &    384 &  228&  0.593750 &  24.699176  \\ \hline
11 &      6~649~279 &   416  &   248&0.596154 & 26.479917 \\ \hline
12 &      8~400~511&   429  &   256& 0.596737 & 26.907006\\ \hline
13 &      63~728~127 &   592  &  357 & 0.603041 & 32.943545 \\ \hline
14 &     3~743~559~068~799 &   966  &  583  & 0.603520 & 33.366656 \\ \hline
15 &     100~759~293~214~567  &   1134  &  686 &0.604938 & 35.169600\\ \hline
?16 &     104~899~295~810~901~231&   1404  &  850&0.605413 & 35.823841 \\ \hline
?17 &     268~360~655~214~719~480~367&   1688  & 1022&0.605450 & 35.885221 \\ \hline
?18 &    6~852~539~645~233~079~741~799 &  1840  & 1115&  0.605978 & 36.595864 \\ \hline
?19 &    7~219~136~416~377~236~271~195&  1848  &   1120&  0.606061 & 36.716918 \\ \hline
\end{tabular}
\caption{Record Values for $\gamma_{\infty}(n)$ and for ones-ratio(n). 
  }
\label{tab21}
\end{table}
%
%
%

In 
\S\ref{sec4} we present a stochastic model which makes a prediction for
the extremal value of $\gamma$.
A quite different model is discussed in
\S\ref{sec7}, 
which makes exactly the same prediction! 
For both 
models 
the analogue of the constant 
 $\gamma:=\limsup \gamma_{\infty}(n) $ exists and equals a constant 
 which numerically is approximately $41.677647$, with  corresponding ones-ratio of
 about $0.609091$.
 Compare these predictions with the data in Table \ref{tab21}. 
%
%

\subsection{$3x+1$ Size Statistics: Maximum Excursion Values}

Another interesting  statistic is the maximum value attained in a trajectory, which
we call the maximum excursion value. \\
%
%
%

\begin{defi}\label{de241}
{\em 
The {\em maximum excursion value} $t(n)$
 is the maximum value occurring in the forward iteration of
the integer $n$, i.e.
\beql{241a}
t(n) := \max (T^{(k)}(n):~ k \ge 0),
\eeq
with $t(n)= +\infty$ if the trajectory is divergent.
}
\end{defi}

The quantity $t(n)$ will be finite for all $n$ if and only if there are no
divergent trajectories for the $3x+1$ problem 
(but does not exclude the possibility of as yet unknown loops)
.\\

We define the following extremal statistic for 
maximum excursions. 
%
%
%

\begin{defi}~\label{de242}
{\em 
Let 
the $3x+1$ {\it maximum excursion ratio} be given by
\beql{242aa}
\rho(n)
:=
{
\log t(n) 
\over
 \log n
}
.
\eeq
%
%
Then the $3x+1$ {\em maximum excursion constant} is the quantity
\beql{242a}
\rho
:=
 \limsup_{n \to \infty}\rho(n)
=
 \limsup_{n \to \infty} \frac{\log t(n)}{\log n}
.
\eeq
}
\end{defi}

%
%
%
  \begin{table}\centering
\renewcommand{\arraystretch}{.85}
\begin{tabular}{|r|r|r|r|r|}
\hline
\multicolumn{1}{|c|}{$k$} &
\multicolumn{1}{c|}{$\mbox{\#k-th  record} ~n_k^{\ast}$}
&\multicolumn{1}{c|}{$t(n_k^{\ast})$} &
\multicolumn{1}{c|}{$r(n_k^{\ast}) $} &
\multicolumn{1}{c|}{$\rho(n_k^{\ast}) $} \\ \hline
1 &          2 &    2 & ~0.500 &          1.000 \\ \hline
 2 &         3 &     8 & 0.889 &       1.893      \\ \hline
 3 & 7 & 26 & 0.531 &                 1.674       \\ \hline
 4 & 15 & 80 & 0.356 &               1.618         \\ \hline
 5 &          27&    4 616&    6.332 &   2.560  \\ \hline
 6 & 255 & 6 560 & 0.101  &             1.586   \\ \hline
 7 & 447 & 19 682 & 0.099  &           1.620    \\ \hline
 8 & 639 & 20 782 & 0.051  &             1.539      \\ \hline
 9 & 703 & 125 252 & 0.253  &             1.792    \\ \hline
 10 & 1 819 &  638 468 & 0.193  &         1.781     \\ \hline
 11 & 4 255 & 3 405 068 & 0.188  &          1.800    \\ \hline
 12 & 4 591 & 4 076 810 & 0.193  &             1.805     \\ \hline
 13 & 9 663 & 13 557 212 & 0.145  &              1.790    \\ \hline
 14 & 20 895 & 25 071 632 & 0.057 &            1.712         \\ \hline
 15 & 26 623 & 53 179 010 & 0.075  &                 1.746     \\ \hline
 16 & 31 911 & 60 506 432  & 0.059  &                  1.728   \\ \hline
  17 & 60 975 & 296 639 576 & 0.080  &             1.771  \\ \hline
 18 & 77 671 & 785 412 368  & 0.130  &               1.819    \\ \hline
 19 & 113 383 & 1 241 055 674 & 0.097  &              1.799    \\ \hline
 20 & 138 367 & 1 399 161 680 & 0.073  &                1.779    \\ \hline
 21 & 159 487 & 8 601 188 876 & 0.338  &                 1.861 \\ \hline
 22 & 270 271 & 12 324 038 948 & 0.169  &            1.858  \\ \hline
 23 & 665 215 & 26 241 642 656 & 0.059 &              1.789   \\ \hline
 24 & 704 511 & 28 495 741 760 & 0.057  &              1.788  \\ \hline
 25 & 1 042 431 & 45 119 577 824 & 0.042  &            1.770  \\ \hline

\end{tabular}
\caption{Seeds $n$ giving record heights for $3x+1$ maximum excursion value $t(n)$.}
\label{tab22}
\end{table}

%
%
%

\begin{figure}
\centering
\includegraphics[width=3.5in]{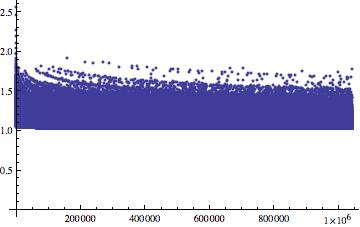}
\caption{
A plot of $n$ versus the maximal excursion ratio $\rho(n)$ for $3\le n\le 1\, 042\, 431$ and odd, cf. \eqn{242aa}. The only seeds $n$ in this range with $\rho(n)>2$ are $n=27,~ 31,~ 41,~ 47,~ 55,$ and $63$ (which all look at this scale as if they are on the $y$-axis).
}
\label{fignRho}
\end{figure}

%
%
%

%
%
%
\begin{table}\centering
\renewcommand{\arraystretch}{.85}
\begin{tabular}{
|r|r|r|r|}
\hline
\multicolumn{1}{|c|}{$
n
$}
&\multicolumn{1}{c|}{$t(n
)$} &
\multicolumn{1}{c|}{$r(n
) $} 
&
\multicolumn{1}{c|}{$\rho(n
) $} 
\\ \hline

27& 4 616& 6.332& 2.560 \\ \hline

319 804 831& 707 118 223 359 971 240&
  6.914& 2.099 \\ \hline
  
 1 410 123 943& 3 562 942 561 397 226 080& 1.792&
  2.028 \\ \hline
  
  3 716 509 988 199& 103 968 231 672 274 974 522 437 732& 7.527& 
  2.070 \\ \hline
  
  9 016 346 070 511& 126 114 763 591 721 667 597 212 096& 1.551&
  2.015 \\ \hline
  
   1 254 251 874 774 375& 1 823 036 311 464 280 263 720 932 141 024&
  1.159& 2.004 \\ \hline
  
  1 980 976 057 694 848 447& 
  32 012 333 661 096 566 765 082 938 647 132 369 010& 
  8.158& 2.050 \\ \hline

\end{tabular}
\caption{Values of $n
$
for which the maximal excursion ratio $\rho(n)={\log t(n)\over \log n}>2$ (equivalently, $r(n)= t(n)/n^2>1$),
culled from Oliveira e Silva's \cite[Table 8]{OeS09} 
record $t(n)$ values.}
\label{tab42}
\end{table}

%
%
%

The maximal excursion constant will be $+\infty$ if there is a divergent trajectory. 
The fact that the logarithmic scaling used in defining this constant is
the ``correct" scaling is justified by empirical data given in 
Oliveira e Silva \cite
{OeS09} (in this volume) and by the
predictions of the stochastic model given in \S
\ref{sec3}.
As explained in \S
\ref{sec4p3}
, the stochastic model prediction
for the maximum excursion constant is $\rho=2.$ \\

In Table~\ref{tab22}  we give 
the set of initial champion values for the maximum excursion,
extracted from data of Oliveira e Silva \cite{OeS99}.
 For comparison we give for each the ratio
$r(n) := \frac{t(n)}{n^2}$ and the value of the maximal excursion ratio $\rho(n)={\log t(n)\over \log n}.$
It is also useful to examine the larger table to $10^{18}$ given in Oliveira e Silva \cite{OeS09}.\\

While record values of $t(n)$ have received tremendous computational attention, there has not been a substantial amount of effort put into congregating those $n$ with large  $\rho(n)$ 
(the difference being that the former seeks seeds $n$ with large values of $t(n)$, whereas the latter seeks large values of $t(n)$ {\it relative} to the size of  $n$). 
We have computed 
that the only seeds $n<10^{6}$ for which $\rho(n)>2$ 
are: $n\in\{27, 31, 41, 47, 55, 63\}$, cf. Figure \ref{fignRho}. 
\\

Nevertheless, some ``large'' values of $\rho(n)$ 
already appear in tables of large $t(n)$'s.
In Table \ref{tab42} 
 we
 extract 
 from a table of $t(n)$ champions computed by 
 Oliveira e Silva  \cite{OeS09} 
 the subset of seeds $n$  for which $\rho(n)>2$, i.e. 
 $$
 r(n) = \frac{t(n)}{n^2}>1.
 $$
 Only seven
 such values appear. 
 This
 data  seems to (however weakly) support Conjecture \ref{conj42}.\\

%
%
\subsection{$3x+1$ Count Statistics: Inverse Iterate Counts}\label{sec2p5}

In considering backwards iteration of the $3x+1$ function, we
can ask: given an integer $a$ how many numbers $n$ 
have 
$T^{(k)}(n)=a$,
that is,
iterate forward
to $a$ after exactly $k$ iterations?  \\

The set of backwards iterates of a given number $a$ can be pictured as a tree; we call
these {\em $3x+1$ trees} and describe their structure in 
\S\ref{sec6}. Here
$N_k(a)$ counts the number of leaves at depth $k$ of a tree with root node $a$, 
and $N_k^{\ast}
(a)$ counts the number of leaves in a {\em pruned $3x+1$ tree},
in
which all nodes with label $n \equiv 0 ~(\bmod ~3)$ have
been removed. The definitions are as follows.\\

%
%
%

\begin{defi}~\label{de251}
{\em 
 (1) Let $N_k(a)$ count the number of integers that forward iterate under the
 $3x+1$ map $T(n)$ to $a$ after exactly $k$ iterations, i.e. 
\beql{251a}
N_k(a) := |\{ n: ~T^{(k)} (n) =a \}|.
\eeq

(2) Let $N_k^{\ast} (a)$ count the number of integers not divisible by $3$ that forward iterate under the
 $3x+1$ map $T(n)$ to $a$ after exactly $k$ iterations, i.e. 
\beql{251b}
N_k^{*}(a) := |\{ n: ~T^{(k)} (n) =a, ~n \not\equiv 0 (\bmod~3) \}|.
\eeq
}
\end{defi}

The case $a=1$ is of particular interest, since the quantities then count integers that
iterate to $1$
.
We 
set 
$$N_k := N_k(1), ~~~~~~~N_k^{\ast} := N_k^{\ast}(1).
$$
The secondary quantity $N_k^{\ast}(a)$ is introduced because it is somewhat more convenient
for analysis. It satisfies the monotonicity properties 
$N_{k}^{\ast}(a) \le N_{k+1}^{\ast}(a)$ and 
$$
N_k^{\ast}(m) \le N_k(a) \le \sum_{j=0}^{k} N_j^{\ast}(m) \le (k+1) N_k^{\ast}(a).
$$

We have the 
trivial
exponential upper bound 
\beql{253}
N_k(a) \le 2^k.
\eeq
since each number has at most $2$ preimages. 
We 
are
interested
in 
the exponential growth rate of $N_k(a)$. 

%
%
%

\begin{defi}\label{d252}
{\em
(1) For a given $a$
the  $3x+1$ {\em tree growth constant} $\delta_3(a) $  
is given by
\beql{255}
\delta_3(a) :=\limsup_{k \to \infty} \frac{1}{k}\left( \log N_k(a)\right).
\eeq

(2) The  $3x+1$ {\em 
universal 
 tree growth constant} 
 is
 $\delta =\delta_3= \delta_3(1).$
}
\end{defi}

The constant $\delta_3(a)$ exists and is finite, as follows from the upper bound \eqn{253}.
 It is easy to prove
unconditionally that
$\delta_3(3a) =0,
$
because the only preimages of a number $3a$ are $2^k 3 a$ and $N_k( 3a) =1$
for all $k \ge 1$.  The interesting case is when $a \not \equiv 0~(\bmod~3)$. \\

Applegate and Lagarias \cite{AL95a} determined
by computer the maximal and minimal number of leaves in 
pruned
$3x+1$ trees of depth $k$ for $k 
\le
30$.
The maximal and minimal number of leaves in such trees at level $k$ is given by
$$
N_k^{+} := \max \{  N_k^{\ast}( a):  ~a~(\bmod~3^{k+1}) ~\mbox{with}~ a \not\equiv 0~(\bmod~3) \}
$$
and
$$
N_k^{-} := \min\{ N_k^{\ast}( a):  ~a~(\bmod~3^{k+1}) ~\mbox{with} ~a \not\equiv 0~(\bmod~3) \},
$$
respectively. 
Counts for the number of leaves
in maximum and minimum size trees
of various depths $k$ are given  in the following table, taken from Applegate and Lagarias (\cite{AL95a},
\cite{AL95c}).  It is known that the average number of leaves at depth $k$ (averaged over $a$)
is proportional to $\left(\frac{4}{3}\right) ^{k}$, therefore in Table \ref{tab23} below we include
the value $(\frac{4}{3})^k$ and the scaled statistics
$$
D_k^{\pm} := N_k^{\pm} \left(\frac{4}{3}\right) ^{-k}. 
$$
This table also 
gives the number of distinct types of trees of each depth 
(there are some symmetries which speed up the calculation)%
.\\

%
%
%

\begin{table}\centering
\renewcommand{\arraystretch}{.85}
\begin{tabular}{|r|r|r|r|r|r|r|}
\hline
\multicolumn{1}{|c|}{$k$} &
\multicolumn{1}{c|}{$\begin{array}{c}\mbox{\# tree types}
\end{array}$} &
\multicolumn{1}{c|}{$N_k^{-}$} &
\multicolumn{1}{c|}{$N_k^{+}$} &
\multicolumn{1}{c|}{\rule[-2ex]{0cm}{6ex}$\left(\frac{4}{3}\right)^k$} &
\multicolumn{1}{c|}{$D_k^{-}$} &
\multicolumn{1}{c|}{$D_k^{+}$}\\ \hline
 1 &         4 &    1 &    2 &    1.33 & 0.750 & 1.500 \\ \hline
 2 &         8 &    1 &    3 &    1.78 & 0.562 & 1.688 \\ \hline
 3 &        14 &    1 &    4 &    2.37 & 0.422 & 1.688 \\ \hline
 4 &        24 &    2 &    6 &    3.16 & 0.633 & 1.898 \\ \hline
 5 &        42 &    2 &    8 &    4.21 & 0.475 & 1.898 \\ \hline
 6 &        76 &    3 &   10 &    5.62 & 0.534 & 1.780 \\ \hline
 7 &       138 &    4 &   14 &    7.49 & 0.534 & 1.869 \\ \hline
 8 &       254 &    5 &   18 &    9.99 & 0.501 & 1.802 \\ \hline
 9 &       470 &    6 &   24 &   13.32 & 0.451 & 1.802 \\ \hline
10 &       876 &    9 &   32 &   17.76 & 0.507 & 1.802 \\ \hline
11 &      1638 &   11 &   42 &   23.68 & 0.465 & 1.774 \\ \hline
12 &      3070 &   16 &   55 &   31.57 & 0.507 & 1.742 \\ \hline
13 &      5766 &   20 &   74 &   42.09 & 0.475 & 1.758 \\ \hline
14 &     10850 &   27 &  100 &   56.12 & 0.481 & 1.782 \\ \hline
15 &     20436 &   36 &  134 &   74.83 & 0.481 & 1.791 \\ \hline
16 &     38550 &   48 &  178 &   99.77 & 0.481 & 1.784 \\ \hline
17 &     72806 &   64 &  237 &  133.03 & 0.481 & 1.782 \\ \hline
18 &    137670 &   87 &  311 &  177.38 & 0.490 & 1.753 \\ \hline
19 &    260612 &  114 &  413 &  236.50 & 0.482 & 1.746 \\ \hline
20 &    493824 &  154 &  548 &  315.34 & 0.488 & 1.738 \\ \hline
21 &    936690 &  206 &  736 &  420.45 & 0.490 & 1.751 \\ \hline
22 &   1778360 &  274 &  988 &  560.60 & 0.489 & 1.762 \\ \hline
23 &   3379372 &  363 & 1314 &  747.47 & 0.486 & 1.758 \\ \hline
24 &   6427190 &  484 & 1744 &  996.62 & 0.486 & 1.750 \\ \hline
25 &  12232928 &  649 & 2309 & 1328.83 & 0.488 & 1.738 \\ \hline
26 &  23300652 &  868 & 3084 & 1771.77 & 0.490 & 1.741 \\ \hline
27 &  44414366 & 1159 & 4130 & 2362.36 & 0.491 & 1.748 \\ \hline
28 &  84713872 & 1549 & 5500 & 3149.81 & 0.492 & 1.746 \\ \hline
29 & 161686324 & 2052 & 7336 & 4199.75 & 0.489 & 1.747 \\ \hline
30 & 308780220 & 2747 & 9788 & 5599.67 & 0.491 & 1.748 \\ \hline
\end{tabular}
\caption{Normalized 
extreme 
values for $3x+1$ 
trees of depth $k$}
\label{tab23}
\end{table}

%
%
%

Applegate and Lagarias \cite[Theorem 1.1]{AL95a} proved 
the following result
by an easy induction 
using this table.
%
%

\begin{theorem}\label{th24}
{\em ($3x+1$ Tree Sizes)}
For any fixed $a \not\equiv 0~ (\bmod ~3)$ and for all sufficiently large $k$,
\beql{258}
(1.302053)^k \le N_k
^{\ast}
(a) \le (1.358386)^k.
\eeq
In consequence, for any $a \not\equiv 0 (\bmod ~3)$, 
\beql{259}
\log (1. 302053) \le \delta_3(a) \le \log (1.358386).
\eeq
\end{theorem}


We describe probabilistic models for $3x+1$ inverse iterates in 
\S\ref{sec6}.
The models are Galton-Watson processes 
for the number of leaves in the tree,
and branching random walks 
for the sizes of the labels in the tree. 
The model prediction is that $\delta_3(a) = \log \left( \frac{4}{3} \right)$ for all
$a \not\equiv 0 ~(\bmod~3)$.

%
%

\subsection{$3x+1$ Count Statistics: Total Inverse Iterate Counts}

In considering backwards iteration of the $3x+1$ function from an integer $a$, 
complete data  is the set of integers that  contain $a$ in their forward orbit. 
The $3x+1$ problem concerns exactly this question for $a=1$. The following
function describes this set. 

%
%
%

\begin{defi}~\label{de261}
{\em 
Given an integer $a$,  the {\em inverse iterate counting function}
$\pi_a(x)$ counts the number of integers $n$
with $|n|\le x$ that contain $a$ in their forward orbit under the $3x+1$ function.
That is, 
\beql{261a}
\pi_a(x) := \#\{ n:~ |n| \le x~~\mbox{and 
} ~T^{(k)}(n) = a
~~\mbox{for~some}~k \ge 0\}.
\eeq
}
\end{defi}

It is possible to obtain rigorous lower bounds for this counting function. 
For $a \equiv 0 ~(\bmod ~3)$ the set of inverse iterates is exactly
$\{ 2^k a:~ k \ge 0\}$ and $\pi_{a} ( x) = \lfloor \log_2 (\frac{2x}{|a|})\rfloor$
grows logarithmically.
 If $a \not\equiv 0~(\bmod~3)$ then $\pi_a(x)$ satisfies
a bound $\pi_a(x) > x^c$ for some positive $c$, as was first shown by
Crandall \cite{Cra78} in 1978.  The strongest method currently known to
obtain lower bounds on $\pi_a(x)$ was initiated
by Krasikov \cite{Kra89} in 1989, and extended in \cite{AL95b}, \cite{KL03}.
It gives the following result.

%
%

\begin{theorem}\label{th25}
{\em (Inverse Iterate Lower Bound)}
For each  $a \not\equiv~0~(\bmod~3)$, there is a positive constant $x_0(a)$
such that for all $x \ge x_0(a)$, 
\beql{265}
\pi_a(x) \ge x^{0.84}.
\eeq
\end{theorem}

\paragraph{Proof.} This is proved in Krasikov and Lagarias \cite{KL03}. 
The proof uses systems of difference inequalities $(\bmod~3^k)$,
analyzed in Applegate and Lagarias \cite{AL95b}, and by
increasing $k$ one gets better exponents. The exponent above was obtained
by computer calculation using $k=9$.
\hfil $\bsq$\\

  The following statistics measure the size of the inverse iterate set in
  the sense of fractional dimension.
  
%
%
%

\begin{defi}~\label{de262}
{\em 
Given an integer $a$,  the {\em upper and lower $3x+1$ growth exponents} for
$a$ are given by 
$$
\eta_3^{+}(a) := \limsup_{ x \to \infty}  \frac{ \log \pi_a(x)} {\log x},
$$
and
$$
\eta_3^{-}(a) := \liminf_{ x \to \infty}  \frac{ \log \pi_a(x)} {\log x}.
$$
If  these quantities are equal, we define the {\em $3x+1$ growth exponent} $\eta_3(a)$
to be  $\eta_3(a) = \eta_3^{+}(a)= \eta_3^{-}(a)$.
}
\end{defi}

We clearly have $\eta_3(a)=0$ if $a \equiv 0~~(\bmod ~3)$.  For the remaining values
Applegate and Lagarias made the following conjecture.

%
%
\begin{conj}~\label{conj21}
{\em ($3x+1$ Growth Exponent Conjecture)} 
For all integers $a \not\equiv 0~(\bmod ~3)$, the $3x+1$ growth exponent
$\eta_3(a)$ exists, with 
\beql{266}
\eta_3(a) = 1.
\eeq
\end{conj}

The truth of the $3x+1$ Conjecture would imply  that $\eta_3(1)=1$; however it
does not seem to determine $\eta_3(a)$ for all such $a$.
Applegate and Lagarias \cite[Conjecture A]{AL95a} made the stronger conjecture
that for each $a \not\equiv 0~(\bmod~3)$ $\pi_a(x)$ grows linearly, i.e.
there is a constant $c_a>0$
such that $\pi_{a}(x) > c_a x$ holds for all $x \ge 1$.\\

Note that Theorem~\ref{th25} shows that $\eta_3^{-}(a) \ge 0.84$ when
$a \not\equiv 0 ~(\bmod~3).$ Thus the lower bound in 
Conjecture \ref{conj21} thus  seems approachable.
A stochastic model in 
\S\ref{sec6p5} makes the prediction that $\eta_3(a)=1$.\\
%
%
%
\section{$3x+1$ Forward Iteration: Random Product and Random Walk Models}\label{sec3}
\hsp

In this section we formulate stochastic models intended to
predict
the behavior of iterations 
of the $3x+1$ map $T(\n)$
on a ``random" starting value 
$n$. These models are exactly
analyzable. We describe results obtained for these models,
which can be viewed as predictions for the ``average" behavior of
the $3x+1$ function. \\

%
%

\subsection{Multiplicative Random Product Model and $\lambda$-stopping times}\label{sec3p1}


Recall that the $\gl$-stopping time is defined 
(see \eqn{231}) 
by
$$
\sigma_{\lambda}(n) := \inf \{k \ge 0:  \frac{T^{(k)}(n)}{n}  < \lambda\}.
$$
Rawsthorne \cite{Raw85} and Borovkov and Pfeifer \cite{BP00} obtained 
a probabilistic interpretation of the $\lambda$-stopping time, as follows. 
They  consider  a  stochastic model which studies 
the random products
$$
Y_k := X_1 X_2 \cdots X_k,
$$
in which the  $X_i$ are each
independent identically distributed (i.i.d.) random
variables $X_i$ having the discrete distribution
$$
X_i = \left\{
\begin{array}{cl}
\df{3}{2}  & \mbox{with~probability} ~~\frac{1}{2}, \\
~~~ \\
\df{1}{2} & \mbox{with~probability}~~\frac{1}{2}.\\
\end{array}
\right.
$$
We call this the {\em $3x+1$ multiplicative random product} ($3x+1$ MRP) model.\\

This model does not include the choice of
the
 starting value of the iteration, which would correspond
to $X_0$; the random variable $Y_k$ really models the {\em ratio}
$\frac{T^{(k)}(X_0)}{X_0}$. 
They define for $\lambda > 0$ the {\em $\lambda$-stopping time random variable}
\beql{301c}
V_{\lambda} (\omega) := \inf\{k :~ Y_k \le \lambda\},
\eeq
where $\omega=(X_1, X_2, X_3, \dots)$
denotes a sequence of random variables as above. This random vector $\omega$
will model the effect of choosing a random starting value $n=X_0$ in iteration of the $3x+1$ map. \\

This stochastic model can be used to exactly describe the density of 
$\lambda$-stopping times, as follows. 
Let $\PP[E]$ denote the
probability of 
an
event $E$. \\

%
%

\begin{theorem}~\label{th30} {\em ($\lambda$-Stopping Time Density Formula)}
For the $3x+1$ function $T(n)$ the
 natural density $\DD (S_{\lambda}(k))$ for integers having  $\lambda$-stopping
 time at most $k$  is given exactly  by the formula
\beql{302a}
\DD (S_{\lambda}(k)) = \PP[ V_{\lambda}(\omega)  \le k],
\eeq
in which $V_{\lambda}$ is the $\lambda$-stopping time random variable
in the $3x+1$ multiplicative random product (MRP) model.
\end{theorem}

\paragraph{Proof.} 
In 1985 Rawsthorne \cite[Theorem 1]{Raw85} proved a  
 weaker version of this result, with $\DD (S_{\lambda}(k)) $ replaced by the 
lower
asymptotic density $\underline{\DD} 
(S_{\lambda}(k))$.
The result, using
natural density, is a special case of  Borovkov and Pfeifer \cite[Theorem 3]{BP00}.
$~~~\bsq$\\

It is natural to apply the 
$3x+1$ MPR 
model 
with an initial condition added, 
 which is 
 a proxy for the 
expected  behavior of the total stopping time.
To do this  we must
allow variable $\lambda$ (as a function of $n$),
in a range of parameters where there is  no rigorous proof 
that
the
model behavior agrees with
that of  iteration of the map $T(n)$, namely for 
$\lambda = \alpha \log n$ with 
various $\alpha>1$. 
What is missing is a result saying that it accurately matches the
behavior of iteration of the $3x+1$ map.  \\

The  behavior of the resulting probabilistic model is rigorously analyzable,
as we 
discuss
in the next subsection, 
cf.~Theorem~\ref{th31} below. 

%
%

\subsection{Additive  Random Walk Model and Total Stopping Times}

The $3x+1$ iteration takes $x_0=n$ and $x_k= T^{(k)}(n).$ 
In studying the iteration, it is often more convenient to use a logarithmic scale
and set $y_k= \log x_k$ (natural logarithm) so that 
$$
y_k= \log x_k :=\log T^{(k)} (n).
$$
Then we have
\beql{211}
y_{k+1} =
\left\{
\begin{array}{cl}
y_k + \log \frac{3}{2}  + e_k & \mbox{if}~ x \equiv 1~~ (\bmod ~2 ) ~,  \\ 
~~~ \\
y_k + \log \frac{1}{2} & \mbox{if} ~~x \equiv 0~~ (\bmod~2) ~,
\end{array}
\right.
\eeq
with
\beql{212}
e_k:= \log \left( 1+ \frac{1}{3 x_k}\right).
\eeq
Here $e_k$ is small as long as $|x_k|$ is large.\\

Theorem ~\ref{th21} implies that if an integer is drawn at random from $[1, 2^k]$ then
its $k$-truncated parity sequence will be uniformly distributed in $\{ 0, 1\}^k$. In consequence,
equations \eqn{211} and \eqn{212} show that 
the quantities $\log T^{(k)}(n)$ (natural logarithm) 
can be modeled by
a random walk starting at
initial position $y_0 = \log n$ and taking
steps of size $\log \frac{3}{2}$ if the parity value is odd, and $\log \frac{1}{2}$ if it
is even. \\

%
%

The  
MRP
model considered before 
is converted
to an additive model by
making a logarithmic change of variable, taking new
random variables $W_k := \log X_k.$ The additive model  considers the 
random variables $Z_k$ which are a sum of random variables 
$$
 Z_k  := Z_0 + \log Y_k = Z_0+ W_1+ W_2 + \cdots + W_k
.
 $$
 Here $Z_0$ is a specified initial starting point, and 
  $Z_k$ is  the result of a (biased) random walk, taking steps of 
 size either $ \log \frac{3}{2}$ or $\log \frac{1}{2}$ with equal probability.  In
 terms of these variables, 
 the $\lambda$-stopping time random variable above is
 $$
 V_{\lambda}(\omega) = \inf\{ k: Z_k - Z_0\le \log \lambda\}.
$$

We consider the approximation of this iteration process by the
following stochastic model, which we term the {\em $3x+1$ Biased 
Random 
Walk 
Model }
({\em $3x+1$ BRW Model}).
For   an integer  $n \ge 1$ it separately makes a random walk which takes 
steps of size $\log \frac{1}{2}$ half the time and $\log \frac{3}{2}$ half the time. We can
write such a random variable as
$$
\xi_k  := -\log 2+ \delta_k  \log 3
,
$$
in which $\delta_k$
are  independent  Bernoulli zero-one random variables.
The random walk positions  $\{ Z_k: k \ge 0\},$ are then random variables having starting value
$Z_0 = \log n$, 
and with 
$$
Z_k := Z_0 + \xi_1 + \xi_2 + \cdots + \xi_k.
$$
The $Z_k$ define a biased random walk, whose expected drift $\mu$ is given by 
\beql{312a}
\mu: = E[\xi_k)] = - \log 2 + \frac{1}{2} \log 3= \frac{1}{2} \log \left( \frac{3}{4} \right) \approx -0.14384.
\eeq
The variance $\sigma$ of each step is given by
$$
\sigma:= {\rm Var}[ \xi_k] = \frac{1}{2} \log 3 \approx 0.54930.
$$

In
the addive model we 
 associate to a random walk  a {\em total stopping time random variable} 
$$
S_{\infty}(n) := \min \{ k >0: Z_k \le 0, ~\mbox{given}~Z_0= \log n\},
$$
which detects  when the walk  first crosses $0$ (this corresponds in the multiplicative 
model to reaching $1$). 
The expected number of steps to reach a nonpositive value starting from
$Z_0= \log n$ is
$$
E[ S_{\infty}(n) ]=\frac{1}{|\mu|} 
\log a
=  \frac{1}{\frac{1}{2} \log(\frac{4}{3})} \log n  \approx  6.95212 \log n.
$$

 As noted in \S\ref{sec2}, Borovkov and Pfeifer \cite{BP00} 
 consider the  multiplicative stochastic model obtained
by exponentiation of the positions of
the biased 
random walk above, from a given starting value
$X_0= e^{n_0}$.  They conclude the following result \cite[Theorem 5]{BP00}.
%
%

\begin{theorem}~\label{th31} {\em ($3X+1$ BRW Gaussian Limit Distribution)}
In the Biased Random Walk 
Model, 
for each fixed $n \ge 2$  define the  normalized  random variable
$$
Z_{\infty}(n) := \frac{S_{\infty}(n) - \frac{1}{\mu} \log n}{\mu^{- \frac{3}{2} }\sigma \sqrt{\log n}}, 
$$
which has cumulative distribution function $P_n(x):={\rm Prob}[ Z_{\infty}(n) < x]$. 
Here $\mu=|\frac12\log\frac34|$, and $\gs=\frac12\log3$.
Then 
for each 
fixed
 real $x$, allowing $n$ to vary, one has
$$
P_n(x):= 
{\rm
Prob
}
[ Z_{\infty}(n) < x] \longrightarrow  \Phi(x), ~~~{as}~~n \to \infty,
$$
where $\Phi(x)=\frac{1}{\sqrt{2 \pi}} \int_{- \infty}^x   e^{- \frac{1}{2} t^2} dt$ 
is the cumulative distribution function  of the standard normal 
distribution $N(0,1)$. 
\end{theorem} 


Borovkov and Pfeifer
note further that  the rate of convergence of the normalized distribution $P_n
(x)$
with fixed $n$ to the limiting normal distribution as $n \to \infty$ is uniform in $x$, but 
is quite slow.
They assert that for all $n \ge 2$ and all $-\infty < x < \infty$,
\beql{312b}
|P_n(x) - \Phi(x)| = O \left(  (\log n)^{- \frac{1}{2}} \right),
\eeq
where the implied constant in the O-symbol  is absolute.\\


%
%
%

\begin{figure}
\centering
$
{
\includegraphics[width=3in]{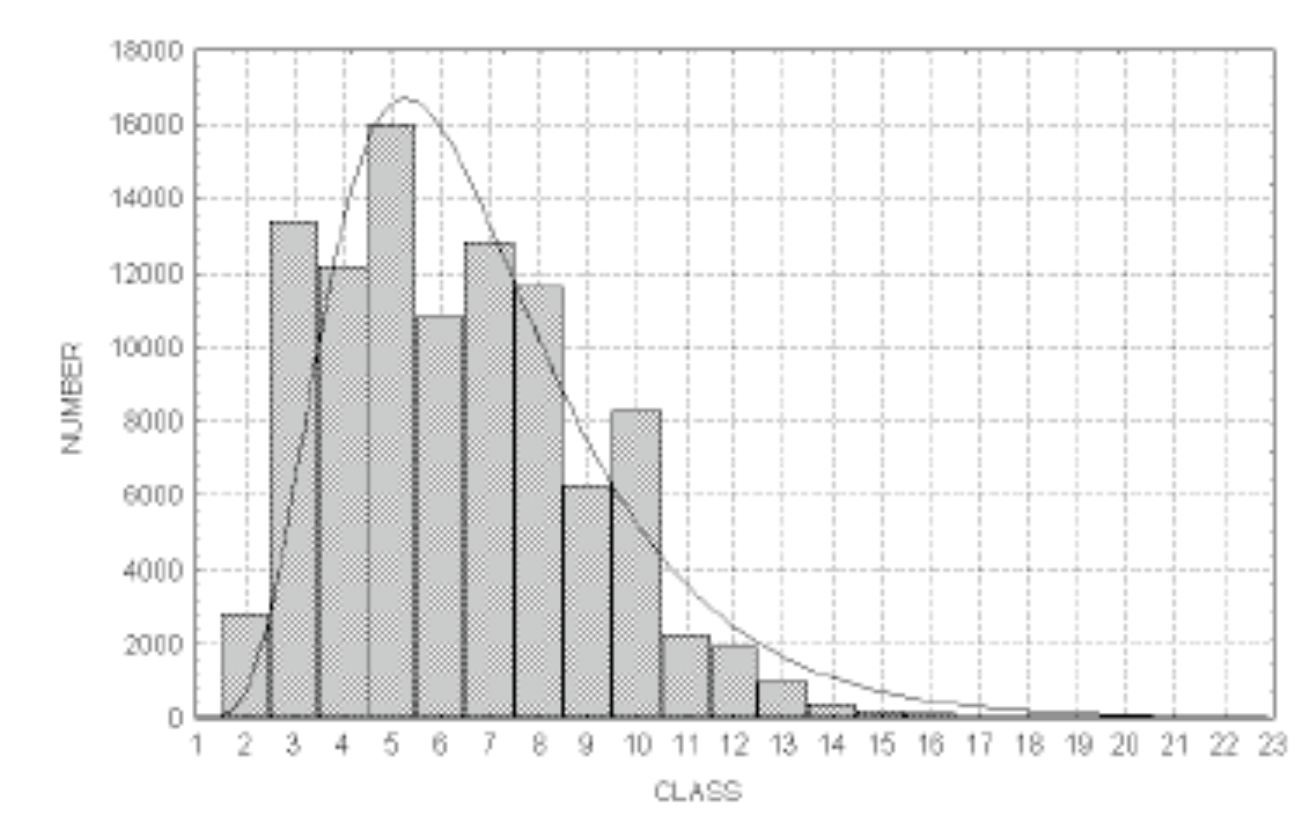}
\atop
x_{0} \ \in \ [0.95\times 10^{6}, \ 1.05\times 10^{6})
}
\qquad
{
\includegraphics[width=3in]{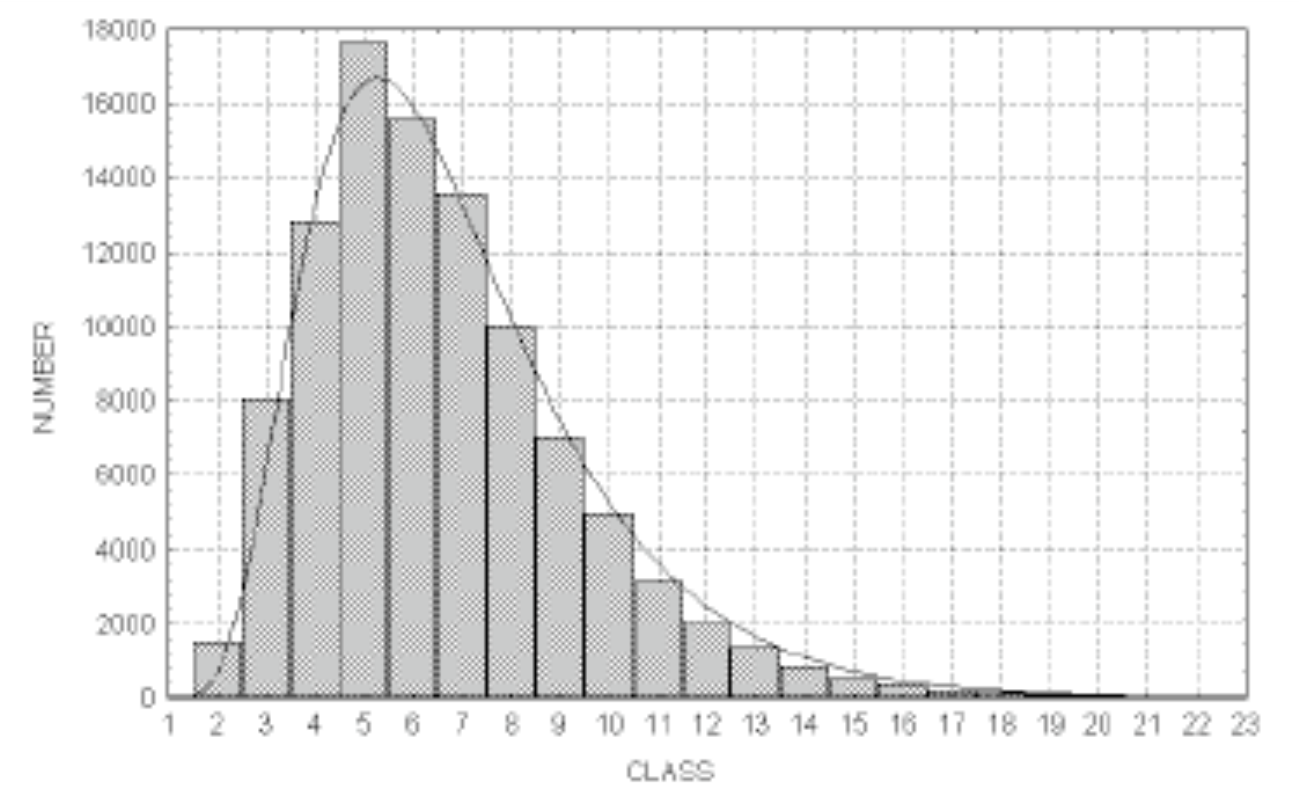}
\atop
n \ =  \ 10^{6}
}
$

\vskip.4in

$
{
\includegraphics[width=3in]{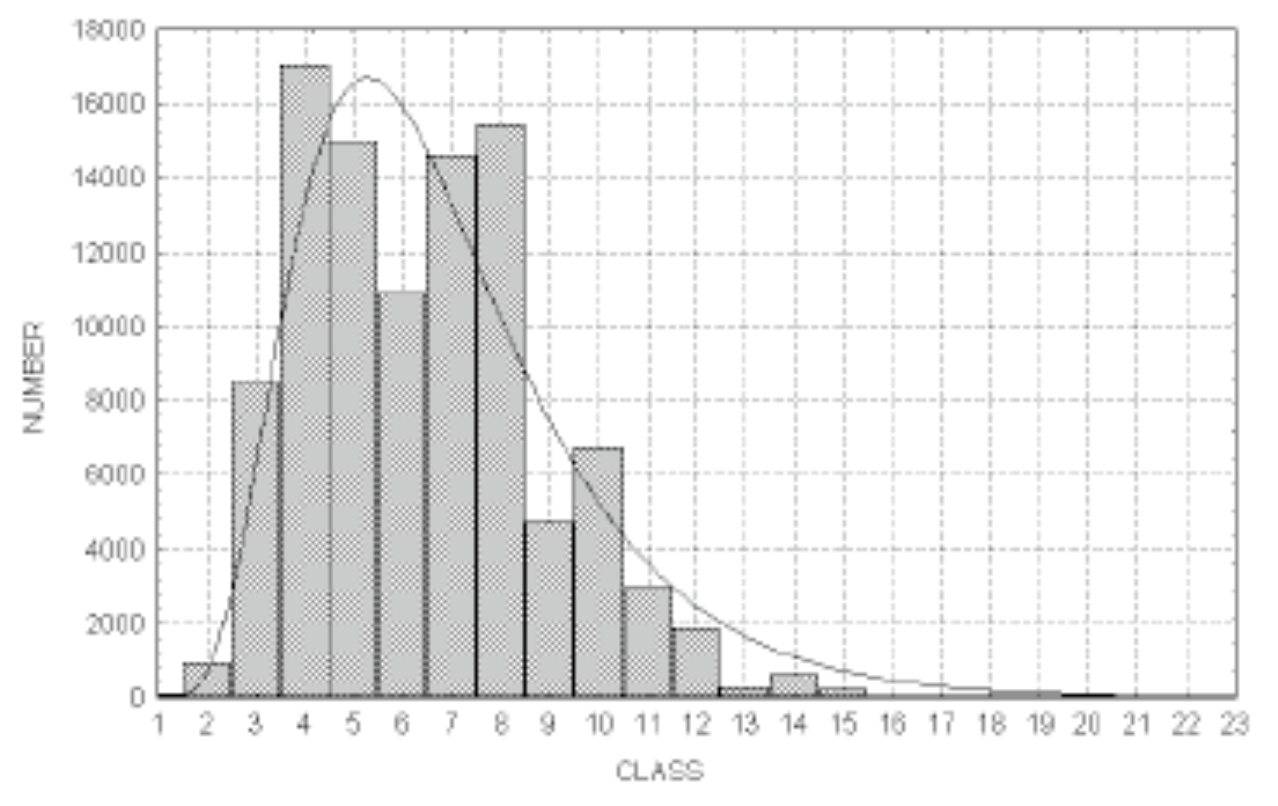}
\atop
x_{0} \  \in \ [1.95\times 10^{6}, \ 2.05\times 10^{6})
}
\qquad
{
\includegraphics[width=3in]{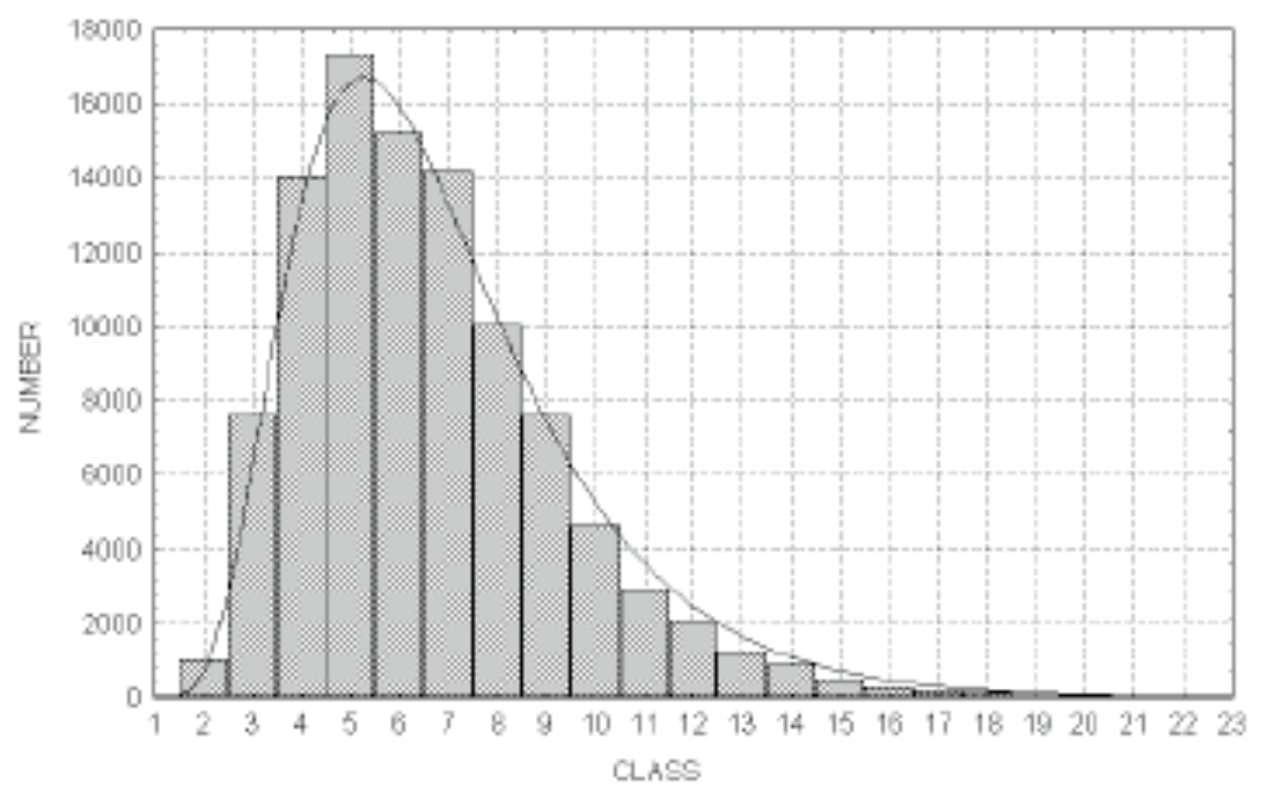}
\atop
n \ = \  2\times10^{6}
}
$

\vskip.4in

$
{
\includegraphics[width=3in]{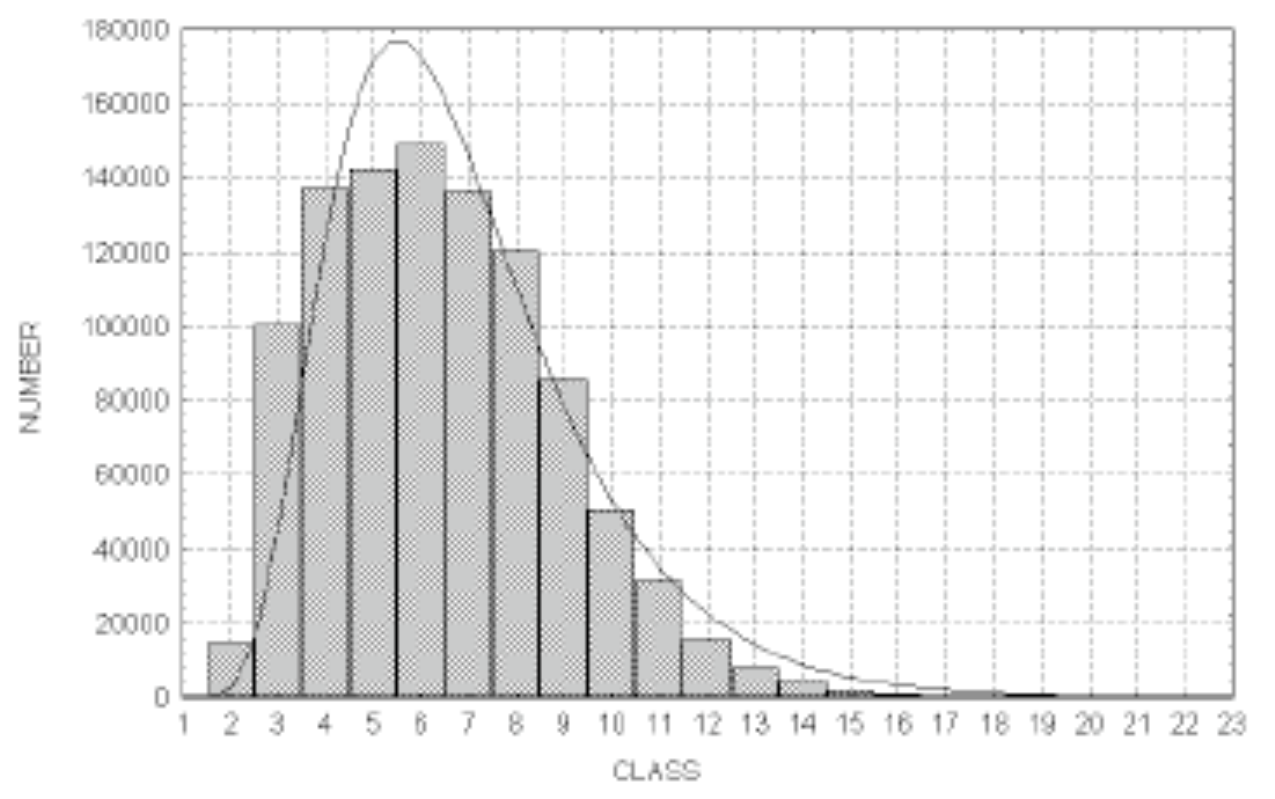}
\atop
x_{0} \ \in \ [3.5\times 10^{6}, \ 4.5\times 10^{6})
}
\qquad
{
\includegraphics[width=3in]{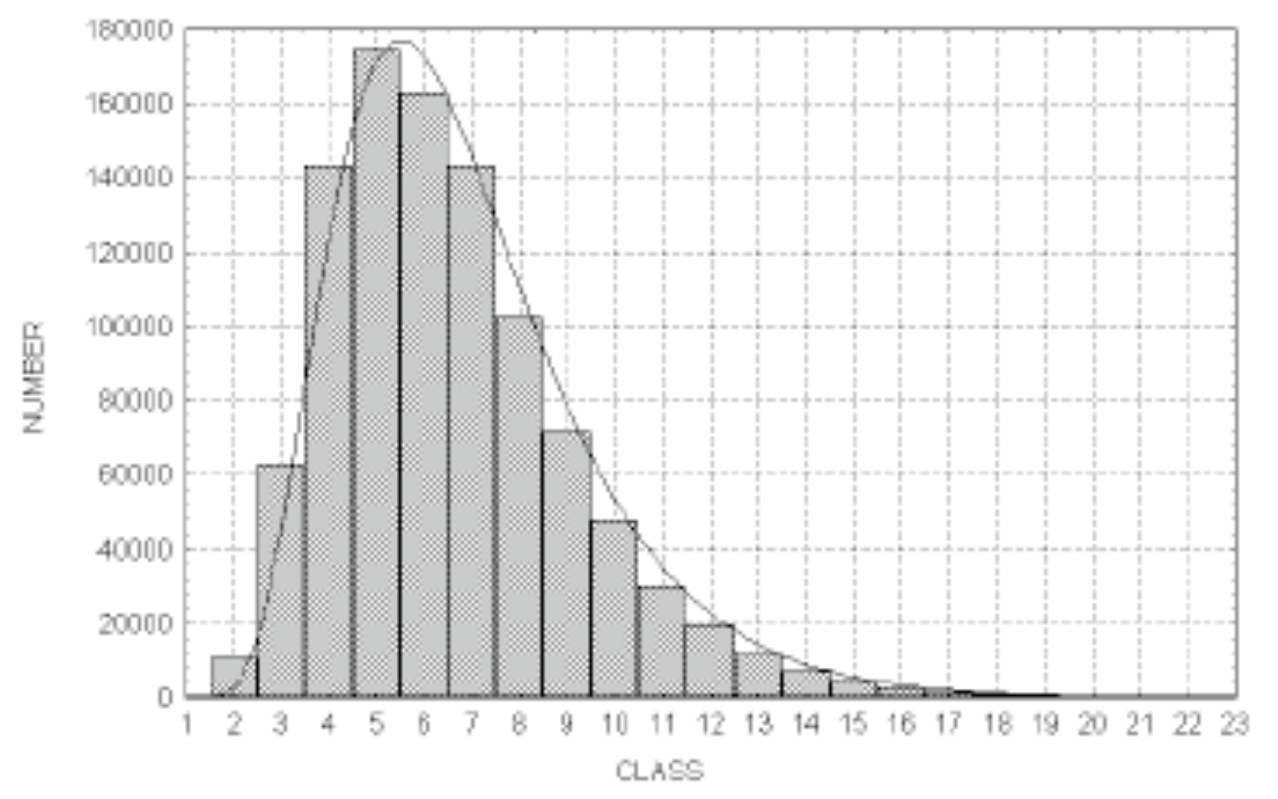}
\atop
n\ = \ 4\times 10^{6}
}
$

\vskip.4in

\caption{Histograms for $\sigma_{\infty}(x_{0})/ \ln x_{0}$ and 
its stochastic analog 
$T(n)/ \ln n$ with fitted density. 
Taken from Borovkov-Pfeifer \cite{BP00}.
}
\label{figBP}
\end{figure}


They also propose  a better approximation to the
distribution of the total stopping time of a random integer of size near $n
$,
reflecting the fact that it is  nonnegative random variable. 
They assert that the rescaled variable 
$$
Y_{\infty}(n) := \frac{S_{\infty}(n)}{\log n}
$$
should have a good second order approximation given by 
the nonnegative random variable $\tilde{Y}(n
)$ 
having the  distribution function 
$$
\Psi_n(x) =C_n  \frac{\sqrt{\log n}}{\sigma} \int_{0}^x  \frac {1}{\sqrt{2 \pi t^3}}
e^{- \frac{ (\mu t -1)^2 \log n} {2 \sigma^2 t}}dt ,  ~~x >0. 
$$
in which $C_n$ is a normalizing constant (\cite[eqn. (25)]{BP00}).\\

They view the random variable $S_{\infty}(n)$ as providing a model for the total stopping time
$\sigma_{\infty}(n)$ of the $3x+1$ function, where one compares the ensemble
of values $\{ \sigma_{\infty}(n):  x \le n \le c_1 x\}$   with $c_1>1$ fixed 
with independent samples of
values $S_{\infty}(n)$. The result above (with error term $O \left(\frac{1}{\sqrt{\log n}}\right)$)
predicts that for any $\epsilon >0$ the number
of values that do not satisfy
$$
\left(\frac{1}{\mu} - \frac{1}{ (\log x)^{\frac{1}{2} - \epsilon} } \right) \log x \le \sigma_{\infty}(n) \le 
 \left(\frac{1}{\mu} + \frac{1}{ (\log x)^{\frac{1}{2} - \epsilon} } \right) \log x
 $$ 
is $o(x)$, as $x \to \infty$.
They compare the distribution of  $\tilde{Y}(n)$ with numerical data $\frac{\sigma_{\infty}(n)}{\log n}$
 for the $3x+1$ function for $n \approx 10^6$
and find fairly good agreement. \\

%
%
%
\section{$3x+1$ Forward Iteration: Large Deviations and Extremal Trajectories}\label{sec4}
\hsp

 Lagarias and Weiss \cite{LW92} formulated and studied stochastic models which
 are intended to give predictions for the 
 extremal behavior of iteration of the $3x+1$ map $T(n)$. \\
 
%
%
%

 \subsection{$3X+1$ Repeated Random Walk Model}

 Lagarias and Weiss  studied the following {\em 
 Repeated 
 Random 
 Walk 
 Model} 
 ({\em $3x+1$ RRW 
 Model}).
For each integer $n\ge 1$, 
independently 
run 
a
 $3x+1$ biased random walk model  trial with starting
value $Z_{0,n} = \log n$. That is, 
generate 
an infinite
 sequence of independent random walks
$ \{Z_{k, n}: k \ge 0\}$, with one walk generated  for each value of $n$.
The model data
is the countable set of random variables 
\beql{400a}
\omega:= \{ Z_{k,n}:~ n\ge 1,~ k \ge 0 \},
\eeq
in which the initial starting points  $Z_{0,n}:= \log n$ are deterministic, and all other
random variables stochastic.
From this data,  one can  form random variables that are functions of $\omega$, corresponding
to the total stopping times and the maximum excursion values above. \\

The $3x+1$ RRW model is exactly analyzable, and makes predictions for
the value of the scaled stopping time constant, and for the maximum
excursion constant.  A  subtlety of the RRW model is the fact that
there are exponentially many trials with inputs of  a given length $j$, namely for those
$n$ with $e^j \le n < e^{j+1}$, which have initial
condition $j \le Z_{0,n} < j+1$, so that the theory of large deviations becomes relevant to
the analysis. \\

%
%
%

 \subsection{$3x+1$ RRW Model Prediction: Extremal Total Stopping Times}
 
 The $3x+1$ RRW model can be used to produce statistics analogous
to the scaled 
total
stopping time  $\gamma
_{\infty}
(n)$ and the 
$3x+1$ scaled stopping time constant $\gamma$
,
cf. \eqn{401} and \eqn{404}.
 \\

For a given trial $\omega$ it yields  an infinite sequence of total stopping time
random variables
$$
S_{\infty}(\omega):= (S_{\infty}(1), S_{\infty}(2), S_{\infty}(3),\dots, S_{\infty}(n),  \dots),
$$
where $S_{\infty}(n)$ is
computed using the  individual random walk $\sR_n$. 
Thus we can compute the scaled statistics $\frac{S_{\infty}(
n
)}{\log n
}$ for
$n
\ge 2$, and set
$$
\gamma(\omega) :=   \limsup_{n \to \infty} \frac{S_{\infty}(n)}{\log n}.
$$
as a stochastic  analogue of the quantity $\gamma$.\\

The $3x+1$ RRW model has the
following asymptotic limiting behavior for this statistic, given by
Lagarias and Weiss \cite[Theorem 2.1]{LW92}.
%
%

\begin{theorem}~\label{th41}
{\em ($3x+1$ RRW Scaled Stopping Time Constant)}
For the $3x+1$ 
RRW model, with probability one the scaled stopping time  
$$
\gamma(\omega) :=   \limsup_{n \to \infty} \frac{S_{\infty}(n)}{\log n}
$$
is finite and equals
a  constant
$$
\gamma_{\RW} \approx 41.677647, 
$$
which is the unique real number $\gamma> 
\left(\frac{1}{2} \log \frac{4}{3} \right)^{-1}\approx 6.952$
of the fixed point equation
\beql{412}
\gamma ~g\left(\frac{1}{\gamma}\right)=1.
\eeq
Here
the rate function $g(a)$ is given by
\beql{413}
g(a) := \sup_{\theta \in \RR} 
\left( \theta a - \log  M_{\RW}(\theta)\right),
\eeq
in which
\beql{423a}
M_{\RW}(\theta) := \frac{1}{2}\left(  2^{\theta} + \left(\frac{2}{3}\right)^{\theta} \right)
\eeq
is a moment generating function associated to the random walk.
\end{theorem}


Lagarias and Weiss also obtain a density result on the number 
of $n$ getting values
close to the extremal constant, as follows (\cite[Theorem 2.2]{LW92}).

%
%

\begin{theorem}~\label{th42} 
{\em ($3x+1$ RRW Scaled Stopping Time Distribution)}
For the $3x+1$
RRW model, and for any constant $\alpha$ satisfying
\beql{421}
\left(\frac{1}{2} \log \frac{4}{3}\right)^{-1} < \alpha < \gamma_{\RW},
\eeq
one has the bound
\beql{422}
E \left[ |\{ n \le x:~ \frac{S_{\infty}(n)}{\log n} \ge \alpha\}| \right] 
\le \left( 1 - \alpha ~g\left(\frac{1}{\alpha}\right) \right)^{-1} x^{1 - \alpha g(1/\alpha)}.
\eeq
In the reverse direction, for any $\epsilon>0$ this expected value satisfies
\beql{423}
E \left[ |\{ n \le x:~ \frac{S_{\infty}(n)}{\log n} \ge \alpha\}| \right] 
\ge x^{ 1- \alpha g(1/\alpha)- \epsilon}
\eeq
for all sufficiently large  $x \ge x_0(\epsilon).$
\end{theorem}


This theorem says that not only is there an upper bound $\gamma_{\RW}$
on the asymptotic limiting value of the stopping ratio, but the set of $n$ 
for which
one gets a value above $\alpha$ becomes very sparse 
(in the logarithmic sense) as $\alpha$
approaches $\gamma_{\RW}$ from below. Theorem~\ref{th42}
is analogous to obtaining a multifractal spectrum for this problem.
This result is well-suited for
comparison with experimental data on $3x+1$ iterates.\\

This analysis suggest the  following prediction, which we state as a conjecture.\\

%
%
%
\begin{conj} \label{conj41}
{\em ($3x+1$ Scaled Stopping Constant Conjecture)} 
 The $3x+1$ scaled stopping constant $\gamma$  is finite and is given by 
\beql{411}
\gamma  = \gamma_{\RW} \approx 41.677647.
\eeq
\end{conj}


The large deviations model does more than predict an extremal value, it also
predicts that the numbers that approach the extremal value must have
a trajectory of iterates whose graph have a specified shape, 
which is a linear function when properly scaled. 
In Figure \ref{figBP1}
we
graph the  set of scaled points
$$
\left\{\left(\frac{k}{\log n}, \frac{\log T^{(k)}(n)}{\log n}\right) : 0 \le k \le \sigma_{\infty}(n) \right\}.
$$
The  predicted large deviations extremal trajectory  in this scaling 
has graph a straight line connecting the
points $(0,1)$ and $(\gamma_{\RW}, 0)$. Figure ~\ref{figBP1} shows the 
scaled trajectories with starting seeds $n_{k}$ taken from Table \ref{tab21}, i.e. those with record
values for $\g_{\infty}(n)$. Compare to Lagarias and Weiss \cite[Figure 3]{LW92}.

%
%
%

\begin{figure}
\centering
\includegraphics[width=3.5in]{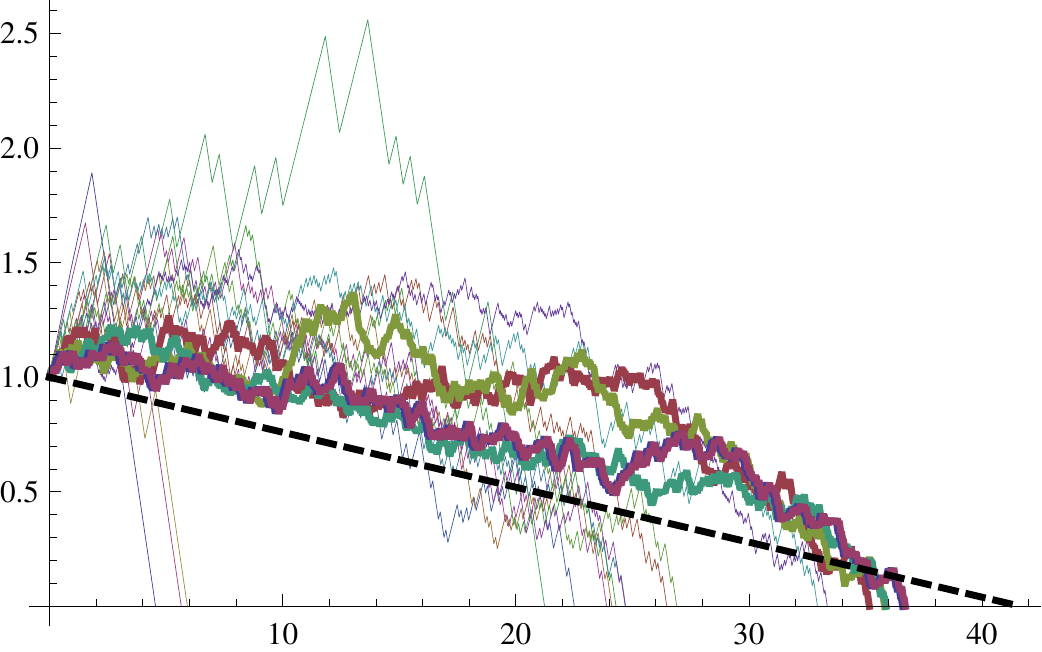}
\caption{
Scaled trajectories of $n_{k}$ maximizing $\gamma(n)$ 
for record values 
from Table \ref{tab21} (thin for $1\le k\le 10$; regular for $11\le k\le15$; thick for $16\le k\le 19$),
plotted against the predicted trajectory.
}
\label{figBP1}
\end{figure}

%
%
%

\subsection{$3x+1$ RRW Model Prediction: Maximum Excursion Constant}\label{sec4p3}

For the $3x+1$ RRW Stochastic Model, an appropriate statistic 
for a single trial that
corresponds to the maximum excursion value  is 
$$
t (n; \omega) := \sup ( e^{ Z_{k, n}} : k \ge  0)
.
$$

The $3x+1$ RRW model behavior
for extremal behavior of maximum excursions $t(n; \omega)$
is given in the following result \cite[Theorem 2.3]{LW92}.\\

%
%

\begin{theorem}~\label{th43} 
{(\em $3x+1$ RRW Maximum Excursion Constant)}
 For the $3x+1$
RRW model, with probability one   the quantities $t(n, \omega)$
are finite for every $n \ge 1$. In addition, with probability one
the random quantity 
\beql{431}
\rho(\omega):= \limsup_{n \to \infty} \frac{\log t(n; \omega)}{\log n}
= \limsup_{n \to \infty} \left( \sup_{k \ge 0} \frac{Z_{k, n}}{\log n} \right)
\eeq 
equals the constant
\beql{432}
\rho_{\RW} =  2.
\eeq
\end{theorem}


Lagarias and Weiss also 
prove
\cite[Theorem 2.4]{LW92} a result
permitting a quantitative comparison with data.

%
%

\begin{theorem}~\label{th44} 
{\em ($3x+1$  RRW Maximum Excursion Density Function)}
For the $3x+1$
RRW model, for any fixed $0 < \alpha <1$, the expected value
\beql{441}
E\left[ | \{ n \le x:~ \frac{ \log t(n; \omega)}{\log n} \ge 2 - \alpha\}|\right] = x^{\alpha( 1- o(1))},
\eeq
as $x \to \infty$.
\end{theorem}

These  theorems suggest formulating the  following conjecture.\\
%
%
%

\begin{conj}\label{conj42}
 The $3x+1$  maximum excursion constant $\rho$ 
 defined in \eqn{242a}
 is finite and is given by 
\beql{417}
\rho=2.
\eeq
\end{conj}

The large deviations model also makes a prediction on the graphs of the
trajectories achieving maximum excursion, when plotted as the scaled data
points 
$$
\left\{
\left(\frac{k}{\log n}, {\log T^{(k)}(n)\over \log n}
\right)
:~ 0 \le k \le \sigma_{\infty}(n) 
\right\}.
$$


It asserts that extremal large deviation trajectories should approximate 
two line
segments, the first with vertices $(0,1)$ to $(7.645, 2)$ and then from
this vertex to  $(21.55, 0).$ The slope of the first line segment
is $\frac{3}{4} \log 3 - \log 2 \approx 0.1308$
and that of the second line segment is $(\frac{1}{2} \log \frac{3}{4})^{-1} \approx -0.1453$.
This prediction is
shown as a dotted black line
in Figure \ref{figBP2};
 it shows substantial agreement with the empirical evidence. 

%
%
%

\begin{figure}
\centering
\includegraphics[width=3.5in]{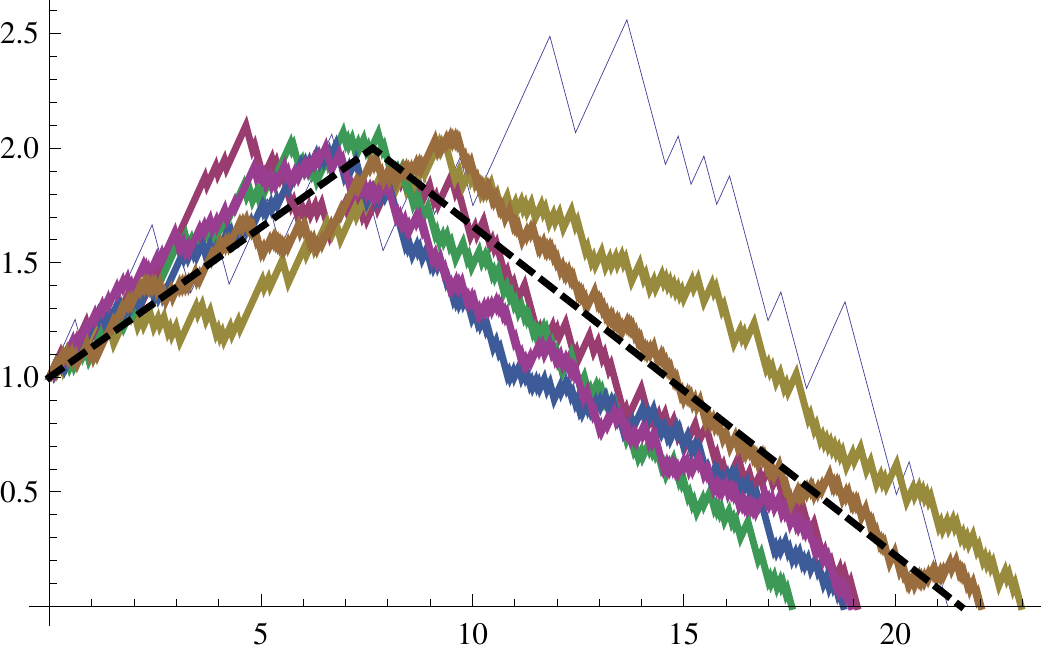}
\caption{
Scaled trajectories of seeds $n$ 
from Table \ref{tab42},
plotted against the predicted trajectory. The trajectory of $n=27$ is thin, while the others are thick.
}
\label{figBP2}
\end{figure}


%
%
%
\subsection{$3x+1$  RRW Model: Critique}\label{sec4p4}

 The $3x+1$
 repeated random walk model has the feature  that  random walks for different $n$ 
 are {\em independent. } However the actual $3x+1$ map certainly 
 has  a  great deal of dependency
 built in, due to the fact that trajectories coalesce under  forward iteration. 
 For example, trajectories of numbers $8n+4$ and $8n+ 5$ always coalesce after
 $3$ iterations of $T$. After coalescence,
 the trajectories are completely correlated. 
  In fact, the $3x+1$ Conjecture predicts that 
 all trajectories of positive integers $n$ reach the orbit $\{1, 2\}$ and then
 cycle,  whence they all should coalesce into exactly  two classes, 
 namely those that reach $1$ in an
 odd number of iterations of $T$, and those that reach this orbit under an even number
 of iterations. \\
 
 For this reason, it is not apparent a priori whether the prediction 
 in Conjecture \ref{conj41} above of the constant $\gamma =\gamma_{\RW}$  is  reasonable.
 Our faith in Conjecture \ref{conj41}  relies on the fact that first, 
 the same prediction is made using a branching random walk model that incorporates
 dependency in the model, see Theorem~\ref{th73} in \S\ref{sec6}, and second, on 
  comparison with empirical data in Table~\ref{tab21}
  .\\


%
%
%
\section{$3x+1$ Accelerated Forward Iteration : Brownian Motion}\label{sec5}
\hsp

Now we consider the accelerated $3x+1$ function $U
$. 
Recall that $U$ is defined on odd integers, and removes all powers of $2$ in one fell swoop. 
Iterates of the accelerated function $U$ are 
of course equivalent (from the point of view of the main conjecture) to those 
of $T$, but there are some subtle differences which make studying 
both points of view appealing.\\


%
%
%

For an odd integer $\n$, we let $\ord(\n)$ count the number of powers of $2$ dividing
$3n+1$, so that
\beql{500a}
\ord(\n):={\rm ord}_2(3\n+1). 
\eeq
Then the accelerated $3x+1$ function $U$ is given by:
\beql{500b}
U(\n):={3\n+1\over 2^{\ord(\n)}}.
\eeq

In analogy with the (truncated) parity sequence, cf. Definition \ref{de21}, we make the following
definition, giving a symbolic dynamics for the accelerated $3x+1$ map.
%
%
%
\begin{defi}~\label{de52}
{\em 
(i) For an odd integer  $\n$, define
the  {\em $\ord$-sequence} of $\n$ to be
\beql{501}
V(\n) := (\ord_{1}(\n) , \ord_{2}(\n), \ord_{3}(\n), \dots)
\eeq
where
$$
\ord_{\kk}(\n):= \ord(U^{(\kk)}(\n)),
$$
and  $U^{(\kk)}(\n)$ denotes the $\kk$-th iterate of $U$, as usual.
This is an infinite vector of positive integers.\\

(ii) For $\kk \ge 1$ the {\em $\kk$-truncated $\ord$-sequence} 
of $\n$ is:
\beql{502}
V^{[\kk]}(\n):= (\ord_{1}(\n) , \ord_{2}(\n),  \dots, \ord_{\kk}(\n))
\eeq
i.e. a vector giving the  initial segment of $\kk$ terms of $V(\n)$. 
}
\end{defi}

%
%
%
\begin{defi}~\label{de53}
{\em 
For an odd integer  $\n$ and $\kk\ge1$, let the {\em $k$-size} 
$\s_{\kk}(\n)$ be the sum of the entries in 
$V^{[\kk]}(\n)$, that is
$$
\s_{\kk}(\n):= \ord_{1}(\n) + \ord_{2}(\n)+ \cdots+\ord_{\kk}(\n).
$$
}
\end{defi}

%
%
%
\subsection{The Structure Theorem}

Notice that $U(\n)$ is not only odd, but is also relatively prime to $3$. Hence we lose no generality by restricting the domain for $U$ from
$\ZZ$
to   the (more natural) set $\Pi$ of positive integers prime to $2$ and $3$, i.e.
\beql{500f}
\Pi : = \{ n \in \ZZ:  \gcd(n, 6)=1 \}.
\eeq
Moreover, $\Pi$ is the disjoint union of $\Pi^{(1)}$ and $\Pi^{(5)}$, where $\Pi^{(\vep)}$ consists of positive integers congruent to $\vep~(\bmod ~6)$, $\vep=1$ or $5$.\\

%
%
%
\begin{defi}~\label{de54}
{\em 
Given $\vep=1$ or $5$, $\kk\ge1$, and a vector $
(\ord_{1},\dots,\ord_{\kk})$ of positive integers, let 
$$
\gS^{(\vep)}(\ord_{1},\dots,\ord_{\kk})
$$
 be the set of all $\n\in\Pi^{(\vep)}$ with $V^{[\kk]}(\n)=(\ord_{1},\dots,\ord_{\kk})$.
}
\end{defi}

The result analogous to Theorem \ref{th21} is given by Sinai \cite{Si03a} and 
Kontorovich-Sinai \cite{KS02}.
%
%

\begin{theorem}~\label{th51}
{\em (Structure Theorem for $\ord$ Symbolic Dynamics)}
Fix $\vep=1$ or $5$, and let
 $\n\in\Pi^{(\vep)}$.

(i) The $\kk$-truncated $\ord$-sequence $V^{[\kk]}(\n)$ of the first $\kk$ iterates of the accelerated map $U(\n)$
is periodic  in $\n$. Its period is $6\cdot 2^{\s}$, where 
$$\s=\s_{\kk}(\n) = \ord_1(n) + \ord_2(n) + \cdots + \ord_k(n).$$ 

(ii) For any 
$\kk\ge1$
and 
$\s\ge\kk$,
each of the 
$\left({\s-1\atop \kk-1}\right)$ 
possible vectors $(\ord_{1},\cdots,\ord_{\kk})$ with $\ord_{j}\ge1$ and  $\ord_{1}+\cdots+\ord_{\kk}=\s$ occurs
exactly once as $V^{[\kk]}(\n)$ for some $\n\in\Pi^{(\vep)}$  in the initial segment \break $1\le \n<6\cdot 2^{\s}$.   \\

(iii) The least element $\n_{0}\in\gS^{(\vep)}(\ord_{1},\dots,\ord_{\kk})$ satisfies $\n_{0}<6\cdot 2^{\s}$; moreover 
 $$
\gS^{(\vep)}(\ord_{1},\dots,\ord_{\kk})
=
\bigg\{
6\cdot 2^{\s}\cdot m+\n_{0}
\bigg\} _{m=0}^{\infty}.
$$

\end{theorem}

\paragraph{Proof.} This is proved as part one of the Structure Theorem in Kontorovich-Sinai \cite{KS02}. Here (iii) follows immediately from (i) and (ii). $~~~\bsq$  \\

Again one easily shows that an integer $\n$ is uniquely determined by the $\ord$-sequence $V(\n)$ of its forward $U$-orbit.\\

Moreover, the following result  shows that the image under the iterated map $U^{(\kk)}$ of 
$\n\in  \gS^{(\vep)}(\ord_{1},\dots,\ord_{\kk})$ is also a nice arithmetic progression!

%
%

\begin{theorem}~\label{th52}
{\em (Iterated Structure Theorem)}
Fix $\vep=1$ or $5$, $\kk\ge1$,  a vector $(\ord_{1},\cdots,\ord_{\kk})$, and let
$\s= \ord_1+ \cdots + \ord_{\kk}$.
Suppose $1 \le \n_{0}<6\cdot 2^{\s}$ is the least element of 
$
\gS^{(\vep)}(\ord_{1},\dots,\ord_{\kk})
$.
Then there is a $\gd_{\kk}=1$ or $5$ and an $r_{\kk}\in\{0,1,2,\dots,3^{\kk}-1\}$, both depending only on $\vep$ and $(\ord_{1},\dots,\ord_{k})$,
 such that, for each  positive integer $
m$, 
\beql{521a}
U^{(\kk)}(6\cdot2^{\s}\cdot m+\n_{0})=
6(3^{\kk}\cdot m+r_{k})+\gd_{k}.
\eeq
Moreover, $\gd_{k}$ is determined by the congruence
\beql{delk}
\gd_{k}\equiv 2^{\ord_{k}}(\bmod~3).
\eeq
\end{theorem}

\paragraph{Proof.} This is part two of the Structure Theorem in Kontorovich-Sinai \cite{KS02}. 
Note that $m$ is the same number on both sides of \eqn{521a}; this equation says that  an arithmetic
progression with common difference $6 \cdot 2^{\s}$ 
mapped under $U^{(k)}$
to one with
common difference $6 \cdot 3^{k}$.
$~~~\bsq$  \\

%
%
%
\subsection{Probability Densities}

We first tweak the notion of 
natural density defined in \eqn{220b} on subsets of the natural numbers,
by restricting to just elements of  our domain $\Pi$. For a subset $\gS\subset\Pi$, let
the {\em $\Pi$-natural density} be
$$
\DD_{\Pi}(\gS):=
\lim_{t\to\infty}
\frac3t
\left|
\bigg\{
\n\in\gS:\n\le t
\bigg\}
\right|
=
\lim_{t\to\infty}
{
\left|
\bigg\{
\n\in\gS:\n\le t
\bigg\}
\right|
\over
\left|
\bigg\{
\n\in\Pi:\n\le t
\bigg\}
\right|
}
,
$$
provided that the limit exists. (The factor $3$ appears because $\Pi$ contains two residue classes modulo $6$.)\\

For a vector 
$
(\ord_{1},\dots,\ord_{\kk})
$%
, let
$$
\gS
(\ord_{1},\dots,\ord_{\kk})
:=
\gS^{(1)}
(\ord_{1},\dots,\ord_{\kk})
\
\cup
\
\gS^{(5)}
(\ord_{1},\dots,\ord_{\kk})
.
$$
Recall that a random variable $X$ is {\em geometrically distributed} with parameter $0<\rho<1$ if 
$$
\PP[X=m]=
\rho^{m-1}(1-\rho)
\qquad
\qquad
\mbox{
for $m=1,2,3,\dots$
}
$$ 

%
%
\begin{theorem}~\label{th52b}
(Geometric Distribution)

(1) The sets $ \gS(\ord_{1},\dots,\ord_{\kk})$   
have a $\Pi$-natural density given by 
\beql{503}
\DD_{\Pi}
\left(
\gS
(\ord_{1},\dots,\ord_{\kk})
\right)
=2^{-\s}=
2^{-\ord_{1}}
\cdot
2^{-\ord_{2}}
\cdots
2^{-\ord_{\kk}}.
\eeq

(2) This natural density matches the probability density of the distribution
for independent geometrically 
distributed random variables $(\pp_1, ..., \pp_k)$ with
parameter $\rho= \frac{1}{2}$, which have
$$
\mu_{\ord}
:=
\EE[\pp_{j}]=2,
\qquad
\mbox{ and }
\qquad
\gs_{\ord}
:=
Var[\pp_{j}]=2.
$$
That is,
\beql{prob}
\PP[(\pp_{1}= \ord_1,\dots,\pp_{\kk}= \ord_k)] 
=
\DD_{\Pi}
\left(
\gS
(\ord_{1},\dots,\ord_{\kk})
\right).
\eeq
\end{theorem}

\paragraph{Proof.} (1) The existence of a
natural density 
is
 automatic, since these sets are finite unions of arithmetic progressions.
For $\vep=1$ or $5$, we easily compute from Theorem \ref{th51} that
$$
\DD_{\Pi}
\left(
\gS^{(\vep)}(\ord_{1},\dots,\ord_{\kk})
\right)
=
3\cdot\frac1{6\cdot 2^{\ord_{1}+\cdots+\ord_{\kk}}}
=
\frac12\cdot
2^{-\s},
$$
and hence \eqn{503} follows.

(2) The identity \eqn{prob} 
is immediate from
independence and \eqn{503}. 
$~~~\bsq$\\

We now deduce the following result.

%
%
%
\begin{theorem}~\label{th53}
{\em (Central Limit Theorem)}
For the accelerated $3x+1$ map $U$, with symbolic iterates $(\ord_1, \ord_2, \dots)$,
the scaled ordinates 
satisfy
$$
\lim_{k\to\infty}
\DD_{\Pi}
\left[
\n
:
{
\ord_{1}(n)%
+\cdots+\ord_{k}(n)
-\mu_{\ord}k
\over
\sqrt {\gs_{\ord}k}
}
<
a
\right]
=
{1\over\sqrt{2\pi}}
\int_{-\infty}^{a}
e^{-u^{2}/2}
du
.
$$
\end{theorem}

\paragraph{Proof.} This  follows immediately from the argument above and the Central Limit Theorem
for geometrically distributed random variables.
~~~ $\bsq$  \\

Compare the above to Theorem \ref{th31}. 
The rate of convergence 
is  %
again
quite slow 
(this feature 
is shared by
Borovkov-Pfeifer; see \eqn{312b}).\\


%
%
%
\subsection{Brownian Motion}

Consider some starting value $x_{0}=\n\in\Pi$,  denote its iterates by $x_{\kk}:=U^{(\kk)}(x_{0})$, and take logarithms $y_{\kk}:=\log x_{\kk}$. 
As in \eqn{211}
, the multiplicative behavior of $U$ is converted via logarithms to an additive behavior.
Normalize the above by
\beql{504}
\gw_{k}
:=
{
y_{k} -y_{0}
-k\log (\frac34)
\over
\sqrt{2k}
\log 2
} 
.
\eeq

Then we have the following scaling limits for ``random" accelerated trajectories,
chosen in the sense of density. 

%
%

\begin{theorem}~\label{th54}
{\em (Geometric Brownian Motion Increments)}
Fix 
a partition of
 the interval $[0,1]$ as $0=t_{0}<t_{1}<\cdots<t_{r}=1$. 
 Given an integer $k$, set  $k_{j}=\lfloor t_{j}k\rfloor$, with $j=1,\dots,r$.
Then for  any $a_{j}<b_{j}$, 
$$
\DD_{\Pi}
\left[
x_{0}:
a_{j}<
{
\gw_{k_{j}}
-
\gw_{k_{j-1}}
}
<
b_{j}
,
\mbox{ for all }
j=1,2,\dots,r
\right]
\quad
\to
\quad
\prod_{j=1}^{r}
\bigg(
\Phi(b_{j})-\Phi(a_{j})
\bigg)
,
$$
 as $k\to\infty$,
where
recall that $\Phi(a)$ is the cumulative distribution function for the standard normal 
 distribution:
$$
\Phi(a)
=
{1\over\sqrt{2\pi}}
\int_{-\infty}^{a}
e^{-u^{2}/2}
du
.
$$
\end{theorem}

\paragraph{Proof.} This  appears as Theorem 5 in Kontorovich-Sinai \cite{KS02}. See Figure \ref{fig:BM}.
~~~ $\bsq$  \\

\begin{figure}
\begin{center}
\includegraphics{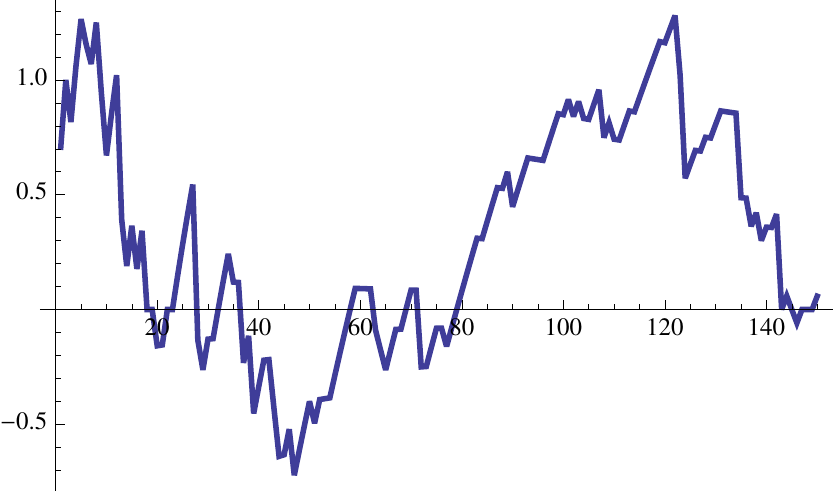}
\end{center}
\caption{A sample path of the $3x+1$ map. Here we took the starting value $x_{0}=123\,456\,789\,135\,791\,113\,151\,719$, computed $150$ iterates of $U$, and plotted $\gw_{k}$. 
}
\label{fig:BM}
\end{figure}

The interpretation of the above result is that  
the paths of the accelerated $3x+1$ map, when
properly scaled, approach  
 those of a geometric Brownian motion, that is, a stochastic process whose logarithm is a Brownian motion, or a Weiner process. 

\paragraph{Remark.}
There are in fact {\bf two} limits taken in the above theorem, whose order is highly non-interchangeable! The first limit is hidden inside the definition of density, that is, first we take the limit as $x\to\infty$ of the set of all $x_{0}<x$ satisfying the given condition with the number $k$ of iterates of $U$ fixed, and 
only
then 
do
we let $k\to\infty$. If $x_{0}$ were to be fixed and $k$ allowed to grow, there would be nothing stochastic at all happening, since we believe the $3x+1$ Conjecture!


\paragraph{Remark.}
The drift, as given in \eqn{504}, is $\log(\frac34)\approx-0.28768$. Compare this to  \eqn{312a}, where the drift of the Biased Random Walk model is computed to be $\frac12\log(\frac34)\approx-0.14384$.
While it is not surprising that the accelerated map $U$ should have a more aggressive pull to the origin, it is curious that it is exactly twice as fast (on an exponential scale) as the $3x+1$ function $T$. 

\paragraph{Remark.}
Given that the drift of the (logarithm of the) accelerated $3x+1$ function $U$ is $\mu=\log(\frac34)$,
one expects that the typical 
total stopping time 
of a seed $\n$ is roughly 
$$
{1\over |\mu|}\log \n\approx 3.476\log \n
.
$$

%
%
%
\subsection{Entropy}

%
%
%
\begin{defi}~\label{de55}
{\em 
The entropy of a random variable $X$ taking values in $[M]:=\{1,2,\dots,M\}$ is given by
$$
H:=-\sum_{m=1}^{M}\PP[X=m]\log\PP[X=m].
$$
}
\end{defi}
The following facts are classical:
\begin{enumerate}
\item[(i)] If $X$ is distributed uniformly in $[M]$ then $H=\log M$.
\item[(ii)] The entropy $H$ is maximized by the uniform distribution.
\end{enumerate}

The first is an elementary exercise, while the second is proved easily using, e.g., Lagrange's multiplier method.\\

In light of Theorem \ref{th52}, for any fixed $k\ge1$, the value $0\le r_{k}\le3^{k}-1$ is a function of the values $\vep$ and $(\ord_{1},\dots,\ord_{k})$, 
and hence has
 a natural density. 
For a fixed $\frak r
\in[0,3^{k}-1]$ we write
$$
\DD_{\Pi}
[x_{0}:r_{k}(x_{0})= 
\frak r
]\qquad\mbox{ to mean }\qquad 
\sum_{(\ord_{1},\dots,\ord_{k}), \ \vep\in\{1,5\}\atop r_{k}(\vep,\ord_{1},\dots,\ord_{k})=\frak r
}
\DD_{\Pi}[\gS^{(\vep)}(\ord_{1},\dots,\ord_{k})
]
.
$$
One might hope that $ r_{k}$ (which is a deterministic function but can be thought of as a ``random variable'') is close to being uniformly distributed in $\{0,1,\dots,3^{k}-1\}$; then one could attempt to ``bootstrap'' iterations of $U$ to one-another to have better quantitative control on various asymptotic densities with $k$ not too large. 
Were this to be the case, the entropy (defined for this using $\DD_{\Pi}$ in place of $\PP$) would be $\log 3^{k}=k\log 3$.

%
%
\begin{theorem}~\label{th55}
{\em (Entropy of $r_{k}$)}
There is some constant $c>0$ such that the entropy $H$ of $r_{k}$ satisfies:
$$
H\ge k\log 3 - c \log k.
$$
\end{theorem}

\paragraph{Proof.} This  statement is  Theorem 5.1 in Sinai \cite{Si03a}.
~~~ $\bsq$  \\


The function 
$r_{k}$ in Theorem \ref{th52} is accompanied by the residue class $\gd_{k}\in\{1,5\}$, which satisfies, cf. \eqn{delk}, 
$$
\gd_{k}\equiv 2^{\ord_{k}}(\bmod 3).
$$
It follows immediately from the fact that $\ord_{k}$ is geometrically distributed with parameter $1/2$, that
$$
\DD_{\Pi}
[x_{0}:\gd_{k}(x_{0})=1]
=
\DD_{\Pi}
[x_{0}:\ord_{k}\mbox{ is even }] = \frac13,
$$
and hence of course, $
\DD_{\Pi}
[x_{0}:\gd_{k}(x_{0})=5]
=\frac23$. \\

Moreover, if $r_{k}$ is uniformly distributed, then so are the digits $h_{k}(j)\in\{0,1,2\}$ in its $3$-adic expansion:
$$
r_{k}
=
h_{k}(k-1)\cdot 3^{k-1}
+h_{k}(k-2)\cdot 3^{k-2}
+\cdots
+ h_{k}(1) \cdot 3 
+h_{k}(0)
.
$$

Note that only the first few leading digits $h_{k}(k-1), h_{k}(k-2), \dots , h_{k}(k-t)$ are needed to specify that location of $r_{k}/3^{k}$, to within an error of $1/3^{t}$.

%
%
\begin{theorem}~\label{th56}
{\em (Joint Uniform Distribution)}
The joint distributions of $(r_{k}/3^{k},\gd_{k})$ converge weakly to the uniform distribution, that is, for 
any fixed $t\ge1$ and  $\frak h_{1},\dots,\frak h_{t}\in\{0,1,2\}$, 
as $k\to\infty$,
$$
\DD_{\Pi}
[x_{0}:h_{k}(k-1)=\frak h_{1},h_{k}(k-2)=\frak h_{2},\dots,h_{k}(k-t)=\frak h_{t},\gd_{k}(x_{0})=1]
\to
{
1\over
3^{t}
}
\cdot
{
1\over
3
}
,
$$
and
$$
\DD_{\Pi}
[x_{0}:h_{k}(k-1)=\frak h_{1},h_{k}(k-2)=\frak h_{2},\dots,h_{k}(k-t)=\frak h_{t},\gd_{k}(x_{0})=5]
\to
{
1\over
3^{t}
}
\cdot
{
2
\over
3
}
.
$$
\end{theorem}

\paragraph{Proof.} This appears as Theorem 1 in Sinai \cite{Si03b}. See also \cite{Si04}.
~~~ $\bsq$  \\









%
%
%
\section{$3x+1$ Backwards Iteration: $3x+1$ Trees}\label{sec6}
\hsp

One can also model  backwards iteration of the $3x+1$ map  $T(x)$.\\

Backwards iteration 
 is described by a tree of inverse iterates, and 
there  are either one or two inverse iterates.
 Here
 $$
 T^{(-1)}(n) =
 \left\{
\begin{array}{cl}
\{ 2n\} & \mbox{if}~ n \equiv 0, 1 ~(\bmod~3),  \\
~~~ \\
\{ 2n, \frac{2n-1}{3} \}& \mbox{if} ~~n \equiv 2  (\bmod~3).
\end{array}
\right.
$$
Starting from a root node labelled $a$ we can grow an infinite tree
$\sT(a)$ of all the inverse iterates of $a$.  Each node in the tree
is labelled by its associated $3x+1$ function value. To 
a node labelled $n$ we add  either one or two (directed) edges from
the elements of $T^{-1}(n)$ to $n$, and we label these two edges by 
the value of this element. \\

%
%
%
\subsection{Pruned $3x+1$ Trees}

Next we note  that any $a \equiv 0 ~(\bmod~3)$ has exactly one inverse iterate,
which itself is $0~(\bmod~3)$. Thus if $a \equiv 0~(\bmod~3)$
the  set of inverse iterates forms a single 
branch that never divides. However if $a \not\equiv 0~(\bmod~3)$
then the tree grows exponentially in size. 
It is convenient therefore to restrict to numbers $a \not\equiv 0~(\bmod~3)$
and furthermore to prune such a tree to remove all nodes
$n \equiv 0~(\bmod~3)$. This produces  an (infinite depth)  {\em pruned tree} ${\sT}^{\ast}(a)$ 
which is described by inverse iterates of the modified map
 \beql{602}
 \tilde{T}^{(-1)}(n) =
 \left\{
\begin{array}{cl}
\{ 2n\} & \mbox{if}~ n \equiv 1, 4, 5 ~\mbox{or} ~7 ~(\bmod~9),  \\
~~~ \\
\{ 2n, \frac{2n-1}{3} \}& \mbox{if} ~~n \equiv 2  ~\mbox{or}~8 ~ (\bmod~9),
\end{array}
\right.
\eeq
applied starting with root node labelled $n_0:=
a
$. The pruning operation
is depicted in Figure ~\ref{fig1}, with root node assigned depth $0$.

%
%
%

\begin{figure}[hbp]\centering
\setlength{\tabcolsep}{0.125in}
\vspace{.1in}
\begin{tabular}{ccc}
\setlength{\unitlength}{0.012in}
\begin{picture}(153,172)(57,436)
\thicklines
\put( 70,600){\makebox(0,0)[rb]{16}}
\put(210,600){\makebox(0,0)[lb]{1}}
\put(170,520){\makebox(0,0)[lb]{1}}
\put(170,560){\makebox(0,0)[lb]{2}}
\put(170,600){\makebox(0,0)[lb]{4}}
\put(160,560){\line( 1, 1){ 40}}
\put(160,560){\line( 0, 1){ 40}}
\put( 80,480){\line( 2, 1){ 80}}
\put(160,520){\line( 0, 1){ 40}}
\put( 70,520){\makebox(0,0)[rb]{4}}
\put(130,600){\makebox(0,0)[lb]{5}}
\put( 70,560){\makebox(0,0)[rb]{8}}
\put( 70,480){\makebox(0,0)[rb]{2}}
\put( 70,440){\makebox(0,0)[rb]{1}}
\put( 80,440){\circle*{6}}
\put(120,600){\circle*{6}}
\put( 80,600){\circle*{6}}
\put( 80,560){\circle*{6}}
\put( 80,520){\circle*{6}}
\put( 80,480){\circle*{6}}
\put( 80,560){\line( 1, 1){ 40}}
\put(160,520){\circle*{6}}
\put( 80,560){\line( 0, 1){ 40}}
\put( 80,520){\line( 0, 1){ 40}}
\put( 80,480){\line( 0, 1){ 40}}
\put( 80,440){\line( 0, 1){ 40}}
\put(200,600){\circle*{6}}
\put(160,600){\circle*{6}}
\put(160,560){\circle*{6}}
\end{picture}
&
\setlength{\unitlength}{0.012in}
\begin{picture}(113,172)(57,436)
\thicklines
\put(120,520){\line( 1, 1){ 40}}
\thinlines
\put(167,593){\line( 0, 1){ 14}}
\put(167,607){\line(-1, 0){ 14}}
\put(153,607){\line( 0,-1){ 14}}
\put(153,593){\line( 1, 0){ 14}}
\put(167,553){\line( 0, 1){ 14}}
\put(167,567){\line(-1, 0){ 14}}
\put(153,567){\line( 0,-1){ 14}}
\put(153,553){\line( 1, 0){ 14}}
\thicklines
\put(120,560){\line( 0, 1){ 40}}
\put( 80,480){\line( 1, 1){ 40}}
\put( 70,600){\makebox(0,0)[rb]{64}}
\put(130,510){\makebox(0,0)[lb]{5}}
\put(110,560){\makebox(0,0)[rb]{10}}
\put(110,600){\makebox(0,0)[rb]{20}}
\put(170,600){\makebox(0,0)[lb]{6}}
\put(170,560){\makebox(0,0)[lb]{3}}
\put(120,520){\line( 0, 1){ 40}}
\thinlines
\put(120,520){\circle{14}}
\thicklines
\put(120,600){\circle*{6}}
\put( 70,560){\makebox(0,0)[rb]{32}}
\put( 70,520){\makebox(0,0)[rb]{16}}
\put( 70,480){\makebox(0,0)[rb]{8}}
\put( 80,440){\circle*{6}}
\put( 80,480){\circle*{6}}
\put( 80,520){\circle*{6}}
\put( 80,560){\circle*{6}}
\put( 80,600){\circle*{6}}
\put( 70,440){\makebox(0,0)[rb]{4}}
\put( 80,440){\line( 0, 1){ 40}}
\put( 80,560){\line( 0, 1){ 40}}
\put( 80,520){\line( 0, 1){ 40}}
\put( 80,480){\line( 0, 1){ 40}}
\put(160,560){\circle*{6}}
\put(160,560){\line( 0, 1){ 40}}
\put(160,600){\circle*{6}}
\put(120,560){\circle*{6}}
\put(120,520){\circle*{6}}
\end{picture}
&
\setlength{\unitlength}{0.012in}
\begin{picture}(73,172)(57,436)
\thinlines
\put(120,520){\circle{14}}
\thicklines
\put(120,560){\line( 0, 1){ 40}}
\put(120,520){\line( 0, 1){ 40}}
\put( 80,480){\line( 1, 1){ 40}}
\put(130,510){\makebox(0,0)[lb]{5}}
\put(110,600){\makebox(0,0)[rb]{20}}
\put(110,560){\makebox(0,0)[rb]{10}}
\put( 70,600){\makebox(0,0)[rb]{64}}
\put( 70,480){\makebox(0,0)[rb]{8}}
\put( 70,520){\makebox(0,0)[rb]{16}}
\put( 70,440){\makebox(0,0)[rb]{4}}
\put( 70,560){\makebox(0,0)[rb]{32}}
\put( 80,440){\circle*{6}}
\put( 80,600){\circle*{6}}
\put( 80,560){\circle*{6}}
\put( 80,520){\circle*{6}}
\put( 80,480){\circle*{6}}
\put( 80,560){\line( 0, 1){ 40}}
\put(120,600){\circle*{6}}
\put( 80,520){\line( 0, 1){ 40}}
\put( 80,480){\line( 0, 1){ 40}}
\put( 80,440){\line( 0, 1){ 40}}
\put(120,520){\circle*{6}}
\put(120,560){\circle*{6}}
\end{picture}
\\
(i) $\sT_4(1)$ & (ii) $\sT_4(4)$ & (iii) $\sT_{4}^{*}(4)$
\end{tabular}
\caption{$3x+1$ trees $\sT_k(
a
)$ and ``pruned'' $3x+1$ tree $\sT_k^{*}(
a
)$, with $k=4$.}
\label{fig1}
\end{figure}
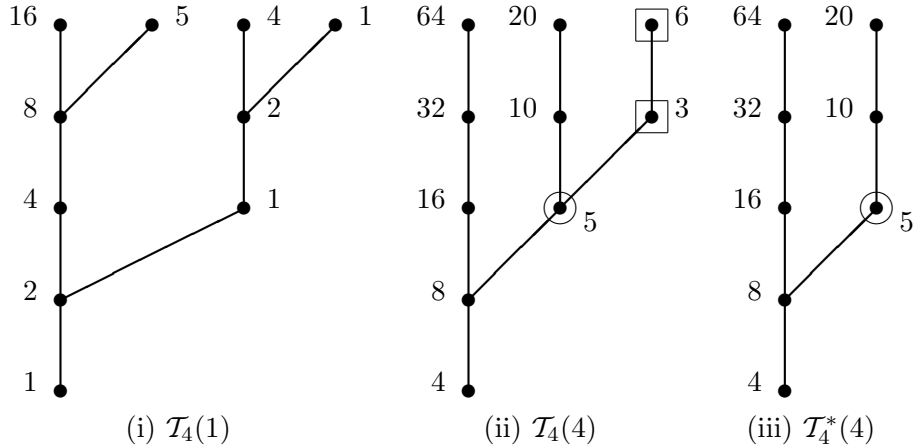

We obtain a reduced tree ${\overline{\sT}}^{\ast}(a)$ obtained by
labelling each node with the $(\bmod~2)$ residue class of the $3x+1$
value assigned to that node. (One may also think of this as labelling the
directed edge leaving this node, with the exception of the root node.)\\

We let ${\sT}_k^{\ast}(a)$ denote the pruned tree with root node $n_0=a$,
cut off at depth $k$, and we let ${\overline{\sT}}_k^{\ast}(a)$ denote the
same tree, keeping only the node labels $(\bmod~2),$ for all nodes
except the root node, where no data is kept.
Let $N^{\ast}(k; a)$ count the number of depth $k$
leaves in this tree.  Then we have
\beql{610a}
N^{\ast}(k, a) := |\{ n: ~n \not\equiv 0~(\bmod~3) ~\mbox{and} ~~T^{(k)}(n)=a \}|
.
\eeq
We have $N^{\ast}(k, a) \le 2^k$ as a consequence of the fact that each $3x+1$ tree
has at most two upward branches at each node.\\

The following result gives information on the sizes of depth $k$
trees over all possible tree types (\cite[Theorem 3.1]{LW92}).\\
%
%
%

\begin{theorem}~\label{th61} {\em (Structure of Pruned $3x+1$ Trees)}

(1) For $k \ge 1$ and $
a
\not\equiv 0~ (\bmod~3)$, the
 structure of the pruned level $k$ tree ${\overline{\sT}}_k^{\ast}(a)$, and
 hence the number $N^{\ast}(k, a)$,  is completely 
determined by $a ~(\bmod~3^{k+1})$. \\

(2) There are $2 \cdot 3^k$ residue classes $
a
~(\bmod ~3^{k+1})$ with
$
a
\not\equiv 0 ~(\bmod~3)$. For these 
\beql{611}
\sum_{{a~(\bmod~3^{k+1})}\atop{m \not\equiv 0~(\bmod~3)}} N^{\ast}(k, a)= 2 \cdot 4^k.
\eeq 
It follows that if a residue class $a~(\bmod~3^{k+1})$ with $a \not\equiv 0 (\bmod~3)$ is
picked with the uniform distribution, the expected number of leaves in 
the random tree  ${\overline{\sT}}_k^{\ast}(a)$ is exactly $\left(\frac{4}{3}\right)^k$.
\end{theorem}


We now consider the complete set of numbers having total stopping time $k$.
Set
\beql{620}
N_k := |\{ n :  ~\sigma_{\infty}(n)=k\}|.
\eeq
Recall from 
\S\ref{sec2p5} that $N_k= N_k(1)$, where $N_k(a)$ counts the number of integers
that iterate to $a$ after exactly $k$ iterations of the $3x+1$ map $T$. 
We defined there the $3x+1$ tree growth constants
$$
\delta_3(a) :=\limsup_{k \to \infty} \frac{1}{k}\log N_k(a).
$$
Theorem~\ref{th61}
suggests the following conjecture for these tree growth constants,
 made by Lagarias and Weiss \cite{LW92}.\\

%
%
%

\begin{conj} \label{conj61}
 For each $a \not \equiv 0 ~(\bmod~3)$, the $3x+1$ tree growth constant  $\delta_3(a)$ is given by 
\beql{714}
\delta_3(a) = \log \left(\frac{4}{3}\right).
\eeq
\end{conj}

Applegate and Lagarias \cite{AL95a} determined
by computer the maximal and minimal number of leaves in $3x+1$ trees of depth $k$ for $k \ge 30$.
The maximal and minimal number of leaves in such trees at level $k$ is given by
$$
N_k^{+} := \max \{  N_k^{\ast}( a):  ~a~(\bmod~3^{k+1}) ~\mbox{with}~ a \not\equiv 0~(\bmod~3) \}.
$$
and
$$
N_k^{-} = \min\{ N_k^{\ast}( a):  ~a~(\bmod~3^{k+1}) ~\mbox{with} ~a \not\equiv 0~(\bmod~3) \},
$$
respectively. Figure \ref{fig2} pictures maximal and minimal  trees for depth $k=5$. 
(Circled nodes indicate
an omitted inverse iterate under  $T^{-1}$ that is $\equiv 0~(\bmod~3)$.)\\

%
%
%

\begin{figure}[htbp]\centering
\setlength{\tabcolsep}{0.125in}
\vspace{.1in}
\begin{tabular}{cc}
\setlength{\unitlength}{0.0085in}
\begin{picture}(79,212)(51,436)
\thinlines
\put( 80,480){\circle{14}}
\put( 70,600){\makebox(0,0)[r]{112}}
\put(130,640){\makebox(0,0)[l]{74}}
\put( 70,640){\makebox(0,0)[r]{224}}
\thicklines
\put( 80,600){\line( 0, 1){ 40}}
\put(120,600){\line( 0, 1){ 40}}
\put( 80,560){\line( 1, 1){ 40}}
\put( 70,480){\makebox(0,0)[r]{14}}
\put(130,600){\makebox(0,0)[l]{37}}
\put( 70,560){\makebox(0,0)[r]{56}}
\put( 70,520){\makebox(0,0)[r]{28}}
\put( 70,440){\makebox(0,0)[r]{7}}
\put( 80,440){\circle*{6}}
\put( 80,600){\circle*{6}}
\put( 80,560){\circle*{6}}
\put( 80,520){\circle*{6}}
\put( 80,480){\circle*{6}}
\put( 80,560){\line( 0, 1){ 40}}
\put( 80,640){\circle*{6}}
\put( 80,520){\line( 0, 1){ 40}}
\put( 80,480){\line( 0, 1){ 40}}
\put( 80,440){\line( 0, 1){ 40}}
\put(120,640){\circle*{6}}
\put(120,600){\circle*{6}}
\end{picture}
&
\setlength{\unitlength}{0.0085in}
\begin{picture}(313,212)(217,436)
\thinlines
\put(320,640){\circle{14}}
\put(400,600){\circle{14}}
\put(440,640){\circle{14}}
\put(360,600){\circle{14}}
\thicklines
\put(400,480){\line( 0, 1){ 40}}
\put(360,440){\line( 1, 1){ 40}}
\put(400,560){\line( 0, 1){ 40}}
\put(400,520){\line( 0, 1){ 40}}
\put(400,520){\line( 1, 1){ 40}}
\put(280,600){\line( 0, 1){ 40}}
\put(320,560){\line( 0, 1){ 40}}
\put(320,600){\line( 0, 1){ 40}}
\put(360,520){\line( 0, 1){ 40}}
\put(360,560){\line( 0, 1){ 40}}
\put(360,600){\line( 0, 1){ 40}}
\put(400,600){\line( 0, 1){ 40}}
\put(270,600){\makebox(0,0)[r]{35}}
\put(230,640){\makebox(0,0)[r]{23}}
\put(530,640){\makebox(0,0)[l]{7}}
\put(490,600){\makebox(0,0)[l]{11}}
\put(450,560){\makebox(0,0)[l]{17}}
\put(410,520){\makebox(0,0)[l]{26}}
\put(410,480){\makebox(0,0)[l]{13}}
\put(310,560){\makebox(0,0)[r]{53}}
\put(440,560){\line( 0, 1){ 40}}
\put(440,600){\line( 0, 1){ 40}}
\put(440,560){\line( 1, 1){ 40}}
\put(480,600){\line( 0, 1){ 40}}
\put(350,520){\makebox(0,0)[r]{80}}
\put(350,480){\makebox(0,0)[r]{40}}
\put(350,440){\makebox(0,0)[r]{20}}
\put(480,600){\line( 1, 1){ 40}}
\put(320,560){\line(-1, 1){ 40}}
\put(280,640){\circle*{6}}
\put(320,600){\circle*{6}}
\put(320,640){\circle*{6}}
\put(360,560){\circle*{6}}
\put(360,600){\circle*{6}}
\put(360,640){\circle*{6}}
\put(240,640){\circle*{6}}
\put(360,440){\circle*{6}}
\put(360,480){\circle*{6}}
\put(400,480){\circle*{6}}
\put(360,520){\circle*{6}}
\put(320,560){\circle*{6}}
\put(280,600){\circle*{6}}
\put(400,520){\circle*{6}}
\put(440,600){\circle*{6}}
\put(440,640){\circle*{6}}
\put(360,440){\line( 0, 1){ 40}}
\put(360,480){\line( 0, 1){ 40}}
\put(360,520){\line(-1, 1){ 40}}
\put(280,600){\line(-1, 1){ 40}}
\put(480,640){\circle*{6}}
\put(400,560){\circle*{6}}
\put(400,600){\circle*{6}}
\put(400,640){\circle*{6}}
\put(440,560){\circle*{6}}
\put(480,600){\circle*{6}}
\put(520,640){\circle*{6}}
\end{picture}
\\
(i) $\sT_5^*(7)$ attains $N^-(5) = 2$ & (ii) $\sT_5^*(20)$ attains $N^+(5) = 8$
\end{tabular}
\caption{Maximal and Minimal depth 5 pruned $3x+1$ Trees}
\label{fig2}
\end{figure}
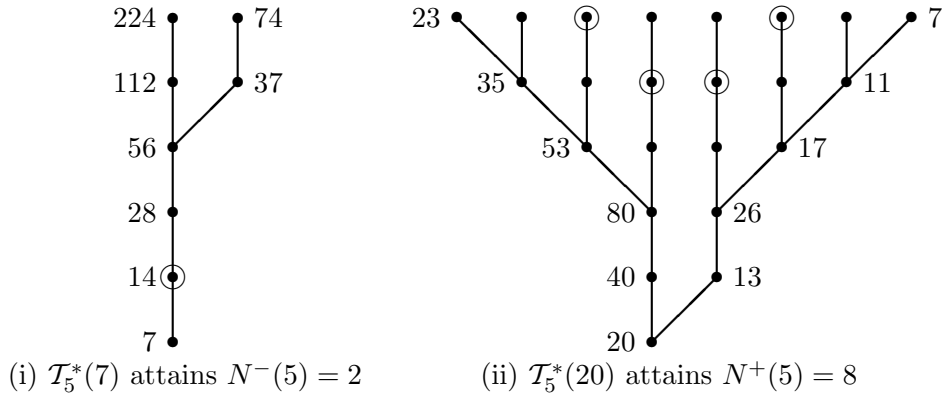

The data on these counts $N_{\pm}(k)$ was presented already in \S\ref{sec2p5}, cf. Table~\ref{tab23}.
Based on this data, Applegate and Lagarias \cite[Conjecture C]{AL95c}
 formulated the following  strengthened conjecture, which implies Conjecture \ref{conj61}.

%
%
%

\begin{conj}~\label{conj62}
 The maximal and minimal number of leaves of $3x+1$ trees satisfy, as $k \to \infty$, 
\beql{641}
N_k^{-} = \left(\frac{3}{4}\right)^{k+ o(k)}
\eeq
\beql{642}
N_k^{+} = \left(\frac{3}{4}\right)^{k+ o(k)}.
\eeq
\end{conj}

%
%
%
\subsection{$3x+1$ Backwards Stochastic Models:  Branching Random Walks}

Lagarias and Weiss \cite{LW92} formulated
stochastic models for the growth of $3x+1$  trees that were 
 {\em multi-type branching processes}. Such models grow  a random
tree, with nodes marked as several different kinds of individuals.
In this case the number of nodes of each type at each depth $k$ 
(also called generation $k$) can be viewed as
the output of the branching process. The particular branching processes they
used are {\em multi-type Galton-Watson processes}. \\

Lagarias and Weiss also 
modeled 
the size of preimages
of elements in  a (pruned) $3x+1$ tree. This size 
is
specified
by a real number attached to each node. 
Branching process models which  attach to each node in the tree
a real number giving the 
position of those individuals on a line, according to some (possibly random) rule,
are models called {\em multi-type branching random walks}. Here the location
of the individuals on the line give the random walk aspect; offspring nodes at level $k$ are shifted
in position from their parent ancestor at level $k-1$ by a point process. The
process starts with a root node giving a single progenitor at level $0$ (generation $0$).\\


Lagarias and Weiss defined a hierarchy of 
branching random walk models, which they denoted $\sB[3^j]$, for
each $j \ge 0$. 
These branching random walk models, having several kinds of
individuals,  model the backwards iteration viewed $~(\bmod ~ 3^j)$.
The model for $j=0$ is simpler than the other models.\\

 {\em $3x+1$ Branching Random walk ~$\sB[3^0]$.} There is one type of
 individual. With probability $\frac{2}{3}$ an individual has a single offspring
 located at a position shifted by $\log 2$ on the line from its progenitor, and with probability $\frac{1}{3}$
 it has two offspring located at positions shifted $\log 2$ and $\log \frac{2}{3}$ on
 the line from their progenitor. If the progenitor is in generation (or depth) $k-1$, the offspring are
 in generation $k$. The tree is grown from a single individual at generation $0$,
 with specified location $\log a$.\\

The more general models  for $j \ge 1$ are given as follows.\\

 {\em $3x+1$ Branching Random walk ~$\sB[3^j], (j \ge 1)$.} There are $p=2 \cdot 3^{j-1}$
 types of individuals, indexed by residue classes $
 a
 ~(\bmod~3^j)$ with 
 $
 a
 \not\equiv 0~(\bmod~3)$. The distribution of offspring of an individual of
 type $
 a
 ~(\bmod~3^j)$, at any given depth $k$ in
 the branching, is determined as follows: Regard $
 a
 ~(\bmod~3^{j})$
 labelling a node at depth $d-1$. Regard it 
 as being, with probability $\frac{1}{3}$ each, one of the three possible residue
 classes $\tilde{
a 
}~(\bmod~3^{j+1})$ consistent with it. The tree (of depth $1$)
 with $\tilde{
 a
 }$ as root node, 
  given by $(T^{\ast})^{-1}(\tilde{
  a
  })$ has either one or
 two progeny, at depth $1$ and their node labels are well-defined classes $(\bmod~3^j)$,
 either $2 \tilde{
 a
 }$ or, if it legally occurs, $\frac{2 \tilde{
 a
 }-1}{3} (\bmod~3^j)$. 
 The branching random walk then
 produces an individual of type $2\tilde{
 a
 }$ 
 at generation $k+1$ whose position  is additively shifted by $\log 2$ from that of the
 generation $k$ progenitor node, 
 plus, if legal, another labelled
   $\frac{2 \tilde{
   a
   }-1}{3}$, which is shifted in position by
    $\log(\frac{2}{3})$ on the line from that of the generation $k$-node. 
    The tree is grown from a single individual at depth $0$, with
    specified type and location $\log a$.  \\

In these models, the behavior of the random walk part of the model can be
completely reconstructed from knowing  the type of each node. 
This is a very special property of these branching random walk models, which does not hold for
general branching random walks.\\
 
 In such models, one may think of the nodes as representing individuals,
with individuals at level $k$ being children of a particular individual at level
$k-1$; the random walk aspect indicates position in space of these individuals. \\
 
 Let $\omega$ denote a single realization of such a branching random walk 
 $\sB[3^j]$  which starts from a single individual $\omega_{0,1}$ of type $1~(\bmod~3^j)$
 at depth $0$, with initial position labeled $\log 
 a
 $. Here $\omega$ describes
 a particular infinite tree. We let $N_k(\omega)$
 denote the number of individuals at level $k$ of the tree. We let
 $S( \omega_k, j)$
 denote the position of the $j$-th individual at level $k$ in the tree, for
 $1 \le j \le N_k(\omega)$.  \\
 
 These models are all {\em supercritical branching processes}
in the following very strong sense. In {\em every} random realization $\omega$,
the number of
nodes at level $d$  grows exponentially in  $d$, and there are no 
extinction events.\\

 Lagarias and Weiss \cite{LW92}
  observed that the predictions of these models stabilized for all $j \ge 1$, 
 as far as the behavior of 
 asymptotic statistics related to $3x+1$ trees is concerned. 
 This is illustrated in the following theorems.

%
%
%
\subsection{$ 3x+1$ Backwards Model Prediction: Tree Sizes }

 Concerning the number of nodes $N_k(\omega)$ in a realized tree at depth $k$,
Lagarias and Weiss  proved the following result \cite[Corollary 3.1]{LW92}. \\

%
%
%
\begin{theorem}~\label{th71}
{\em ($3x+1$ Stochastic Tree Size)} 
For all $j \ge 0$, 
a realization $\omega$ of a tree grown in
the $3x+1$ branching random walk model $\sB[3^j]$ 
has
\beql{711}
\lim_{k \to \infty} \frac{1}{k} \left( \log N_k(\omega)\right) = \log \left(\frac{4}{3}\right), ~~~~\mbox{
for almost every $\gw$.}
\eeq
\end{theorem}

This result only uses the Galton-Watson structure built into the process $\sB[3^j].$
Its prediction is consistent with the rigourous results on average tree size for pruned $3x+1$ trees
given in Theorem~\ref{th61},
and 
it 
also 
 supports 
 Conjecture \ref{conj61}.\\

%
%
%
\subsection{$3x+1$ Backwards Model Prediction: Extremal Total Stopping Times}

Next, as a statistic that corresponds to an extremal trajectory, consider the
{\em first birth in generation $k$}, which is the leftmost individual on the line at
depth $k$ in the branching random walk. Denote the location of this individual
by $L_k^{\ast}(\omega)$, for a given realization $\omega$ of the random walk.
Lagarias and Weiss  \cite[Theorem 3.4]{LW92} proved the following result.

%
%
%
\begin{theorem}~\label{th72} 
{\em (Asymptotic First Birth Location)}
For any $3x+1$ branching random walk model $\sB[3^j]$ with $j \ge 2$, 
 there is a constant $\beta_{BP}$ such that 
 for all $j \ge 
0
 $,
 the branching random walk $\sB[3^j]$ has asymptotic first birth (leftmost birth)
 \beql{721}
 \lim_{k \to \infty} L_k^{\ast}(\omega) = \beta_{BP} ~~~~~\mbox{
for almost every $\gw$.
 }
 \eeq
 This constant $\beta_{BP} \approx 0.02399$ is determined uniquely
 by the properties that it is the unique 
$\gb>0$
 that
 satisfies
 \beql{722}
 \tilde{g}(\beta)=0
 \eeq
 where
 \beql{723}
 \tilde{g}(a) := - sup_{\theta \le 0} \left( a \theta - \log  M_{BP}(\theta) \right) 
 .
 \eeq
 Here
  $M_{BP}(\theta)$ is the branching process moment generating function
 \beql{723a}
 M_{BP}(\theta) := 2^{\theta} + \frac{1}{3}(\frac{2}{3})^{\theta}.
 \eeq
 \end{theorem}
 
 Since the first birth individual at depth $k$ corresponds to taking $k$ iterations to
 reach the root node, we can define
 a  {\em branching process scaled stopping limit} $\gamma_{BP}(\omega)$. 
 This is
 the BP model's prediction for the scaled stopping constant $\g$ from \eqn{404},
defined 
  by
 $$
 \gamma_{BP}(\omega)  := \limsup_{k \to \infty} \frac{k}{L_k^{\ast} (\omega)}   .
$$
Theorem \ref{th72} implies that this value is constant (almost surely independent of $\gw$), and takes 
 the value
\beql{730f}
\gamma_{BP} = (\beta_{BP})^{-1}.
\eeq
Note that since $\gb_{BP}\approx 0.02399$, we have $1/\gb_{BP}\approx41.7$. 
At this point we have two completely different predictions for the scaled stopping constant $\g$, one from the RRW model (cf. Theorem \ref{th41}) which 
approximates forward iterations, and another from the BP models which estimate backwards iterations. 
Applegate and Lagarias then prove \cite[Theorem 4.1]{LW92}  the following 
striking 
identity.

%
%
%
\begin{theorem}~\label{th73} 
{\em ($3x+1$ Random Walk-Branching Random Walk Duality)}
The $3x+1$ repeated random walk (RRW) stochastic model scaled stopping time limit 
$\gamma_{\RW}$ and the $3x+1$ branching random walk (BP)
model $\sB[3^j]$ with $j 
\ge
0$, scaled stopping time limit $\gamma_{BP}$ are
identical! I.e.,
\beql{731aa} 
\gamma_{\RW}= \gamma_{BP}.
\eeq
\end{theorem}

\paragraph{Proof.} 
This is a consequence  \cite{LW92}  of an identity
relating the moment generating functions associated to the two models, which is
$M_{BP}(\theta) = M_{\RW}(\theta +1)$; compare \eqn{423a} and 
\eqn{723a}.
$~\bsq$\\

\paragraph{Remark.} 
Recall 
the critique of the RRW model given in 
\S\ref{sec4p4},
that various trajectories coalesce in their forward iterates.
But the BP
 models, by their tree construction, 
completely take into account the dependence
caused by coalescing trajectories! 
Since both models predict the same exact value for $\g$, it appears the critique has been thwarted off.

%
%
%
\subsection{$3x+1$ Backwards Model Prediction: Total Preimage Counts}\label{sec6p5}

We next consider what the branching process models have to say
about the number of integers below $x$ that eventually iterate to
a given integer $a$. 

The following result gives, for the simplest branching random walk model, 
 an almost sure asymptotic  of
 the number of inverse
iterates of size below a given bound (\cite[Theorem 4.2]{LW92}).

%
%
%
\begin{theorem}~\label{th75} 
{\em (Stochastic Inverse Iterate Counts)}
For a realization $\omega$ of the branching random walk $\sB[1]$, let
$I^{\ast}(t; \omega)$ count the number of progeny located at positions  $S(\omega_{k,j}) \le x$, i.e.
\beql{681}
I^{\ast}(x; \omega) := \# \{ \omega_{k,j}:  S(\omega_{k,j}) \le x, \mbox{for~any}~~k\ge 1, ~ 1 \le j \le 
N_k(\omega)\}.
\eeq
Then the asymptotic estimate
\beql{682}
I^{\ast}(x; \omega) = x^{1+o(1)}~~~\mbox{as}~~ x \to \infty
\eeq
holds almost surely.
\end{theorem}


The model  statistic $I^{\ast}(x; \omega)$ functions as  a proxy for the function $\pi_{a}(x)$,
where $\log 
a 
$ gives the position of the root node of the branching random walk. 
This result is the stochastic analogue of Conjecture~\ref{conj21} about
the $3x+1$ growth exponent.

%
%
%

 \section{The $5x+1$ Function: Symbolic Dynamics and Orbit Statistics}\label{sec7}

 We now turn for comparison to the $5x+1$ iteration.
Some features of the dynamics of this iteration are similar to that
of the $3x+1$ problem, and some are different. 
 Here the dynamics of iteration in the long run are expected to be quite
 different globally from the $3x+1$ problem; most trajectories
 are expected to diverge.  In this section we formulate several orbit
 statistics for this map, some the same as for the
 $3x+1$ map, and some changed. We review basic results on them.  \\
%
%
%

 \subsection{$5x+1$ Forward Iteration: Symbolic Dynamics}
 
 The basic features of the $5x+1$ problem  are similar to the
 $3x+1$ problem. We introduce 
 the
 parity sequence
 \beql{801}
 S_5(n):= ( \n ~(\bmod~2), T_5(\n)~(\bmod~2), T_5^{(2)}(\n) ~(\bmod~2), ...) .
\eeq
 The symbolic dynamics is similar to the $3x+1$ map:
 all finite initial symbol sequences of length $k$ occur, 
 each one for a single residue class $
 (\bmod~2^k)$.

%
%

\begin{theorem}~\label{th800}
{\em ($5x+1$ Parity Sequence Symbolic Dynamics)}
The $k$-truncated parity sequence $S_5^{[k]}(n)$ of the first $k$ iterates of the $5x+1$ map $T(x)$
is periodic  in $n$ with period $2^k$. Each of the $2^k$ possible $0-1$ vectors occurs
exactly once in the initial segment $1 \le n \le 2^k$.
\end{theorem}

\paragraph{Proof.} The proof of 
this  result exactly parallels that of Theorem~\ref{th21}.$~~~\bsq$\\

 As before,
 the parity
sequence of an orbit of $x_0$ uniquely determines $x_0$.\\

Analysis of this recursion, assuming even and odd iterates are equally likely,
as prescribed by Theorem~\ref{th800}, we find  the logarithms of iterates grow in size
on the average. 
%
%
%

 \subsection{$5x+1$ Forward Iteration: $\lambda^{+}$-
Stopping Times}

Most $5x+1$ iteration sequences  grow on average, rather than
shrinking on average.  An appropriate notion of stopping time for this situation is
as follows. 
 
%
%
%

\begin{defi}~\label{de71a}
{\em
For fixed $\lambda \ge 1$,   the {\em $\lambda^{+}$-stopping time} 
$\sigma_{\lambda}^{+}(n)$ of a map $T_{5}: \ZZ \to \ZZ$  for input $n$ is  the
minimal value of $k \ge 0$ such that
$T_{5}^{(k)}{n} > \lambda n$, e.g.
\beql{731a}
\sigma_{\lambda}^{+}(n) := \inf \{k \ge 0:  \frac{T_{5}^{(k)}(n)}{n}  > \lambda\}.
\eeq
If no such value $k$ exists, we set  $\sigma_{\lambda}^{+}(n) = + \infty.$
}
\end{defi}

One now has the following result, which parallels Theorem~\ref{th22}
for the $3x+1$ map, except that here iterates  grow in size rather than shrink in size. 
 
%
%
%
\begin{theorem}~\label{th72a} 
{\em ($\lambda^{+}$-Stopping Time Natural  Density)}

(i) For the $5x+1$ map $T_{5}(n)$, and fixed $ \lambda \ge 1$
and $k\ge1$
,  the set 
$ S_{\lambda}^{+}(k)$
of 
integers having
$\lambda^{+}$-stopping time at most $k$ 
has a well-defined natural density 
$\DD( S_{\lambda}^{+}(k))$. \\

(ii) This natural density satisfies 
\beql{721a}
\lim_{k \to \infty} \DD (S_{\lambda}^{+}(k)) = 1.
\eeq
In particular, the set of numbers with finite $\lambda^{+}$-stopping time has
natural density $1$.
\end{theorem}

\paragraph{Proof.} 
Claim
(i) follows using the Parity Sequence Theorem~\ref{th800}. Here the set is a finite union
of arithmetic progressions $(\bmod ~2^k)$, except a finite number of initial elements
may be omitted from each such progression.

The result (ii) can be established by a similar argument to that used for
the $3x+1$ problem  in Theorem~\ref{th22}. 
$~~~\bsq$\\

Here we note a surprise: there are infinitely many exceptional integers $n$ that have
$\lambda^{+}$-stopping time equal to $+\infty$! This occurs because the
$5x+1$ problem has a periodic orbit $\{ 1, 3, 8, 4, 2\}$, and infinitely 
many positive 
seeds $n_{0}$
eventually enter this orbit, e.g. $n 
_{0}
= \frac{2^{4k} -1}{5}$
for 
any
$k \ge 2$. All of these integers have $\sigma_{\lambda}^{+} 
(n_{0})
= + \infty.$ 
Nevertheless Theorem~\ref{th72a} asserts 
that such
integers have natural density zero.\\

%
%
%

 \subsection{$5x+1$ Stopping Time Statistics: Total Stopping Times}
 
The $5x+1$ problem has a finite orbit containing $1$, and
we may define total stopping time as for the $3x+1$ function. 

%
%
%
 
 \begin{defi}~\label{de731}
 {\em 
For $n \ge 1$ the {\em total stopping time} $\sigma_{\infty}(n; T_5)$ of the
$5x+1$ function is given by 
\beql{731c}
\sigma_{\infty}(n; T_5)  := \inf \{k \ge 1:~ T_5^{(k)} (n) =1\}.
\eeq
We set $\sigma_{\infty}(n; T_5)= +\infty$ if no finite $k$ has this property.
}
\end{defi}

 Here we  expect that the vast majority of positive 
$n$ will belong to divergent trajectories, and only  a small minority of $n$ have
a well-defined finite value $\sigma_{\infty}(n; T_5)< \infty$.
It is an open problem to 
prove that {\it even a single trajectory} (such as that emanating from the starting seed $n_{0}=7$)
is divergent!\\

The best we can currently show unconditionally is 
a lower bound on the size of the extremal total stopping 
time  that grows proportionally to $\log n$.

%

\begin{theorem}\label{th731a}
{\em (Lower Bound for $5x+1$ Total Stopping Times)}
There are infinitely many $n$ whose total stopping time satisfies
\beql{739a}
\sigma_{\infty}(n, T_5) \ge  \left( \frac{\log 2 + \log 5}{(\log 2)^2}\right)\log n \approx 4.79253\log n.
\eeq
\end{theorem}

\paragraph{Proof.}
The Parity Sequence
Theorem~\ref{th71} implies there is at least one odd number $n_k$
with $1 \le n_k < 2^k$ whose
first $k-1$ iterates are also odd, so that $T_5^{(k)}(n_k) \ge (\frac{5}{2})^k n_k.$ Since 
a single step can divide by at most  $2$, we necessarily have
(using 
$\log n_k \le k \log 2$), 
$$
\frac{\sigma_{\infty}(n_k, T_5)}{\log n_k} \ge  \frac{k}{\log n_k} +
\left(\frac{k \log\frac{5}{2} +\log n_k}{\log 2}\right)\frac{1}{\log n_k} 
\ge \frac{2}{\log 2} + \frac{ \log \frac{5}{2}} {(\log 2)^2}  \approx 4.79253.
$$
We do not know if these numbers $n_k$ have  a finite total stopping time. $~~~\bsq$\\

The methods of Applegate and Lagarias \cite{AL95c}
for $3x+1$ trees can potentially be applied to this problem, to further improve this
lower bound, and to establish  it for numbers $n$ having a finite total stopping time.\\

An interesting challenge is whether one can   show  for each $c>0$ that 
{\em only a density zero  set of $n$ have $\frac{\sigma_{\infty}(n; T_5)}{\log n}< c$.}
A stochastic model in \S\ref{sec8p9} 
predicts
that all but finitely many trajectories having  $\sigma_{\infty}(n) > 85 \log n $ 
will necessarily have
$\sigma_{\infty}(n) = +\infty$, so establishing this for $c=85$ would be consistent
with the prediction that only a density zero set of $n$ have $1$ in their forward
orbit under $T_5$. 

 
%
%
%

\subsection{$5x+1$ Size Statistics: Minimum Excursion Values}

In the topsy-turvy world of the $5x+1$ problem, 
since most trajectories get large, 
 our substitute for the maximum excursion constant  is the following 
 reversed notion.

%
%
%
 
 \begin{defi}~\label{de741}
 {\em 
For
an
 integer $n$ the {\em minimal excursion value } $t^{-}(n) $ of the
$5x+1$ function is given by 
\beql{741a}
t^{-}(n)   := \inf \{ |T_5^{(k)} (n)| : k \ge 0 \}.
\eeq
}
\end{defi}

We have $t^{-}(0)=0$, 
while  infinitely many $n$ will have minimum excursion value equal to $1$.\\

%
%
%
 
 \begin{defi}~\label{de741b}
 {\em 
For $n \ge 1$ the {\em minimal excursion constant} $\rho_5^{-}(n) $ of the
$5x+1$ function is given by 
\beql{741b}
\rho_5^{-}(n)   :=  \liminf_{n \to \infty} \frac{\log t^{-}(n)}{\log n}.
\eeq
}
\end{defi}

We now immediately have the following result.
%
%
%

\begin{theorem}~\label{th743a} 
The $5x+1$ minimum excursion constant is given by
\beql{741c}
\rho_5^{-} = 0.
\eeq
\end{theorem}

\paragraph{Proof.} The inverse orbit of $n=1$ for $T_5$ contains $\{ 2^j: j \ge 1\}$,
whence $t^{-}(2^j)= 1$. 
$~~~\bsq$\\

We state this easy result as a theorem,
 because it has the
remarkable feature, among all the constants associated
to these $3x+1$ and $5x+1$ maps, of being unconditionally proved!
It also has the interesting feature that the
stochastic models below make an incorrect prediction
in this case, cf. Theorem~\ref{th824}.
%
%
%

\subsection{$5x+1$ Count Statistics: $5x+1$ Tree Sizes}

In considering backwards iteration of the $5x+1$ function, we
can ask: given an integer $a$ how many numbers $n$ iterate forward
to $a$ after exactly $k$ iterations, 
that is,
$T_5^{(k)}(n)=a$?  \\

The set of backwards iterates of a given number $a$ can 
again
be pictured as a tree; we call
these {\em $5x+1$ trees}. 
Now
$N_k(a)$ counts the number of leaves at depth $k$ of 
the
tree with root node $a$, 
and $N_k^{\ast}$ counts the number of leaves in a {\em pruned $5x+1$ tree},
which is one from which all nodes with label $n \equiv 0 ~(\bmod ~5)$ have
been removed. The definitions are as follows.\\

%
%
%

\begin{defi}~\label{de751}
{\em 
 (1) Let $N_k(a; T_5)$ count the number of integers that forward iterate under the
 $5x+1$ map $T_5(n)$ to $a$ after exactly $k$ iterations, i.e. 
\beql{751a}
N_k(a; T_5) := |\{ n: ~T_5^{(k)} (n) =a \}|.
\eeq

(2) Let $N_k^{\ast} (a; T_5)$ count the number of integers not divisible by $5$ that forward iterate under the
 $5x+1$ map $T_5(n)$ to $a$ after exactly $k$ iterations, i.e. 
\beql{751b}
N_k
^{\ast}
(a; T_5) := |\{ n: ~T_5^{(k)} (n) =a, ~n \not\equiv 0 (\bmod~5) \}|.
\eeq
}
\end{defi}

The case $a=1$ is of particular interest, since the quantities then count integers that
iterate to $1$, and in this case we let 
$$N_{k, 5} := N_k(1; T_5), ~~~~~~~N_{k, 5}^{\ast} := N_k^{\ast}(1; T_5).
$$

%
%
%

\begin{defi}\label{d752}
{\em
(1) For a given $a$
the  $5x+1$ {\em tree growth constant} $\delta_5(a) $  for $a$ is given by
\beql{755}
\delta_5(a) :=\limsup_{k \to \infty} \frac{1}{k}\left( \log N_k(a; T_5)\right).
\eeq

(2) The  $5x+1$ {\em tree growth constant} $\delta_5= \delta_5(1).$
}
\end{defi}

The constant $\delta_5(a)$ exists and is finite, as follows from the 
same 
upper bound
as in 
 \eqn{253}. \\

The following result gives information on the sizes of depth $k$ pruned 
$5x+1$ trees over all possible tree types.\\
%
%
%

\begin{theorem}~\label{th753} {\em (Structure of Pruned $5x+1$ Trees)}

(1) For $k \ge 1$ and $a \not\equiv 0 (\bmod~5)$, the
 structure of the pruned level $k$ tree ${\overline{\sT}}_k^{\ast}(a)$, and
 hence the number $N_k^{\ast}(a; T_5)$,  is completely 
determined by $a~(\bmod~5^{k+1})$. \\

(2) There are $4 \cdot 5^k$ residue classes $a~(\bmod ~5^{k+1})$ with
$a \not\equiv 0 ~(\bmod~5)$. For these 
\beql{611aa}
\sum_{{a~(\bmod~5^{k+1})}\atop{a \not\equiv 0~(\bmod~5)}} N_k^{\ast}(a; T_5)= 4 \cdot 6^k.
\eeq 
\end{theorem}
It follows that if a residue class $a~(\bmod~5^{k+1})$ with $a \not\equiv 0 (\bmod~5)$ is
picked with the uniform distribution, the expected number of leaves in 
the random tree  ${\overline{\sT}}_k^{\ast}(a)$ is exactly $\left(\frac{6}{5}\right)^k$.

\paragraph{Proof.} This result is shown by a method exactly similar to
the $3x+1$ tree case (\cite[Theorem 3.1]{LW92}). 
We omit details.
 $~~~\bsq$\\

Theorem~\ref{th753} suggests  the following conjecture.

%
%
%

\begin{conj} \label{conj71}
 For each $a \not \equiv 0 ~(\bmod~5)$, the $5x+1$ tree growth constant  $\delta_5(a)$ is given by 
\beql{754a}
\delta_5(a) = \log \left(\frac{6}{5}\right).
\eeq
\end{conj}

Compare this conjecture with the prediction of Theorem~\ref{th852}.
%
%

\subsection{$5x+1$ Count Statistics: Total Inverse Iterate Counts}

In considering backwards iteration of the $5x+1$ function from an integer $a$, 
the complete data  is the set of integers that  contain $a$ in their forward orbit. 
The following
function describes this set. 

%
%
%

\begin{defi}~\label{de761}
{\em 
Given an integer $a$,  the {\em inverse iterate counting function}
$\pi_{a,5}(x)$ counts the number of integers $n$
with $|n|\le x$ that contain $a$ in their forward orbit under the $3x+1$ function.
That is
\beql{761a}
\pi_{a,5}(x) := \#\{ n:~ |n| \le x~~\mbox{such~that~some} ~T_5^{(k)}(n) = a, ~k \ge 0\}.
\eeq
}
\end{defi}

The inverse tree methods for the $3x+1$ problem carry over to the $5x+1$ problem,
so that one can obtain a result qualitatively of the following type, by similar proofs.

%
%

\begin{theorem}\label{th761a}
{\em (Inverse Iterate Lower Bound)}
There is a positive constant $c_1$ such that the following holds.
For each  $a \not\equiv~0~(\bmod~5)$, there is 
some
$x_0(a)$
such that for all $x \ge x_0(a)$, 
\beql{765}
\pi_{a
,5 
}
(x) \ge x^{c_1}.
\eeq
\end{theorem}

  The following statistics measure the size of the inverse iterate set in
  the sense of fractional dimension.
  
%
%
%

\begin{defi}~\label{de762}
{\em 
Given an integer $a$,  the {\em upper and lower $5x+1$ growth exponents} for
$a$ are given by 
$$
\eta_5^{+}(a) := \limsup_{ x \to \infty}  \frac{ \log 
\pi_{a
,5 
}
(x)
} {\log x},
$$
and
$$
\eta_5^{-}(a) := \liminf_{ x \to \infty}  \frac{ \log 
\pi_{a
,5 
}
(x)
} {\log x}.
$$
If  these quantities are equal, we define the {\em $5x+1$ growth exponent} $\eta_5(a)$
to be  $\eta_5(a) = \eta_5^{+}(a)= \eta_5^{-}(a)$.
}
\end{defi}

In parallel with conjectures for the $3x+1$ map,  we formulate the following 
conjecture.

%
%
\begin{conj}~\label{conj72}
{\em ($5x+1$ Growth Exponent Conjecture)} 
For all integers $a \not\equiv 0~(\bmod ~5)$, the $3x+1$ growth exponent
$\eta_5(a)$ exists, and takes a constant value $\eta_5$ independent of $a$.  This value satisfies
\beql{766a}
\eta_5 < 1.
\eeq
\end{conj}

The stochastic models discussed in \S\ref{sec8} suggest that the constant $\eta_5$ exists
and takes a value strictly smaller than $1$. 
There is some controversy concerning the conjectured value of the constant. 
In \S\ref{sec8} we present a repeated random walk  model and a branching random
walk model that both suggest  the value $\eta_5 \approx 0.649$.
A different branching random walk
model formulated by Volkov \cite{Vo06} suggests the value $\eta_5 \approx 0.678$.
Lower bounds toward this conjecture can be rigorously established, cf. Theorem~\ref{th761a}
above.
We have not bothered to determine $c_{1}$ in \eqn{765}, though we suspect it is well below either of the above predictions, and hence cannot distinguish between them.

%
%
%

 \section{$5x+1$ Function: Stochastic Models and Results}\label{sec8}

We now discuss stochastic models for the $5x+1$ problem paralleling those
for the $3x+1$ problem. These include random walk models for
forward iteration of the $5x+1$ map, analysis of the accelerated $5x+1$
map, and branching random walks for the backwards iteration of
the $5x+1$ map.

%
%
%

 \subsection{$5x+1$ Forward Iteration: Multiplicative Random Product Model}

Concerning forward iteration,
we may formulate a 
multiplicative random product model 
parallel  to that  in 
\S\ref{sec3}.  
Consider 
the random products
$$
Y_k := X_1 X_2 \cdots X_k,
$$
in which the  $X_i$ are each
independent identically distributed (i.i.d.) random
variables $X_i$ having the discrete distribution
$$
X_i = \left\{
\begin{array}{cl}
\df{5}{2}  & \mbox{with~probability} ~~\frac{1}{2}, \\
~~~ \\
\df{1}{2} & \mbox{with~probability}~~\frac{1}{2}.\\
\end{array}
\right.
$$
We call this the {\em $5x+1$ multiplicative random product} (MRP) model.\\

As before,
this model does not include the choice of starting value of the iteration, which would correspond
to $X_0$; the random variable $Y_k$ really models the {\em ratio}
$\frac{T_5^{(k)}(X_0)}{X_0}$. 
We define for $\lambda^{+}  \ge 1$ the {\em $\lambda^{+}$-stopping time random variable}
\beql{801c}
V_{\lambda}^{+} (\omega) := \inf\{k :~ Y_k \ge \lambda\}.
\eeq
where $\omega=(X_1, X_2, X_3, \cdots)$
denotes a sequence of random variables as above. This random vector $\omega$
 models the change in size
 of a random starting value $n=X_0$ that occurs on  iteration of the $5x+1$ map. \\

This stochastic model can be used to exactly account for  the density of 
$\lambda^{+}$-stopping times, as follows. 
\\

%
%

\begin{theorem}~\label{th800a} {\em  ($\lambda^{+}$-Stopping Time Density Formula)}
For the $5x+1$ map $T
_{5 } 
(n)$ and
any fixed
 $\lambda 
 >1$, the 
 natural density $\DD (S_{\lambda}(k))$ for integers having  $\lambda^{+}$-stopping
 time at most $k$  is given exactly  by the formula
\beql{802a}
\DD (S_{\lambda}^{+}(k)) = \PP[ V_{\lambda}^{+}(\omega)  \le k],
\eeq
in which $V_{\lambda}^{+}$ is the $\lambda^{+}$-stopping time random variable
in the multiplicative random product (MRP) model.
\end{theorem}

\paragraph{Proof.} 
This follows by a parallel argument to that in 
 Borovkov and Pfeifer \cite[Theorem 3]{BP00} for the $3x+1$ problem.
$~~~\bsq$\\

Theorem~\ref{th800a} is the stochastic model parallel of Theorem~\ref{th72a}.

%
%
%

 \subsection{$5x+1$ Forward Iteration: Additive Random Walk Model}


We next formulate  additive random walk models, obtained
after logarithmic rescaling of the $5x+1$ iteration.
The $5x+1$ iteration takes $x_0=n$ and $x_k= T_{5}^{(k)}(n).$ 
Using a logarithmic rescaling with 
$y_k= \log x_k$ (natural logarithm) we have 
$$
y_k= \log x_k :=\log T^{(k)} (n).
$$
Then we have
\beql{711b}
y_{k+1} =
\left\{
\begin{array}{cl}
y_k + \log \frac{5}{2}  + e_k & \mbox{if}~ x \equiv 1~~ (\bmod ~2 ),  \\ 
~~~ \\
y_k + \log \frac{1}{2} & \mbox{if} ~~x \equiv 0~~ (\bmod~2),
\end{array}
\right.
\eeq
with
\beql{712b}
e_k:= \log \left( 1+ \frac{1}{5 x_k}\right).
\eeq
Here $e_k$ is small as long as $|x_k|$ is large.\\

We approximate the deterministic process above with the following
 random walk  model with unequal size steps.
 We take random variables 
 $$
W_k  := -\log 2+ \delta_k  \log 5,
$$
in which  $\delta_k$ are i.i.d. Bernoulli random variables. 
The random walk positions  $\{ Z_{ k}: k \ge 0\},$ are 
then random variables having starting value
$Z_0 = \log m$, for some fixed initial condition $m >1$, 
and with 
$$
Z_k = Z_0 + W_1 + W_2 + \cdots + W_k.
$$
The $Z_k$ define a biased random walk, whose expected drift $\mu$ is given by 
$$
\mu: = E[W_k)] = - \log 2 + \frac{1}{2} \log 5= \frac{1}{2} \log \left( \frac{5}{4} \right) \approx 0.11157.
$$
The variance $\sigma$ of each step is given by
$$
\sigma:= {\rm Var}[ W_k] = \frac{1}{2} \log 5 \approx 0.80472.
 $$
 Call this random walk the {\em $5x+1$ Biased 
 Random 
 Walk 
 Model} ({\em $5
 x
 +1$ BRW 
 Model}).\\
 
 Since the mean of this random walk is positive, this biased random walk 
 has a {\em positive drift.} 
 This positive  drift implies that a random trajectory diverges with probability one.

%
%
%
\begin{theorem}~\label{th81b} 
For the $5x+1$ BRW model, with probability one, a  trajectory $\{Z_k: k \ge 0\}$ diverges to $+\infty$. 
 \end{theorem}
 
\paragraph{Proof.} This is an elementary fact about random walks with 
positive drift. $~~\bsq$\\
 
 This result implies that 
a generic
 trajectory 
has
 total stopping time 
 equal to
 $+\infty$.
 That is, 
 starting from $Z_{0}=\log n$,
  the probability $\PP[E_{n}]$
of the event $E_{n}$ 
 that
 for some $k\ge1$,
  the total stopping time condition  $Z_k \le 0$ 
 is satisfied,
 is 
 strictly smaller than $1$,
i.e., $\PP[E_{n}]<1$.
It is positive but decreases to $0$ as $n$
 increases to $+\infty$. 
 (To not confuse this fact with Theorem \ref{th81b}, even if $Z_{k}$ dips below $0$, it charges back up to infinity, almost surely.)
 \\
 
 To obtain a result parallel to those of 
 \S\ref{sec3} on the average behavior of
 numbers $n$ having a finite total stopping time, one needs to condition on the
 set of $n$ that have a finite total stopping time. 
 This  appears an
 approachable problem, but  requires a more complicated
 analysis than that given in \cite{LW92} or Borovkov and Pfeifer \cite{BP00}. \\
 
%
%
%

 \subsection{$5x+1$ Forward Iteration: Repeated  Random Walk Model}

 Next, paralleling 
 \S\ref{sec4}, we  formulate a {\em $5x+1$ Repeated Random Walk (RRW) model} as
 follows. A model trial is the countable set of random variables
 \beql{817a}
  \omega:= \{ Z_{k,n}: k \ge 0, n \ge 1\},
  \eeq
 having initial condition $Z_{0,n} = \log n$, with the individual random walks
 being $5x+1$ biased random walks, as above. In the following 
 subsections 
 we consider 
 other predictions that RRW model makes for various statistics.

%
%
%
\begin{theorem}~\label{th81c} 
For the $5x+1$ RRW model, with probability one, for
every $n \ge 1$ the  trajectory $\{Z_{k,n}: k \ge 0\}$ diverges to 
$+\infty$. 
 \end{theorem}
 
\paragraph{Proof.} This follows 
immediately
from Theorem~\ref{th81b}, since the
complement of this event is a countable union
of measure zero events. $~~\bsq$\\

 One might misinterpret the above as suggesting that the
   $5x+1$ RRW model 
predicts that
{\it all} 
trajectories are unbounded.
 Of course
this is an incorrect prediction. 
The $5x+1$ iteration has some finite cycles, and
 furthermore  there are
 infinite number of integers that eventually enter one of these cycles. 
  The stochastic model above 
  cannot account for such bounded trajectories! 
  Instead we interpret the stochastic model prediction to 
 be that  a 
 {\em density one} set of integers lie on unbounded trajectories. \\

 This should make you {\it very worried} about relying on stochastic models to predict  that $3x+1$ trajectories 
 decay! There could potentially be a set of measure zero escaping to infinity, which the model
 simply cannot see.
 Such a pathological trajectory is the {\it heart and soul} of the $3x+1$ problem, and root cause 
 of its difficulty! \\

%
%
%

\subsection{$5x+1$ RRW Model Prediction: Minimum Excursion Constant}

The $5x+1$ RRW model has   the following analogues of 
 minimal excursion values and of the minimum excursion constant.

%
%
%
 
 \begin{defi}~\label{de822}
 {\em 
For a realization $\omega= \{Z_{k,n}: k \ge 0, n \ge 1\}$ of the $5x+1$ RRW
model,  the {\em minimal excursion value } $t^{-}(n, \omega)$ is given, for each $n \ge 1$, by  
\beql{822a}
t^{-}(n, \omega)   := \inf \{e^{Z_{k,n}}: k \ge 0\}.
\eeq
}
\end{defi}

Theorem~\ref{th81c} implies that with probability one the value $t^{-}(n, \omega)$ is well-defined
and strictly positive.

%
%
%
 
 \begin{defi}~\label{de823}
 {\em 
For  a realization $\omega$ of the $5x+1$ RRW
model,  the {\em minimum excursion constant  } $\rho_5^{-}(\omega)$
is given by 
\beql{823a}
\rho_5^{-}(\omega) 
 := \liminf_{n \to \infty} \frac{\log t^{-}(n, \omega)}{\log n}.
\eeq
}
\end{defi}

Now a large deviations analysis yields the following result.
%
%
%

\begin{theorem}~\label{th824} 
{\em ($5x+1$ RRW Minimum Excursion Constant)}
 For the $5x+1$ 
RRW model, with probability one   the quantities $t^{-}(n, \omega)$
are finite for every $n \ge 1$. In addition, with probability one
the random quantity 
\beql{831a}
\rho_{5,\RW}^{-}(\omega):= \liminf_{n \to \infty} \frac{\log t^{-}(n; \omega)}{\log n}
= \liminf_{n \to \infty} \left( \inf_{k \ge 0} \frac{Z_{k, n}}{\log n} \right)
\eeq 
equals the constant
\beql{832a}
\rho_{5,\RW}^{-}=1- \frac{1}{\theta^{\ast}} \approx-1.86466 ,
\eeq
in which  $\theta^{\ast} \approx 0.3490813
$ is the larger of the two real roots of the equation
$M_{5,\RW}(\theta)=1$, where
$M_{5,\RW}(\theta): = \frac{1}{2}\left( 2^{\theta} + (\frac{2}{5})^{\theta}\right)$
is a  moment generating function associated to  the random walk.
\end{theorem}

\paragraph{Proof.} This is proved by a large deviations argument similar to
that used for the maximum excursion constant for the $3x+1$ problem in
Lagarias and Weiss \cite[Theorem 2.3]{LW92}. We sketch the main computation.
We estimate the probability $P(r, H, x)$ on a single trial starting at $\log x$ of having
$$-Z_{r \log x,\log x} \ge H \log x.$$
We define $a$ by the condition  $H= ar$ and 
find that  the probability is  given by Chernoff's bound as
$$
P(r, H, x) = \exp \left( -g_{5,\RW} (a) r \log x (1+o(1) \right),
$$
in which 
\beql{858a} 
g_{5,\RW} (a) := \sup_{\theta \in \RR} \left( a \theta - \log M_{5,\RW}(\theta)\right)
\eeq
is a large deviations rate function,
which is the Legendre transform of the logarithm of the moment generating function
$M_{5,\RW}(\theta) = \frac{1}{2}\left( 2^{\theta}+ \frac{5}{2})^{\theta}\right)$.
The repeated random walk makes $x$ trials $1 \le n \le x$ so the probability of a success
over these trials is $x P(r, H, x)$,
 and we want this to be at least $x^{\epsilon}$, so that  a success occurs infinitely often
 as $x \to \infty.$ (We also will let $\epsilon \to 0$, so we set it equal to zero in what follows.)
 We want therefore to maximize $H=ar$ subject to the constraint that $g_{5,\RW}(a) r \le 1$.
To maximize we may take $g(a) r=1$, whence $r= \frac{1}{g(a)}$ can be used to eliminate 
the variable $r$.
 We now have the maximization problem to maximize $H:= \frac{a}{g_{5,\RW}(a)}$
 over $0< a< \infty.$ 
 One finds an extremality condition for maximization which
 yields
 $$
 H^{\ast} = \frac{1}{\theta(a^{\ast})},
 $$
 where $a^{\ast}$ achieves the maximum, and $\theta^{\ast}$ is the
 corresponding value in the Legendre transform.  Uniqueness of the maximum
 follows from convexity properties of the function $\log M_{\RW}(\theta)$.
 Detailed error estimates are also needed to verify that
 this  the maximum gives the dominant
 contribution. 
$~~~\bsq$ \\

This constant $\rho_{5,\RW}^{-}$ found in Theorem \ref{th824} is {\em negative}, 
i.e. the minimum excursion in the model reaches a real number much  smaller than $1$!
As  a prediction for the $5x+1$ problem, this disagrees 
with the exact answer for minimum excursion constant for the 
$5x+1$ problem $\rho_{5}^{-}=0$ given in Theorem ~\ref{th743a}. \\

We view this inaccurate prediction as stemming from  the discrepancy that the $5x+1$ 
function  takes
only values on the integer lattice, and that  its additive correction term is not accounted
for in this stochastic model. That is, the stochastic model will not necessarily make
good predictions on behavior of an orbit once an orbit reaches a small value, e.g. $|x|< C$ for any fixed
constant $C$.
We may hope that the $5x+1$ model still makes an accurate prediction
concerns the question: how many integers  reach some small  value, for example
reaching the  interval $|x|<C$.

%
%
%

\subsection{$5x+1$ RRW Model Prediction: Total Stopping Time Counts}

W
e can 
interpret th
e false prediction above
for minimum excursions  in 
a constructive way:  
as soon as a
$5x+1$ 
 trajectory achieves a size
$e^{Z_{k, n} }<1$,  it enters a periodic orbit. 
Therefore 
this condition
can be treated 
as a ``stopping time" condition that detects when a trajectory reaches the value $1$.


%
%
%

\begin{theorem}~\label{th835} 
{\em ($5x+1$ RRW Total Stopping Time Counts)}
 For the $5x+1$ 
RRW model and for a given $\omega$,
let
$$
S_{\infty}(\omega): = \{ n \ge 1:~~e^{Z_{k,n}} < 1  ~\mbox{holds~for some}~~k \ge 1\}
.
$$
Collect those seeds $n$ whose trajectory according to $\gw$ ``reaches 1''. 
Let $\pi_{5}(\,\cdot\,;\gw)$ denote the corresponding counting function,
%
$$
\pi_5( x 
;
 \omega) := \#\{ 1 \le n \le x: n \in S_{\infty}(\omega)\}
.
$$
Then
$$
\lim_{x \to \infty} \frac{ \log \pi_5(x; \omega)} {\log x} = \eta_{5,\RW}
,
\qquad
\mbox{
for almost every }
\gw.
$$
Here
$\eta_{5,\RW} \approx 0.65049$ is given by
$\eta_{5,\RW} = 1 - \theta_{5,\RW}$ where $\theta_{5,\RW}\approx 0.34951$
is the unique 
positive
solution to the equation  
\beql{833a}
M_{5,\RW}(\theta) := \frac{1}{2}\left( 2^{\theta}+ \left(\frac{5}{2}\right)^{\theta}\right) 
=
1
.
\eeq
 \end{theorem}

\paragraph{Proof.} 
This can be proved by a large deviations model similar in nature to those 
considered in Lagarias and Weiss \cite[Theorem 2.4]{LW92}.
We sketch the main estimate.
For $k=r \log x$,
consider
 the probability $P(r, x)$ that 
for a single random walk $e^{Z_{k, \log x} } < 1.$ 
Since we make $x$ draws for $1 \le n \le x$  in the repeated random walk,
the expected number
of such individuals satisfying this property 
will be $x P(r, x)$. 
This probability is estimated using Chernoff's bound to be
$$
P(r,x) = \exp \left( -g_{5,\RW} (a) r \log x (1 + o(1)\right),
$$ 
where $a= \frac{1}{r}$, and $g_{5,\RW}$ is the large deviations rate function
\eqn{858a} in Theorem~\ref{th824}.
We now maximize this probability over $r$. To do this we eliminate $r$ using
$r =\frac{1}{a}$, so we want to determine 
$$
\tau_{5,\RW}:= \min_{ 0 \le a < \infty} { \frac{g_{5,\RW}(a)}{a}}.
$$
Then we obtain $x P(r, x) \le x^{ 1 - \tau_{5,\RW} + o(1)}$ for all $r$, with equality
holding for $r= \frac{1}{a^{\ast}}$ where
$a^{\ast}$ be the value that attains  the maximum of 
 $f(a) := \frac{g_{5,\RW}(a)}{a} $ taken on the positive half-line.
The extremality conditions for the minimum 
leads to the condition $M_{\RW}(\theta(a^{\ast}) )=1,$ where $\theta$ is the
Legendre transform variable, and also to the identity 
$$
\tau_{5,\RW} = \frac{g_{5,\RW}(a^{\ast})}{a^{\ast}}= \theta(a^{\ast}) := \theta_{5,\RW}.
$$
The strict convexity of  the function $\log M_{\RW}(\theta)$ is used to get a unique
minimum, with  $\eta_{5,\RW} = 1- \tau_{5,\RW}$. 
For
a rigorous proof, one must control various error estimates to show the 
dominant contribution to the probability comes from a small region near $a^{\ast}$. 
$~~~\bsq$

\paragraph{Remark.} The value of $\theta_{5,\RW}$
in the minimization problem in the proof of Theorem ~\ref{th835}
turns out to be identical to that in  the maximization problem that is needed for
proving Theorem~\ref{th824}.


%
%
%

\subsection{$5x+1$ Accelerated Forward Iteration:  Brownian Motion.}

Kontorovich and Sinai \cite{KS02} extended the Structure Theorem (that is, Theorems \ref{th51} and \ref{th52}) and the consequences on the Central Limit Theorem (Theorem \ref{th53}) and geometric Brownian motion (Theorem \ref{th54}) to a class of functions which they called $(d,g,h)$-maps. The case $d=2$, $g=5$, and $h=1$ corresponds to the accelerated $5x+1$ function, $U_{5}(n)$.\\

The analogous distribution and Central Limit Theorems are proved in the same way, leading to the following.\\

%
%

\begin{theorem}~\label{th82b}
(Geometric Brownian Motion) 
The rescaled  paths of the
accelerated
 $5x+1$ map are those of a geometric Brownian motion with drift $\log(\frac54)$. By this we mean the following.

For 
an initial seed $x_{0}$ which is relatively prime to both $2$ and $5$,
denote its iterates by $x_{\kk}:=U_{5}^{(\kk)}(x_{0})$, 
let $y_{\kk}:=\log x_{\kk}$ and define the scaled variable
$$
\gw_{k}
:=
{
y_{k}
-y_{0} 
-k\log (\frac54)
\over
\sqrt{2k}\log 2
} 
.
$$
Partition the interval $[0,1]$ as $0=t_{0}<t_{1}<\cdots<t_{r}=1$, and set $k_{j}=\lfloor t_{j}k\rfloor$.
Then for  any $a_{j}<b_{j}$, $j=1,\dots,r$, 
$$
\lim_{k\to\infty}
\PP
\left[
x_{0}:
a_{j}<
{
\gw_{k_{j}}
-
\gw_{k_{j-1}}
}
<
b_{j}
,
\mbox{ for all }
j=1,2,\dots,r
\right]
=
\prod_{j=1}^{r}
\bigg(
\Phi(b_{j})-\Phi(a_{j})
\bigg)
,
$$
where  $\Phi(a)$ is the cumulative distribution function for the standard normal 
 distribution.

\end{theorem}

\paragraph{Proof.} This  is a consequence of  Theorem 5 in Kontorovich-Sinai \cite{KS02}.
~~~ $\bsq$  \\

\paragraph{Remark.}
The accelerated drift, $\log(\frac54)$, is again double that of  the Biased Random Walk model, which predicts a drift of  $\frac12\log(\frac54)$. A  zero-mean, unit-variance Wiener process $W_{t}$ satisfies the ``law of iterated logs'' almost surely, that is:
$$
\limsup_{t\to\infty}{|W_{t}|\over \sqrt{2t\log\log t}} =1,
$$
with probability $1$. Hence the drift being positive implies that almost every $5x+1$ trajectory escapes to infinity (yet 
we emphasize again that 
we do not know how to prove this for a single given trajectory!). 

%
%
%

 \subsection{$5x+1$ Backwards Stochastic Models: Branching Random Walks }
  
 We next formulate branching random walks to model the $5x+1$ iteration in exact
 analogy with the $3x+1$ models. We denote these models $\sB[5^j]$ for $j \ge 0$.\\

{\em $5x+1$ Branching Random walk ~$\sB[5^0]$.} There is one type of
 individual. 
 With probability $\frac{4}{5}$ an individual has a single offspring
 located at a position shifted by $\log 2$ on the line from its progenitor, and with probability $\frac{1}{5}$
 it has two offspring located at positions shifted $\log 2$ and $\log \frac{2}{5}$ on
 the line from their progenitor. If the progenitor is in generation $k-1$, the offspring are
 in generation $k$. The tree is grown from a single individual in generation $0$,
 the root, with specified initial location $\log a$.\\

The more general models  for $j \ge 1$ are given as follows.\\

 {\em $5x+1$ Branching Random walk ~$\sB[5^j], (j \ge 1)$.} There are $p=4 \cdot 5^{j-1}$
 types of individuals, indexed by residue classes $a~(\bmod~5^j)$ with 
 $a \not\equiv 0~(\bmod~5)$. The distribution of offspring of an individual of
 type $a~(\bmod~5^j)$, at any given generation (or depth) $k$ in
 the branching, is determined as follows: Suppose $a~(\bmod~5^{j})$ 
 is the type of a node at depth $k-1$. Now regard it 
 as being, with probability $\frac{1}{5}$ each, one of the five possible residue
 classes $\tilde{a}~(\bmod~5^{j+1})$ consistent with its class $(\bmod~5^j)$. A tree of depth $1$
 having $\tilde{a}$ as root node, 
  then has either one or
 two progeny, at depth $1$,  given by $(T^{\ast})^{-1}(\tilde{a})$, 
 whose node labels are well-defined classes $(\bmod~5^j)$,
 either $2 \tilde{a}$ or, if it legally occurs, $\frac{2 \tilde{a}-1}{3} (\bmod~5^j)$. 
 The branching random walk then
 produces an individual of type $2\tilde{a}$ 
 at generation $k$ whose position  is additively shifted by $\log 2$ from that of the
 generation $k-1$ progenitor node of type $\tilde{a}$
 plus, if legal, another node of type 
   $\frac{2 \tilde{a}-1}{5} (\bmod~ 5^j)$, which is shifted in position by
    $\log(\frac{2}{5})$ on the line from that of the generation $k-1$-node.
    The tree is grown from a single individual at depth $0$, with
    specified type and location $\log a$.  \\

Just as in the $3x+1$ branching random walk models, 
the  behavior of the random walk part of the model can
completely reconstructed from knowing  the type of each node.\\

For the rest of this section, let $\omega$ denote a single realization of such a branching random walk 
 $\sB[5^j]$  which starts from a single individual $\omega_{0,1}$ of type $1~(\bmod~5^j)$
 at depth $0$, with initial position label $\log |a|$. Here $\omega$ describes
 a particular infinite tree. We let $N_k(\omega)$
 denote the number of individuals at level $k$ of the tree. We let
 $S( \omega_k, j)$
 denote the position of the $j$-th individual at level $k$ in the tree, for
 $1 \le j \le N_k(\omega)$.  \\

These models are  {\em supercritical branching processes}
exactly as for the $3x+1$ case: In {\em every} random realization $\omega$,
the number of
nodes at level $d$  grows exponentially in  $d$, and there are no 
extinction events.\\

In terms of growth of trees of inverse iterates, these models will accurately
represent certain features of $5x+1$ trees, and not others. They might
accurately describe tree sizes. However 
these branching random walks very likely
do not accurately model  positions of inverse iterates of the $5x+1$ in certain
crucial ways. Namely, individuals whose branching walk position is negative (corresponding
to a  $5x+1$ iteration value $x$ falling in the interval $ (0, 1)$)
are where the correction term $e_k$ in \eqn{712b} in the $5x+1$ iteration becomes
significant, breaking the size connection of the model iterates and 
the $5x+1$ iterates. \\

We now give some quantities of the trees associated to a
realization $\omega$ of the branching random walk $\sB[5^j]$.
We let $N_k:=N_k(\omega)$
denote the number of individuals in generation  $k$, and   
let $\{ \omega_{k, i} : 1 \le i \le N_{k}(\omega)\}$ denote the set of all individuals in generation $i$,
ordered by their branching random walk locations on the line, denoted 
$$
L(\omega_{k
,
1}) \le L(\omega_{k,2}) \le 
\cdots
 \le L(\omega_{k, N_k}).
$$

The  {\em size} of the element $\omega_{k,i}$, viewed as analogues of the
$5x+1$ iterates, is the exponentiated quantity
 \beql{867a}
Z_{k,i}:= e^{L(\omega_{k, i})}.
\eeq
The branching random walk has the property that
the 
sizes of  most individuals
in a tree will tend to get larger. 
(This initially seems rather surprising, but note that if a forward orbit is unbounded, then necessarily
all backward orbits leading to it must be unbounded as well!)
We are interested in individuals whose size under the $5x+1$ iteration is around 
a given value $x$. The tree models will detect individuals whose size is larger than $x$.\\

In the following subsections we address for the
$5x+1$ branching random walk models the following  questions.\\

1. What is the exponential growth rate of the quantities $N_k(\omega)$, as a function of $k$?\\

2. What is the maximum level $k$ that has some individual $Z_{k,i} \le x?$
This requires analyzing the size of the {\em first birth location} $L(\omega_{k,1})$.\\


3. How  does  the total number of individuals $\pi_{5}(x; \omega)$
in the $5x+1$ ttree having location $Z_{k,i} \le x$ grow 
as a function of $x$?\\

%
%
%

 \subsection{Backwards Iteration Prediction: $5x+1$ Tree Counts}\label{sec8p8}
  
 The size of $5x+1$ trees can be estimated for these models $\sB[5^j]$ , as follows. 

%
%
%
\begin{theorem}~\label{th852}
{\em ($5x+1$ Stochastic Tree Size)} 
For all $j \ge 0$  a realization $\omega$ of a tree grown in
the $5x+1$ branching random walk model $\sB[5^j]$ satisfies 
\beql{866a}
\lim_{k \to \infty} \frac{1}{k} \left( \log N_k(\omega)\right) = \log \left(\frac{6}{5}\right), ~~~~\mbox{almost ~surely.}
\eeq
\end{theorem}

\paragraph{Proof.} This is proved in exactly similar fashion to the $3x+1$ stochastic
model case in Lagarias and Weiss \cite[Corollary 3.1]{LW92} $~~~\bsq$ \\

This result only uses the Galton-Watson process branching
structure built into the branching random walk  $\sB[5^j].$
It does not depend on the sizes of the iterates.\\

The conclusion of Theorem \ref{th852},  viewed
as a prediction of the growth behavior of $5x+1$ trees,
 is consistent with the rigourous results on average tree size for pruned $5x+1$ trees
given in Theorem~\ref{th753}.\\

%
%
%

 \subsection{Backwards Iteration Prediction: Extremal Finite Total Stopping Times}
 \label{sec8p9}
  
As indicated above, most integers for the
$5x+1$ map will not have a finite total stopping time. However it is of interest
to analyze the small subset of integers that do have a total stopping time; these are
exactly the integers in the tree of inverse iterates of $a=1$. 
We analyze what is the maximum generation $k$ that contains an individual
having  size $e^{L(\omega_{k, i})} \le x.$

Denote the location of this first birth individual 
in generation $k$ by $L_k^{\ast}(\omega): =L(\omega_{k,1})$, for a given realization $\omega$ of the random walk.

%
%
%
\begin{theorem}~\label{th871a} 
{\em (Asymptotic $5x+1$ First Birth Location)}
 There is a constant $\beta_{5, BP}$ such that, for all $j \ge 1$,
 the branching random walk model $\sB[3^j]$ has asymptotic first birth (leftmost birth)
 \beql{871a}
 \lim_{k \to \infty} \frac{1}{k}L_k^{\ast}(\omega) = \beta_{5,BP} ~~~~~\mbox{a.~s.}
 \eeq
 This constant $\beta_{5, BP} \approx 0.01179816 
 $ is determined uniquely
 by the properties that it is the unique constant with $\beta>0$ that
 satisfies
 \beql{872a}
 \overline{g}_{5, BP}(\beta)=0,
 \eeq
 where
 \begin{eqnarray}\label{873a}
 \overline{g}_{5, BP}(a) 
 &:= &
 - \sup_{\theta \le 0} \left( a \theta - 
 \log \left( 2^{\theta} + \frac{1}{5}(\frac{2}{5})^{\theta}\right) \right)
.
\end{eqnarray}
 \end{theorem}

%
%

\begin{figure}
\begin{center}
\includegraphics{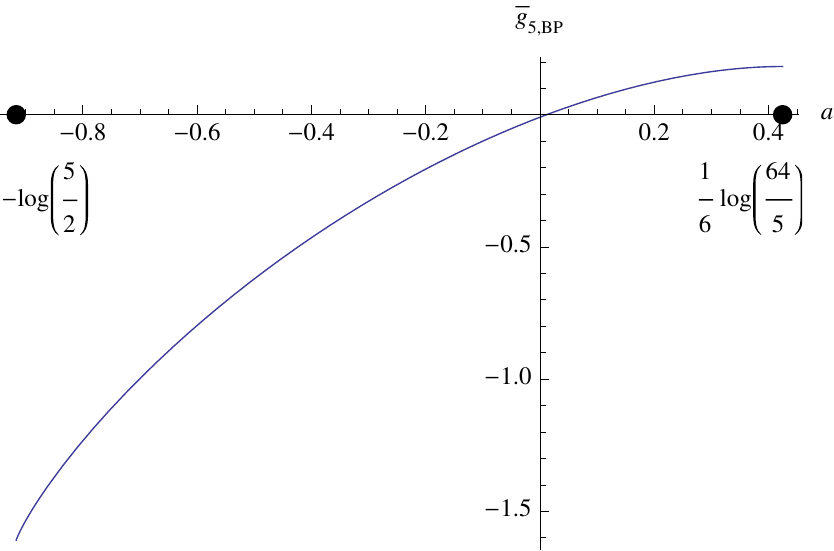}
\end{center}
\caption{
A plot of $a$ versus $\bar g_{5,BP}(a)$, in the range $\log(2/5)<a<\frac16\log(64/5)$.
}
\label{fig81}
\end{figure}

%
%

\begin{figure}
\begin{center}
\includegraphics{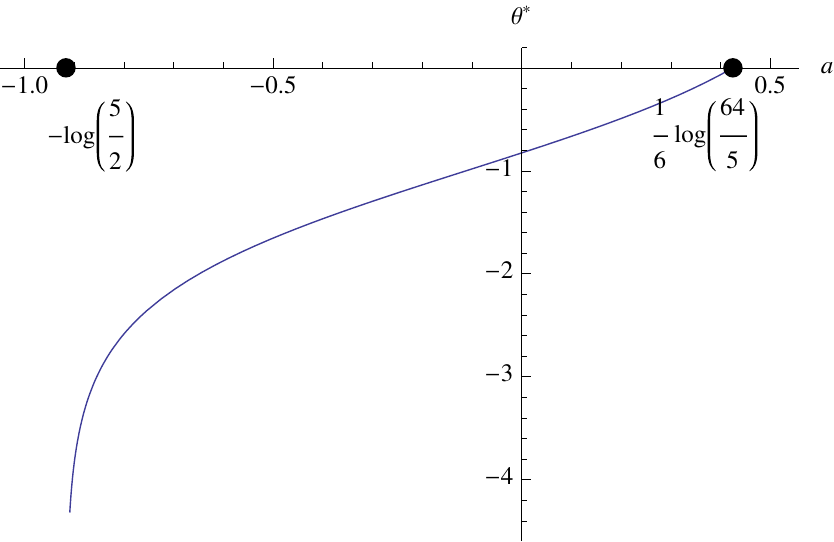}
\end{center}
\caption{
A plot of $a$ versus $\gt^{*}$, in the range $\log(2/5)<a<\frac16\log(64/5)$.
}
\label{fig82}
\end{figure}

\paragraph{Proof.} 
This is proved by 
an argument analogous
to the $3x+1$ case analyzed
in Lagarias and Weiss \cite[Theorem 3.4]{LW92}, cf. Theorem~\ref{th72}.  Here
we use a branching process (inverse) moment generating function
\beql{873ab}
M_{5, BP} (\theta):=2^{\theta} + \frac{1}{5}\left(\frac{2}{5}\right)^{\theta}.
\eeq
in computing the rate function $\overline{g}_{5,BP}(a)$.
We note that $\overline{g}_{5, BP}(a)$ is increasing for 
$\log \frac{2}{5} < a < \frac{1}{6} \log \frac{64}{5}$, 
(see Figure 
\ref{fig81})
and on this range the
value $\theta^{*}:= \theta(a)$ achieving the extremum in \eqn{873a} is 
an increasing 
function of $a$, reaching the value $\theta=0$ at the upper endpoint
(see Figure 
\ref{fig82}). 
We have
$\overline{g}_{5, BP}(a)=
\log (\frac{6}{5}) $ for $\frac{1}{6} \log \frac{64}{5}\le a < \infty$. 
%
%
$~~~\bsq$\\

Now one defines a {\em branching random walk stopping limit}
$$
\gamma_{5, BP}(\omega)  := \limsup_{k \to \infty} \frac{k}{L_k^{\ast} (\omega)}. 
$$
Theorem~\ref{th871a} 
implies that this value is constant almost surely, equaling a value $\gamma_{5, BP}$
given by 
\beql{887}
\gamma_{5, BP} = \frac{1}{\beta_{5, BP}} \approx 84.76012.
\eeq

One can show the constants $\gamma_{5, BP}$ and $\gamma_{5,\RW}$ agree,
just as for the $3x+1$ stochastic models.
%
%
%

\begin{theorem}~\label{th873} 
{\em ($5x+1$ Random Walk-Branching Random Walk Duality)}
The $5x+1$ repeated random walk (RRW)  scaled stopping time limit 
$\gamma_{5,\RW}$ and the branching random walk
stopping limit $\gamma_{5, BP}$ for the $5x+1$ branching random walk (BP)
model $\sB[5^j]$ with $j =0$, are related by 
\beql{731} 
\gamma_{5,\RW}= \gamma_{5, BP}.
\eeq
\end{theorem}

\paragraph{Proof.} 
This result is proved using  a relation between moment generating functions
$$
M_{5, BP}(\theta) = M_{5,\RW}(\theta+1),
$$ 
compare \eqn{833a} and \eqn{873a}.
It is identical in spirit to the proof in Lagarias and Weiss \cite[Theorem 4.1]{LW92}.$~~~\bsq$\\

The analogue of  this result applied to 
the $5x+1$ problem would be the  following heuristic prediction:
{\em  For any constant $\gamma > \gamma_{5, BP}$  all but finitely many
trajectories having total stopping time $\sigma_{\infty}(n) > \gamma \log n$ necessarily have
$\sigma_{\infty}(n) = +\infty$.} We could take $\gamma=85$, for example.

%
%
%

 \subsection{Backwards Iteration Prediction: Total Preimage Counts}

  The following result gives, for the simplest branching random walk model $\sB[5^0]$, 
 an almost sure asymptotic  for
 the number of inverse
iterates of size below a given bound.

%
%
%
\begin{theorem}~\label{th809a} 
{\em (Stochastic Inverse Iterate Counts)}
For a realization $\omega$ of the branching random walk $\sB[1]$, let
$I^{\ast}(t; \omega)$ count the number of progeny located at positions  $Z(\omega_{k,j}) \le x$, i.e.
\beql{896}
I^{\ast}(x; \omega) := \# \{ \omega_{k,j}:  Z(\omega_{k,j}) \le x, \mbox{for~any}~~k\ge 1, ~ 1 \le j \le 
N_k(\omega)\}.
\eeq
This quantity satsfies with probability one  the asymptotic estimate
\beql{897}
I^{\ast}(x; \omega) = x^{\eta_{5, BP}+o(1)}~~~\mbox{as}~~ x \to \infty,
\eeq
in which  $\eta_{5, BP}\approx 0.650919$  is the maximum value of
 $f(a):= \frac{1}{a} {\overline{g}_{5, BP}(a)}$ taken over the interval $0 \le a < \frac{1}{6} \log \frac{64}{5}$.
\end{theorem}

\paragraph{Proof.} This is proved by a large deviations argument similar to that used in 
 Lagarias and Weiss \cite[Theorem 4.2]{LW92}. 
One counts the number of progeny at level $k$ for each level $k$ satisfying the bound,
by estimating the probability that a random leaf satisfies the appropriate
large deviations bound. One
shows that this number peaks for $k \approx 
 \theta_{5, BP} \log x$,
where $\theta_{5, BP} =\frac{1}{a^{\ast}} \approx 9.19963,$
 where 
$a^{\ast}  \approx 0.1087$
is the value of $a$  achieving the maximum 
above. One shows that the right side is an upper bound
for all levels $k$, and that the sum total of levels $k > 100 \log x$ contribute negligibly
to the sum.
$~~~\bsq$ \\

The model  statistic $I^{\ast}(x; \omega)$ functions as  a proxy for the 
$5x+1$ count function $\pi_{a}^{\ast}(x)$,
where $\log |a|$ gives the position of the root node of the branching random walk. 
This result is the stochastic analogue of Conjecture~\ref{conj21} 
about
the $3x+1$ growth exponent.  The argument above also makes the 
 prediction is  that the levels $k$ at which the bulk of
the members of $\pi_{a}(x)$ occur has $k \approx \frac{1}{a^{\ast}} \log x$.

  \paragraph{Remark.} An entirely different set of branching random walk models 
  has been developed
  by S. Volkov \cite{Vo06} to model
  the   $5x+1$ problem. 
  Volkov  models counting  all  non-divergent trajectories of the $5x+1$ problem, 
  which are those which
  enter some  finite cycle, and denotes the number of these below $x$
 by $Q(x)$. Thus $\pi_5(x) \le Q(x)$, and conjecturally these should be of
similar orders of growth. It is expected there are finitely many cycles,
and each should absorb roughly the same number of integers below $x$,
in the sense of the exponent in the power of $x$ involved.

Volkov's   branching process stochastic models  grow a complete
 binary tree, rather than a tree that may have either one or two branches from each
 node, as in the models above. 
 He suggests that the $5x+1$ problem can be modeled by such  trees,
 using an unusual encoding of the iterates (some edges encode several iteration
 steps of the inverse Collatz function).
  In order to do this, his node weights are 
 chosen differently than above. 
 He arrives at a predicted  exponent $\eta_{5, BP}^{\ast} \approx 0.678,$
 which differs from the prediction $\eta_{5, BP}\approx 0.650919$ made
 in Theorem ~\ref{th809a} above. 
 The  empirical data Volkov presents seems insufficient to discriminate between these
 two predicted exponents. It would be interesting for this problem to be investigated
 further.

%
%
%
%

 \section{Benford's Law
 for $3x+1$ and $5x+1$ Maps}\label{sec9}

Another curious statistic satisfied by the $3x+1$ function was discovered by Kontorovich and Miller \cite{KM05}: {\it Benford's Law}. \\

In the late 1800s, Newcomb \cite{New}
noticed a surprising fact
while perusing tables of logarithms
: certain pages were significantly more worn than others.
Numbers whose logarithm started with 1
were being referenced
 more frequently
than other digits.
Instead of observing one-ninth (about 11\%) of entries having a
leading digit of 1, as one would expect if the digits $1, 2,
\dots, 9$ were equally likely, over 30\% of the entries had
leading digit 1, and about 70\% had leading digit less than 5.
Since $\log _{10}2\approx 0.301 $ and $\log _{10}5\approx 0.699$,
Newcomb speculated that the probability of observing a digit less
than $k$ was $\log _{10}k$.
This logarithmic
phenomenon became known as {\it Benford's Law} after Benford \cite{Ben}
collected and in 1938 popularized extensive empirical evidence of this distribution in
diverse data sets. \\

Benford's law 
seems to hold for many sequences of numbers generated
by dynamical systems having an ``expanding" property, see Berger et al \cite{BBH05}
and Miller and Takloo-Bighash  \cite[Chap. 9]{MT06}.
Benford  behavior has been  empirically observed for initial digits of the 
first 
iterates of the $3x+1$ map or accelerated $3x+1$ map
for  a randomly chosen initial number $n$.  Here we survey some rigorous  theorems
quantifying this statement, for initial iterates. Similar Benford results
can be proved for the $5x+1$ function. \\

We emphasize that the Benford law behavior quantifed here concerns behavior on a 
fixed finite set of
initial iterates of these maps. 
Indeed, the  $3x+1$ conjecture
predicts that Benford behavior  cannot hold for the full infinite set of 
forward iterates, since conjecturally they become
periodic! However it remains possible 
that a strong form of Benford behavior could hold on (infinite) divergent orbits of
the $5x +1$ problem.\\

%
%
%
\subsection{ Benford's Law and Uniform Distribution of Logarithms}

To make Benford's  law precise,
we say that the {\it mantissa} 
function $\cM(n)\in[1,10)$ is the leading entry of $n$ in ``scientific notation'', that is, $n=\cM(n)\cdot 10^{{\lfloor \log_{10} n \rfloor}}$. 
Benford's law concerns the distribution of leading digit of the mantissa, while one can also
consider the distribution of the lower order digits of the mantissa.
%
%
%

\begin{defi}\label{def901a}
{\em 
An infinite  sequence $\{n_{1},n_{2},\dots,n_{k},\dots\}$  satisfies the 
{\em strong Benford's Law} (to base $10$)   if  the logarithmic digit frequency holds for any
order digits in the mantissa. That is, for any $a\in[1,10)$,
\beql{900f}
\lim_{x\to\infty}
{\#\{k\le x:\cM(n_{k})<a\}\over x }
= \log_{10}(a) .
\eeq
}
\end{defi}

The strong version of
Benford's law is well known to be equivalent to uniform distribution $\bmod~1$ of the base 10 logarithms
of the numbers in the sequence, cf. Diaconis \cite[Theorem 1]{Dia77}. 
%
%
%

\begin{theorem} 
{\em (Strong Benford  Law Criterion)}
A sequence $\{n_{1},n_{2},\dots\}$ satisfies the strong Benford's Law 
(or ``is strong Benford'') to base $10$
if and only if the sequence $\{\log_{10}  n_{1},\log_{10} n_{2},\dots\}$ is 
equidistributed $(\bmod~1)$, that is, for any $a\in[0,1)$,
\beql{1001}
\lim_{x\to\infty}
{\#\{k\le x:\log_{10} n_{k}(\bmod~ 1)<a\}
\over{x}}=a.
\eeq
\end{theorem} 

The definition and theorem above extend to expansions
in any integer base $B \ge 2$.
This result  suggests the following  general  definition of strong Benford's
Law to any {\em real} base $B >1$.

%
%
%
 
 \begin{defi}~\label{de81a}
 {\em 
Let $B >1$ be a real number. A sequence $\{n_{1},n_{2},\dots,n_{k},\dots\}$ satisfies 
the {\em strong Benford's Law to base $B$}  if and only if the
sequence $\{ \log_B(n_1), \log_B(n_2), ...\}$ is uniformly distributed modulo one.}
\end{defi} 

 This definition is equivalent to  the earlier one for
integers expanded in a radix expansion to 
any
base $B > 
1
$.
One can similarly define the mantissa function to any real base $B>1$, extending Definition \ref{def901a}.
\\

Benford's Law is stated for infinite sequences. However one 
 can obtain approximate results that apply to  finite sequences $\{x_1, x_2, ..., x_k\}$, by
using the following discrepancy  measure of approximation to uniform distribution of such sequences.

%
%
%
 
 \begin{defi}~\label{de102}
 {\em 
 Given a finite set $\cY=\{y_{1},\dots,y_{k}\}$ of size $k$, for each 
 $0\le a<1$, set
$$
\cD(\cY;a)
:=
{
\#\{j\le k: y_{j}(\bmod~ 1)<a\}
\over
k
}
-a
.
$$
The {\em discrepancy} $\cD(\cY)$  is defined by
$$
\cD(\cY)
:=
\sup_{0\le a<1}\cD(\cY;a) 
- 
\inf_{0\le a<1}\cD(\cY;a) 
.
$$
}
\end{defi}

One always has 
$\cD(\cY) \le 1$.
The smallest possible discrepancy of a finite set $\cY$ is $\cD(\cY)=1/k$, attained by equally spaced elements $y_{j}=\frac{j}{k}, ~1 \le j \le k$. \\

A small discrepancy indicates that the set $\cY$ is close to equidistributed modulo $1$.
In particular, for an infinite sequence $\cX = \{ x_j: j \ge 1\}$, if $\cX_k = \{ x_j: 1 \le j \le k\}$ 
 then $\cX$ is uniformly distributed $(\bmod~1)$
if and only if the discrepancies $\cD( \cX_k) \to 0$ as $k \to \infty$.

%
%
%

\subsection{Benford's Law for  $3x+1$ Function Iterates}

Kontorovich and Miller \cite{KM05} considered iterates of the accelerated $3x+1$ function $U(n)$. 
Fix an odd integer $n=n_{0}$, and let $\{n_{1},n_{2},\dots\}$ be the sequence of iterates from the starting seed $n_{0} \in \Pi$, where $\Pi$ consists of all positive integers relatively prime to $6$.
 The main $3x+1$ conjecture asserts that this sequence is eventually periodic, and hence it is impossible for \eqn{1001} to hold!

 The following was their interpretation of (weak)  ``Benford behavior'' for the $3x+1$ function:

%
%
%
\begin{theorem}\label{th101}
For $x_{0}=n \in \Pi$,  denote its accelerated
$3x+1$ iterates by $x_{\kk}:=U^{(\kk)}(x_{0})$. 
Now set $y_{\kk}:=\log_{10} x_{\kk}$ and define the shifted variables 
$$
\gw_{k}:=
{
y_{k}
-
y_{0}
-k\log_{10} \left(\frac34\right)
} 
.
$$
Then, for  any $a\in[0,1)$,
 $$
\lim_{k\to\infty}
 \DD_{\Pi} \bigg[ x_0: ~
\gw_{k}
 (\bmod ~1)<a\bigg]
 =a.
 $$
\end{theorem}

 \paragraph{Proof.} This is established as Theorem 5.3 in Kontorovich and Miller \cite{KM05}.
 $~~~\bsq$\\

Arguably, the 
normalization from $y_{k}$ to $\gw_{k}$  in  Theorem \ref{th101} makes the above result only an approximation to ``true'' Benford behavior, which should be that 
$
\DD_{\Pi}[x_0: ~y_{k}~(\bmod ~1)<a] 
\to
a$
as $k\to\infty$.\\

 Lagarias and Soundararajan \cite{LagS06} were able to use 
 the non-accelerated $3x+1$ function $T$ to show another approximation to Benford behavior, as follows. \\

%
%
%
\begin{theorem}\label{th102}
{\em (Approximate Strong Benford's Law for $3x+1$ Map)}
Let $B>1$ be any integer base. Then
for a given $N \ge 1$ and each $X \ge 2^N$, most initial starting values $x_0$ in
$1 \le x_0 \le X$ have first $N$ initial $3x+1$ iterates $\{ x_k: 1 \le k \le N\}$ that satisfy
the discrepancy bound
\beql{841}
\cD\left( \{\log_B x_k (mod~ 1): 1 \le k \le N\}\right) \le 2 N^{-\frac{1}{36}}.
\eeq
The exceptional set $\sE(X, B)$ of initial seeds $x_0$ in $1 \le x_0 \le X$ that do not
satisfy the bound has cardinality
\beql{842}
|\sE(X, B)| \le c(B) N^{-\frac{1}{36}}
\eeq
where $c(B)$ is a positive constant depending only on the base $B$.
\end{theorem}

 \paragraph{Proof.} This is established as Theorem 2.1 in Lagarias and Soundararajan \cite{LagS06}.
 $~~~\bsq$\\

%
%
%

\subsection{Benford's Law for  $5x+1$ 
Function Iterates}

   The $5x+1$ map also exhibits similar ``
   Benford" behavior
 for its 
  iterates.
The results of \cite{KM05} apply to general $(d,g,h)$-Maps, in particular, to the $5x+1$ function,
giving a direct analogue of Theorem \ref{th101}.\\

The method of proof in \cite{LagS06} 
of Theorem~\ref{th102} should 
also extend to give qualitatively similar results in  the $5x+1$ case. This proof
relied on the Parity Sequence Theorem for the $3x+1$ map which
has an exact analogue for the $5x+1$ map.
The proof  in \cite{LagS06} also used some Diophantine 
approximation results for the transcendental number $\alpha_3:=\log_2 3$, and qualitatively
similar Diophantine approximation results are valid for $\alpha_5:= \log_2 5$
needed in the $5x+1$ case.\\

These rigorous results concern only the initial iterates of $5x+1$ trajectories. However since
the $5x+1$ map  conjecturally has divergent orbits, 
it seems a plausible guess 
that a strong form of Benford behavior might  hold on all infinite divergent orbits of
the $5x +1$ map.

%
%
%

 \section{2-Adic Extensions of $3x+1$ and $5x+1$ Maps}\label{sec10}

 What happens if we put these probabilistic models in a more general context?
 We can obtain a perfect set of symbolic dynamics  if we extend the domain of these
 maps to the $2$-adic integers. Such extensions are possible for both
 the $3x+1$ map $T_3(x)$ and the $5x+1$ map $T_5(x)$. 
 
%
%

\begin{theorem}~\label{th91}
The $3x+1$ map $T_3$ and the $5x+1$ map $T_5$ extend continuously from  maps on
the integers to maps on 
the $2$-adic integers $\ZZ_2$, viewing $\ZZ$ as a dense subset of $\ZZ_2$. 
Denoting the extensions by $\tilde{T}_3$ and $\tilde{T}_5$, respectively, these
maps have the following properties.

(i) Both maps $\tilde{T}_3$ and $\tilde{T}_5$are homeomorphisms of $\ZZ_2$ to itself.  

(ii) Both maps $\tilde{T}_3$ and $\tilde{T}_5$ are measure-preserving maps on $\ZZ_2$
for the standard  $2$-adic  measure $\mu_2$ on $\ZZ_2$.

(iii) Both maps $\tilde{T}_3$ and $\tilde{T}_5$ are strongly mixing with respect
to the measure $\mu_2$, hence ergodic.
\end{theorem}

\paragraph{Proof.}
For the $3x+1$ map, properties (i)-(iii) are stated in Lagarias \cite[Theorem K]{Lag85}.
The property of strong mixing is an ergodic-theoretic notion explained there. 
Akin \cite{Ak04} gives another proof of these facts for the $3x+1$ map.\\

For the $5x+1$ map, properties (i)-(iii) may be established by proofs similar to
the $3x+1$ map case. This is based on the fact that 
 an analogue of Theorem~\ref{th21} holds for the
symbolic dynamics of iterating the $5x+1$ map. 
It is also a corollary of results of Bernstein and Lagarias \cite[Sect. 4]{BL96},
whose results imply
that
 (i)-(iii) hold  more generally  for all $ax+b$-maps.
Here the {\em $ax+b$ map} $T_{a,b}$ is 
$$
T_{a,b}(x) :=
 \left\{
\begin{array}{cl}
\df{ax+b}{2}  & \mbox{if} ~~x \equiv 1~~ (\bmod ~ 2), \\
~~~ \\
\df{x}{2} & \mbox{if}~~ x \equiv 0~~ (\bmod ~2 ),
\end{array}
\right.
$$
where $a$ and $b$ are odd integers.
$~~~\bsq$\\

A much stronger ergodicity result is valid for  the $2$-adic extensions of
these maps. Define the
 {\em $2$-adic shift map} $S: \ZZ_2 \to \ZZ_2$ to be the $2$-to-$1$ map given 
for $\alpha = \sum_{j=0}^{\infty} a_j 2^j= .a_0 a_1 a_2... $, with each $a_j= 0$ or $1$, by
$$
S (\alpha) = S( .a_0a_1a_2 \cdots) := . a_1a_2 a_3 \cdots
$$
That is,
\beql{821}
S(\alpha) = 
 \left\{
\begin{array}{cl}
\df{\alpha-1}{2}  & \mbox{if} ~~\alpha \equiv 1~~ (\bmod ~ 2) \\
~~~ \\
\df{\alpha}{2} & \mbox{if}~~ \alpha \equiv 0~~ (\bmod ~2 )
.
\end{array}
\right.
\eeq
This map has the $2$-adic measure as Haar measure, and is mixing in
the strongest sense. 

%
%

\begin{theorem}~\label{th92} 
 The $2$-adic  extensions $\tilde{T}_3$ of the $3x+1$ map
 and $\tilde{T}_5$ of the $5x+1$ map 
 are each  topologically conjugate to the $2$-adic shift map,
by a conjugacy map $\Phi_3$, resp. $\Phi_5$. That is, these
maps are homeomorphisms of $\ZZ_2$ with
$\Phi_3^{-1} \circ  \tilde{T}_3\circ \Phi_3 = S$ and
$\Phi_5^{-1} \circ \tilde(T)_5 \circ \Phi_5 =S$.  

(1) The maps $\Phi_j$, $j=3$ or $5$,  
are solenoidal, i.e. for each $n \ge 1$ they have the property
$$
x \equiv y ~(\bmod~2^n) \longrightarrow \Phi_j(x) \equiv \Phi_j(y)~(\bmod~2^n).
$$

(2) The inverses of these conjugacy maps are explicitly given by
$$
\Phi_j^{-1}(\alpha) := \sum_{k=0}^{\infty} \left(T_j^{(k)}(\alpha) ~(\bmod~2) \right)2^k,
$$
for $j= 3$  or $5$, and the residue $(\bmod~2)$ is taken to be $0$ or $1$.
 \end{theorem}
 
 \paragraph{Proof.} These results follow from  Bernstein and Lagarias \cite[Sect. 3, 4]{BL96},
 where results are proved for a general class of mappings including both the $3x+1$ map
 and $5x+1$ map.
 $~~~\bsq$\\
 
 Theorem~\ref{th92} immediately gives the following corollary.\\
%
%

\begin{corollary}~\label{cor92a} 
 The $2$-adic  extensions $\tilde{T}_3$ of the $3x+1$ map
 and $\tilde{T}_5$ of the $5x+1$ map 
 are  topologically conjugate and metrically conjugate maps.
 \end{corollary}

 The corollary  
 shows that from the viewpoint of extensions to the $2$-adic integers, the 
 $3x+1$ maps and the $5x+1$ maps have identical ergodic theory properties, i.e. they
 are both conjugate to the shift map. That is, their symbolic dynamics is ``the same"
 in the topological sense, and their dynamics is also identical 
 in the measure-theoretic sense. \\
 
 The original $3x+1$ problem (resp. $5x+1$ problem) concerns their behavior
 when restricted to the  dense set $\ZZ$ inside $\ZZ_2$. This set $\ZZ$ is countable, so has
 $2$-adic measure zero, so the  general properties of ergodic theory allow no
 conclusion to be drawn about behavior of iteration  on these maps on $\ZZ$.
 Indeed 
 empirical data and the stochastic models above show that the dynamics of iteration
 of the $3x+1$ map and $5x+1$ map are ``not the same"  on $\ZZ$. \\
 
 To conclude, we remark 
 that the  two accelerated functions $U_3$ and $U_5$ also make sense $2$-adically, 
 in a restricted domain.
Let $\ZZ_2^{\times} = \{ \alpha \in \ZZ_2:~\alpha \equiv 1~(\bmod~2)\}$. 
We have 
 $U_3: \ZZ_2^{\times} \to \ZZ_2^{\times} \cup \{0\}$
(in the latter case we set $U( -\frac{1}{3})=0$.)
and $U_5: \ZZ_2^{\times} \to \ZZ_2^{\times} \cup \{0\}$
(in the latter case we set $U( -\frac{1}{5})=0$.) 
It might prove worthwhile to find invariant measures
for these functions, and to study their  ergodic-theoretic behavior.

%
%
%

 \section{Concluding Remarks}\label{sec11}
 
 We have presented results on stochastic models simulating
 aspects of the behavior of the $3x+1$ function and $5x+1$ problems.
 These models resulted in  specific predictions about various
 statistics of the orbits of these
 functions under iteration, which can be tested empirically. The
 experimental tests done so far have  generally been consistent with these predictions.\\
 
%
%
%

 \subsection{Comparisons}
 
 We compare and contrast the behavior of these two maps under iteration.
The $3x+1$ map  and $5x+1$ map are similar in the following dimensions.

\begin{enumerate}
\item
({\em Symbolic dynamics})
The allowed symbolic dynamics of even and odd iterates is the same for the
$3x+1$ and $5x+1$ maps. Every  finite symbol sequence is legal.

\item
({\em Periodic orbits on the integers})
Conjecturally, both the $3x+1$ map and $5x+1$ maps have a finite number of
distinct periodic orbits on the domain $\ZZ$.

\item
({\em Periodic orbits on rational numbers with odd denominator})
Every possible symbolic dynamics for 
a periodic orbit is the periodic orbit for some rational starting point, for both the
$3x+1$ map and $5x+1$ map.
That is,  extensions of the maps $T_3$ and $T_5$ to rational numbers with odd denominator
each have $2^p$ periodic orbits of period $p$, for each $p \ge 1$.
Here the period $p$ may not be the minimal period of the orbit, 
so a period $k$ orbit is also counted as a period
$p=kn$ orbit for each $k \ge 1$.

\item
({\em Benford Law behavior})
Both the initial  $3x+1$ function iterates of a random starting  point, and the initial $5x+1$
iterates of a random starting point, with high probability exhibit strong Benford law
behavior to any integer base $B \ge 2$. 

\item
({\em $2$-adic extensions})
The $2$-adic extensions of the two maps are topologically and metrically conjugate.
Therefore they have the same dynamics in the topological sense, and in the
ergodic theory sense, on the domain $\ZZ_2$.
\end{enumerate}
 
The main differences between the $3x+1$ maps and $5x+1$ maps concerns the 
change in size of  their interates.

\begin{enumerate}
\item
({\em Short-term behavior of iterates})
 For the $3x+1$ map, the initial steps of most orbits shrink in size, while for the
$5x+1$ map most orbits expand in size.
This is rigorously   quantified in 
\S\ref{sec2} and \S\ref{sec7}.
\item
({\em Long-term behavior of iterates})
The $3x+1$ and $5x+1$ conjecturally differ greatly in their long-term behavior of
orbits on the integers. For the $3x+1$ map, conjecturally all orbits are bounded.
For the $5x+1$ map, conjecturally a density one set of integers have unbounded orbits.
\end{enumerate}

It is the long term behavior of iterates where all the difficulties connected with the
$3x+1$ and $5x+1$ function lie.

%
%
%

 \subsection{Insights}

  Comparison of the results of these stochastic models, combined
  with  deterministic  results, 
deliver certain insights in
understanding the $3x+1$ and $5x+1$ problem, and suggest topics for 
further work. \\

First, the $2$-adic results indicate that the differences in  of the dynamics
of the $3x+1$ map $T_3$ and $5x+1$ map on the integers 
are invisible at the level of measure theory. Therefore these differences  must 
depend in some way  on number-theoretic features inside the integers $\ZZ$.\\

 Second,  the behavior of the iteration of these function of 
 in $\ZZ$, viewed 
 inside the $2$-adic framework,
  must be encoded  in the specific 
 properties of the  conjugacy maps $\Phi_3$ and $\Phi_5$ identifying these
 maps with the $2$-adic shift map. Here
 we note that there is an explicit formula for the $3x+1$ conjugacy map, obtained by
 Bernstein \cite{Be94}, and there is an analogous formula for the
 $5x+1$ conjugacy map as well. These conjugacy
 maps have an intricate structure, detailed in \cite{BL96}, which might be worthy
 of further investigation. \\
 
 Third, we observe that  the ergodic behavior of the $2$-adic extensions 
 is exactly the behavior that served as a framework to formulate the random walk
 models presented in 
 \S\ref{sec3}, \S\ref{sec5}, 
 and
 \S\ref{sec7}. These random walk models yield
 information by combining these model iterations
 with estimates of  the size of iterates in the standard absolute
 value on the real line $\RR$. That is, they use information from an
 archimedean norm, rather than the non-archimedean norm on the
 $2$-adic integers. Perhaps one needs to consider models that incorporate
 both norms at once, e.g. functions on $\RR \times \ZZ_2$.\\

 Fourth,  a suitable maximal domain, larger than $\ZZ$, on which to understand the
 difference between the $3x+1$ map $T_3$ dynamics  and the $5x+1$ map
 $T_5$ dynamics appears to  be the domain 
 $$
 \QQ_{(2)} := \QQ \cap  \ZZ_2,
 $$
 i.e. the set of rational numbers that are $2$-adic integers.
 The set  $\QQ_{(2)}$ is  exactly the set of rational numbers having an odd 
 denominator, and both $T_3$ and $T_5$ leave the set $\QQ_{(2)}$
 invariant. This set includes all periodic orbits of both $T_3$ and $T_5$,
 and from the viewpoint of existence of periodic orbits, these two maps
 are the same on $\QQ_{(2)}.$
 The difference in the dynamics of these maps on $\ZZ$ seems to
 have something to do with the distribution of these periodic orbits.
 Viewing $\QQ_{(2)}$ as having the topology induced from the
 $2$-adic topology, one may conjecture that $T_3$ and $T_5$ are
 {\em not}  topologically conjugate mappings on this domain.\\
 
 Fifth, the $5x+1$ map exhibits various ``exceptional" behaviors. Although
 almost all of its integer orbits (conjecturally) diverge, nevertheless there 
 exists  an infinite exceptional
 set of integers that have eventually periodic orbits. The density (fractional dimension)
 of such integers is predicted (conjecturally) to be a constant $\delta_{5} \approx 0.649$,
 solving a large
 deviations functional equation. This seems a hard problem to
 resolve rigorously.  
 Now, for the $3x+1$ map, a similar
 prediction is 
 made by the models for the  growth constant $g=1$. It too
 is the solution of a large deviations functional equation. We currently
 know that $1 \ge g \ge 0.84$.
 This analogy suggests that rigorously
 proving that the growth constant $\delta_3 =1$ may turn out to be a much harder problem than
 it seems at first glance.\\
 
    Sixth, we note that there are extensions of the maps for backwards iteration
    to larger domains, to the invertible $3$-adic integers  $\ZZ_3^{\ast}$ for the
    $3x+1$ map, and to the invertible $5$-adic integers $\ZZ_5^{\ast}$ for the
    $5c+1$ map.  In effect the  branching random walk models may fruitfully be extended to allowing
    root  node labels that are invertible $3$-adic integers (resp. 5-adic integers), and this provides
  enough information to grow the entire infinite tree. Various interesting properties of
 the extended $3x+1$ trees obtained this way have been obtained,
  cf. \cite{AL95c}. This is a topic worth further investigation. 

%

\end{document}